\documentclass{siamltex1213}
\usepackage[english]{babel}
\usepackage{amsmath}
\usepackage{amstext}
\usepackage{amssymb,amsfonts}
\usepackage{mathtools}
\usepackage{verbatim}
\usepackage{xcolor}

\usepackage{epstopdf}
\usepackage{subfigure}
\usepackage{dsfont}
\DeclareGraphicsExtensions{.eps}
\usepackage{tikz}
\usetikzlibrary{arrows,positioning}

\usepackage{pgfplots}
\pgfplotsset{compat=1.16}
\usepgfplotslibrary{groupplots}
\usepackage{cancel}
\allowdisplaybreaks

\allowdisplaybreaks
\def\Ebox#1#2{%
\medskip
\begin{center}
  \strut\epsfxsize=#1 \hsize\epsfbox{#2}
\end{center}
\smallbreak}

\usepackage{epsfig}

\def\IR{I \kern-0.35em R\,}

\newcommand{\sy}[1]{{\color{black} #1}}
\newcommand{\syr}[1]{{\color{black} #1}}

\def\ess{{\mathrm{ess}}}
\def\Proj{\mathop{\rm Proj}}

\def\G{\mathcal{G}}

\def\l1{L_1}
\def\l2{L_2}

\newtheorem{example}{Example}
\newtheorem{assumption}{Assumption}
\newtheorem{remark}{Remark}

\newenvironment{proofsketch}{\noindent{\it Sketch of Proof:}\hspace*{1em}}{\qed\bigskip}
\newcommand{\QEDA}{\hfill\ensuremath{\square}}

\allowdisplaybreaks

\UseRawInputEncoding

\begin{document}
\sloppy
\title{Geometry of Information Structures, Strategic Measures and associated Control Topologies} 
\author{Naci Saldi and Serdar Y\"uksel
\thanks{Naci Saldi is with the Department of Natural and Mathematical Sciences, {\"{O}}zye\u{g}in University, \c{C}ekmek\"{o}y, Turkey, Email: \{naci.saldi@ozyegin.edu.tr\}. 
Serdar Y\"uksel is with the Department of Mathematics and
    Statistics, Queen's University, Kingston, Ontario, Canada; email: yuksel@mast.queensu.ca. This research was partially supported by the Natural Sciences and Engineering Research Council of Canada (NSERC). }}

\maketitle

\begin{abstract}
\sy{In many areas of applied mathematics, engineering, and social and natural sciences, decentralization of information is a key aspect determining how to approach a problem. In this review article, we study information structures in a probability theoretic and geometric context. We define information structures, place various topologies on them, and study closedness, compactness and convexity properties on the strategic measures induced by information structures and decentralized control/decision policies under varying degree of relaxations with regard to access to private or common randomness. Ultimately, we present existence and tight approximation results for optimal decision/control policies. We discuss various lower bounding techniques, through relaxations and convex programs ranging from classically realizable and classically non-realizable (such as quantum and non-signaling) relaxations. For each of these, we establish closedness and convexity properties and also a hierarchy of correlation structures. As a second main theme, we review and introduce various topologies on decision/control strategies defined independent of information structures, but for which information structures determine whether the topologies entail utility in arriving at existence, compactness, convexification or approximation results. These approaches, which we term as the strategic measures approach and the control topology approach, lead to complementary results on existence, approximations and upper and lower bounds in optimal decentralized decision and control.}
\end{abstract}


\section{Introduction}

In statistical decision theory, stochastic control theory, information theory, game theory, economics, quantum physics and computer science, information structures determine which unit,decision maker/agent/controller knows what information. Accordingly, information structures is a subject that has a very broad appeal and interpretation depending on the context. 

\begin{figure}[h]
\centering
\Ebox{.50}{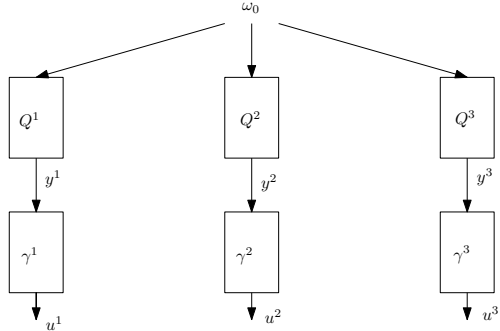}
\caption{ \sy{An example of an information structure. Here, $Q^i(y^i \in \cdot |\omega_0) := P(\eta^i(\omega) \in \cdot | \omega_0)$, $i=1,2,3$.}} \label{InformationStructureFlow}
\end{figure}

To help the reader gain an appreciation of what the rest of the paper entails, consider Figure \ref{InformationStructureFlow}.

Here, three decision makers have access to some local information (e.g. $y^2$) and each one make a action selections (e.g., $u^2$). This selection has to be based only on the available information (and perhaps with some additional independent randomization device). One question of interest in many of the aforementioned disciplines is to study the sets of all possible correlation structures with regard to the random variables $(\omega_0, y^1, y^2, y^3, u^1, u^2, u^3)$. Suppose further that these decision makers wish to minimize a cost criterion of the form $E[c(\omega_0,u^1,u^2,u^3)]$. What can be said about the optimal solutions, their existence, approximations, and numerical methods for exact solutions or upper or lower bounds? 

The operational goal in our paper is, for setups including that depicted in Figure \ref{InformationStructureFlow} but also under much more general {\it information structures}, to develop a probability theoretic and topological approach to information structures and decision/control policies towards arriving at: (i) convexity properties, (ii) continuity and compactness properties, (iii) existence results, (iv) convex relaxations to facilitate numerical and analytical methods, and (v) approximation results for optimal policies under decentralized information structures. 

Different scientific communities have their own notations, terminology, and machinery to study such problems. While in our analysis we will often have a bias towards a stochastic control theoretic angle and language, we will be broadly touching on ideas from various disciplines and often explain parallel contributions in the literature. Accordingly, the scope of our article entails results from a very broad literature across a variety of disciplines (ranging from game theory, quantum physics, information theory, probability, and control). Some of the results presented here have not been published elsewhere, some are based on the authors' prior work, but much of the material is a re-interpretation and unification of related results in the literature (both across various disciplines and across time, meaning that classical results are re-interpreted with more modern findings) whose connections have not been made explicit and reported, to our knowledge.  

\subsection{Summary of results and outline of the article}\label{SectionOutline}

\begin{figure}[h]
\centering
\Ebox{.50}{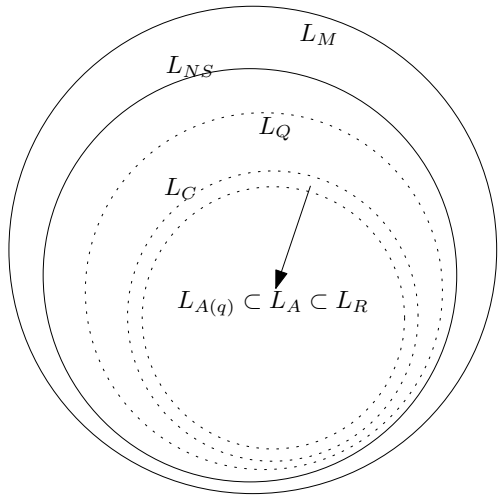}
\caption{ \sy{Hierarchies and a summary of the results reported in the paper on decentralized strategic measures and correlations: (i) $L_{A(q)}(\mu)$, quantized admissible strategic measures, is not closed, not convex, dense in $L_A(\mu)$. (ii) $L_A(\mu)$, set of strategic measures induced by admissible measurable functions, is not closed, not convex (under independent static reduction $L_A(\mu)$ is dense in $L_R(\mu)$ under mild conditions). (iii) $L_R(\mu)$, those generated with local private (relaxed policies) randomness, is typically not closed (under independent static reduction:  $L_R(\mu)$ is closed), not convex. (iv) $L_C(\mu)$, those with arbitrary common randomness, convex, typically not closed. (iv) $L_Q(\mu)$, quantum strategic measures, convex, not closed (but closed under dimension constraints). (v) $L_{NS}(\mu)$, non-signaling strategic measures, is closed and  convex. (vi) $L_M(\mu)$, those with local-Markov conditional independence, is closed and convex.}}
 \label{HierInfoStrategicMeasure}
\end{figure}

The article has two main themes: In the first theme, we will study topological and geometric properties of correlation structures induced on strategic measures under decentralized control/decision policies and information structures. We will establish the following hierarchies depicted, with a summary of the results reported, in Figure \ref{HierInfoStrategicMeasure}.


As our second main theme, we study a {\it control topologies} approach in Section~\ref{TopologyOnPolicies}, where we will review or introduce various topologies on decision/control strategies (defined independently from information structures), but for which information structures determine whether the topologies entail utility in arriving at existence, compactness, convexification or approximation results. These topologies include those generated by product metric of individual controllers viewed as subsets of appropriate probability measures, relaxed control policies and wide-sense admissible control policies, weak$^*$ topology on randomized policies, policies satisfying conditional independence properties leading to a universal dynamic program, and topologies that satisfy a closedness property under weak convergence of finite dimensional marginals. 

These two approaches, which we may term as {\it strategic measures approach} vs. {\it the control topology approach}, lead to complementary results on existence, approximations and upper and lower bounds and solution methods in optimal stochastic control.


In Section \ref{SectionInfoStructures} we introduce information structures in decentralized decision making. Section~\ref{strategic}, we introduce probability measures (i.e., strategic measures) induced on the product spaces of state, measurements and actions under admissible, and some relaxed, control policies. In Section~\ref{regularity}, we study compactness properties of the set of strategic measures and state existence results on optimal policies. In Section~\ref{convexRelax}, we study convex relaxations including non-signaling and quantum information theoretic ones, on strategic measures. In Section~\ref{finiteSupportApp}, we study finite measurement information structures and establish their near optimality. In Section~\ref{TopologyOnPolicies}, we study an alternative approach of placing topologies directly on individual policies. In Section~\ref{relaxedPOMDPs}, we review partially observed Markov decision processes and highlight that when a relaxation of control policies is allowed, this may lead to a strict improvement of optimal performance, thus invalidating the main purpose of relaxations. 

The paper ends with concluding remarks, future directions, and finally, a collection of open problems.

\section{Information Structures}\label{SectionInfoStructures} 

\subsection{Witsenhausen's characterization of information structures}\label{witsenInfoStructureReview}

Hans Witsenhausen's contributions \cite{wit75,wit88,WitsenStandard} to stochastic control theory and information theory, and his characterization of information structures in decentralized stochastic control have been crucial in our modern understanding of decentralized stochastic control and decision theory. In this section, we introduce the characterizations as laid out by Witsenhausen, termed as {\it the Intrinsic Model} \cite{wit75}; see \cite{YukselBasarBook} and \cite{CDCTutorial} for a more comprehensive overview and further characterizations and classifications of information structures. In this model (described in discrete time), any action applied at any given time is regarded as applied by an individual decision maker/agent, who acts only once. One advantage of this model, in addition to its generality, is that the characterizations regarding information structures can be concisely described.

Suppose that in the decentralized system considered below, there is a pre-defined order in which the decision makers act. Such systems are called {\it sequential teams} (for non-sequential teams, we refer the reader to Andersland and Teneketzis \cite{AnderslandTeneketzisI}, \cite{AnderslandTeneketzisII} and Teneketzis \cite{Teneketzis2}, in addition to Witsenhausen \cite{WitsenhausenSIAM71} and \cite[p. 113]{YukselBasarBook}). Suppose that in the following, the action and measurement spaces are standard Borel spaces; that is, Borel subsets of Polish (complete, separable, and metric) spaces. In the context of a sequential system, the {\it Intrinsic Model} has the following components:

\begin{itemize}
\item A collection of {\it measurable spaces} 
$$\left\{(\Omega, {\cal F}),
(\mathbb{U}^i,{\cal U}^i), (\mathbb{Y}^i,{\cal Y}^i), i \in {\cal N}\right\},$$ 
with ${\cal N}:=\{1,2,\cdots,N\}$, specifying the system's distinguishable events, and the control and measurement spaces of decision makers (DMs). Here $N=|{\cal N}|$ is the number of control actions taken, and each of these actions is taken by an individual (different) DM (hence, even a DM with perfect recall can be
regarded as a separate decision maker every time it acts). The pair $(\Omega, {\cal F})$ is a
measurable space (on which an underlying probability may be defined). The pair $(\mathbb{U}^i, {\cal U}^i)$
denotes the Borel space with its Borel $\sigma$-algebra from which the action $u^i$ of DM~$i$ is selected. The pair $(\mathbb{Y}^i,{\cal Y}^i)$ denotes the observation/measurement space with its Borel $\sigma$-algebra for DM~$i$.
\item A {\it measurement constraint} which establishes the connection between the observation variables and the system's distinguishable events. The $\mathbb{Y}^i$-valued observation variables are given by 
$$y^i=\eta^i(\omega,{\bf u}^{[1,i-1]}),$$ 
where ${\bf u}^{[1,i-1]}=\{u^k, k \leq i-1\}$ and $\eta^i$ is a measurable function. Hence, the information variable $y^i$ induces a $\sigma$-field, denoted by $\sigma(y^i)$, over $\Omega \times \prod_{k=1}^{i-1} \mathbb{U}^k$
\item A {\it design constraint} which restricts the set of admissible $N$-tuple control laws 
$$\underline{\gamma}= \{\gamma^1, \gamma^2, \dots, \gamma^N\},$$ 
also called {\it designs} or {\it policies}, to the set of all measurable control functions, so that $u^i = \gamma^i(y^i)$, with $y^i=\eta^i(\omega,{\bf u}^{[1,i-1]})$ and $\gamma^i$ is a measurable function. Let $\Gamma^i$ denote the set of all admissible policies for DM~$i$ and let ${\bf \Gamma} = \prod_{i=1}^N \Gamma^k$.
\end{itemize}
We note that, the intrinsic model of Witsenhausen gives a set-theoretic characterization of information fields, however, for standard Borel spaces, the model above is equivalent to that of Witsenhausen's. Additionally, we can also introduce a fourth component:
\begin{itemize}
\item A {\it probability measure} $P$ defined on $(\Omega, {\cal F})$ which describes the uncertainty on the random events in the model. 
\end{itemize}

Under this intrinsic model, an Information structure (IS) is {\it dynamic} if the information available to at least one DM is affected by the action of at least one other DM. An IS is {\it static}, if the information available at every decision maker is only affected by exogenous disturbances (i.e., state of the Nature) $\omega \in \Omega$; that is no other decision maker can affect the information at any given decision maker. Figure \ref{InformationStructureFlow} is a depiction for a static team problem.


ISs can also be classified as {\it classical}, {\it quasi-classical} or {\it nonclassical}. An IS is {\it classical} if $y^i$ contains all of the information available to DM~$k$ for $k < i$; that is, information is expanding (also known as the {\it perfect-recall} property). An IS is {\it quasi-classical} or {\it partially nested}, if whenever $u^k$, for some $k < i$, affects $y^i$ through the measurement function $\eta^i$, $y^i$ contains $y^k$ (that is $\sigma(y^k) \subset \sigma(y^i)$). An IS which is not partially nested is {\it nonclassical}.

For any $N$-tuple of policies $\underline{\gamma} = \{\gamma^1, \cdots, \gamma^N\}$ let a cost function be defined as:
\begin{eqnarray}\label{lossF}
J(\underline{\gamma}) = E^{\underline{\gamma}}\left[c(\omega_0,{\bf u})\right] = E\left[c(\omega_0,\gamma^1(y^1),\cdots,\gamma^N(y^N))\right],
\end{eqnarray}
for some non-negative measurable loss (or cost) function $c: \Omega_0 \times \prod_{i=1}^N \mathbb{U}^i \to \mathbb{R}_+$. Here, we have the notation ${\bf u}=\{u^i, i \in {\cal N}\}$, and $\omega_0$ may be viewed as the cost function relevant exogenous variable contained in $\omega$. 

\begin{definition}\label{Def:TB1}\index{Optimal team cost}
For a given stochastic team problem with a given information
structure, an $N$-tuple of policies
${\underline \gamma}^*:=({\gamma^1}^*,\ldots, {\gamma^N}^*)\in {\bf \Gamma}$ is
an {\it optimal team decision rule} ({\it team-optimal
decision rule} or simply {\it team-optimal solution}) if
\begin{equation}J({\underline \gamma}^*)=\inf_{{{\underline \gamma}}\in {{\bf \Gamma}}}
J({{\underline \gamma}})=:J^*. \label{eq:5}
\end{equation} 
The expected cost achieved by this strategy $J^*$ is the optimal team cost.
\end{definition}

In the following, we will denote by bold letters the ensemble of random variables across the DMs; that is ${\bf y}=\{y^i, i=1,\cdots,N\}$ and ${\bf u}=\{u^i, i=1,\cdots,N\}$.

\subsection{Independent-measurements reduction of sequential teams}\label{EquivIS}
 
Following Witsenhausen \cite[Eqn (4.2)]{wit88}, as reviewed in \cite[Section 3.7]{YukselBasarBook}, we say that two information structures are equivalent if: (i) The policy spaces are equivalent/isomorphic in the sense that policies under one information structure are realizable under the other information structure, (ii) the costs achieved under equivalent policies are identical, and (iii) if there are constraints in the admissible policies, the isomorphism among the policy spaces preserves the constraint conditions. 

A large class of sequential dynamic team problems admit an equivalent information structure which is static. This is called the {\it static reduction} of a dynamic team problem. 

For some, but not all, results to be presented in our paper, we need to go beyond a static reduction, and we will need to make the measurements independent of each other as well as $\omega_0$. This is not possible for every team which admits a static reduction, for example quasi-classical team problems with LQG models \cite{HoChu} do not generally admit such a further reduction, since the measurements are partially nested: for partially nested (or quasi-classical) information structures, static reduction has been studied by Ho and Chu in the specific context of LQG systems \cite{HoChu} and for a class of non-linear systems satisfying restrictive invertibility properties \cite{ho1973equivalence}. 

Witsenhausen refers to such an information structure as {\it independent static} in \cite[Section 4.2(e)]{wit88}. One can also reduce a static team problem into an independent static form.

Note that the Intrinsic Model is equivalent to the following model \cite{wit88}. The probability space $(\Omega,{\cal F},P)$ is the product of $N+1$ probability spaces $(\Omega_i,{\cal F}_i,P_i)$, $i=0,\ldots,N$; that is, $\Omega  = \prod_{i=0}^N \Omega_i$ and $P(d\omega) = \prod_{i=0}^N P_i(d\omega^i)$. Each DM~$i$ measures $y^i=\eta^i(\omega_0,\omega_i,u^1,\ldots,u^{i-1})$ and the decisions are generated by $u^i=\gamma^i(y^i)$, with $1 \leq i \leq N$. Here $\omega =(\omega_0,\omega_1,\cdots,\omega_N)$ are primitive (exogenous) variables. Under this equivalent model, we can view $\eta^i(\omega_0,\omega_i,u^1,\ldots,u^{i-1})$ as a measurement channel with input $(\omega_0,u^1,\ldots,u^{i-1})$ and output $y^i$, where $\omega_i$ is the noise; that is 
 \[g^i(dy^i | \omega_0, u^1,u^2,\cdots,u^{i-1})\] 
is a (controlled) stochastic kernel (to be defined later). Equivalently, through standard stochastic realization results (see \cite[Lemma 1.2]{gihman2012controlled} or \cite[Lemma 3.1]{BorkarRealization}), we can represent any (controlled) stochastic kernel $g^i(dy^i | \omega_0, u^1,u^2,\cdots,u^{i-1})$ in a functional form $y^i=\eta^i(\omega_0,\omega_i,u^1,u^2,\cdots,u^{i-1})$ for some independent $\omega_i$ and measurable $\eta^i$.

This team model admits an {\it independent static} reduction provided that the following absolute continuity condition holds: 

\smallskip

({\bf AC}): For every $i \in {\cal N}$, there exists a reference probability measure $Q^i$ on $\mathbb{Y}^i$ and a measurable function $f_i$ such that for all Borel $S \subset \mathbb{Y}^i$:
\begin{eqnarray}
&& g^i(y^i \in S | \omega_0,u^1,u^2,\cdots,u^{i-1}) \nonumber \\
&& \quad \quad \quad \quad = \int_{S} f_i(y^i,\omega_0,u^1,u^2,\cdots,u^{i-1}) \, Q^i(dy^i). \nonumber
\end{eqnarray}

Under this absolute continuity condition, since the action of each DM is determined by the measurement variables under a policy, we can write
\begin{eqnarray}
P(d\omega_0,d{\bf y}, d{\bf u}) = P(d\omega_0) \prod_{i=1}^N \bigg(f_i\left(y^i,\omega_0,u^1,u^2,\cdots,u^{i-1}\right) \, Q^i(dy^i) \, 1_{\{\gamma^i(y^i) \in du^i\}}\bigg). \nonumber
\end{eqnarray}
The cost function $J(\underline{\gamma})$ can then be written as
\[J(\underline{\gamma})= \int P(d\omega_0) \prod_{i=1}^N \bigg(f_i\left(y^i,\omega_0,u^1,u^2,\cdots,u^{i-1}\right) \, Q^i(dy^i)\bigg) \, c(\omega_0,{\bf u}),\]
where $u^i = \gamma^i(y^i)$ for $1 \leq i \leq N$. Now, the measurement variables can be regarded as independent from each other and also from $\omega_0$, and by incorporating the density functions $\{f_i\}$ into $c$, we can obtain an equivalent {\it independent static team} problem. Hence, the essential step is to appropriately adjust the probability space and the cost function. The new cost function may now explicitly depend on the measurement values; that is, 
\begin{eqnarray}
c_s(\omega_0,{\bf y}, {\bf u}) = c(\omega_0,{\bf u}) \, \prod_{i=1}^N f_i(y^i,\omega_0,u^1,u^2,\cdots,u^{i-1}). \label{c_sDefn}
\end{eqnarray}

\noindent Here we can reformulate even a static team to one which is, clearly still static, but now with independent measurements which are also independent from the cost relevant exogenous variable $\omega_0$. Such a condition is in general not restrictive. Indeed, as Witsenhausen notes, a static reduction always holds when the measurement variables take values from countable set since a reference measure as in $Q^i$ above can be always constructed on the measurement space $\mathbb{Y}^i$ (e.g., $Q^i(z) = \sum_{j \geq 1} 2^{-j} 1_{\{z = m_j\}}$ where $\mathbb{Y}^i=\{m_j, j \in \mathbb{N}\}$) so that the absolute continuity condition always holds. We refer the reader to \cite{charalambous2016decentralized} for relations with classical continuous-time stochastic control where the relation with Girsanov's classical measure transformation \cite{girsanov1960transforming}\cite{benevs1971existence} is recognized, and \cite[p. 114]{YukselBasarBook} for further discussions. For discrete-time partially observed stochastic control, similar arguments had been presented in \cite{Bor00}, \cite{Bor07}.

%
%

\section{Decentralized Strategic Measures}\label{strategic}

For classical stochastic control problems, strategic measures were defined (see \cite{schal1975dynamic}, \cite{piunovskii1998controlled}, \cite{dynkin1979controlled} and \cite{feinberg1996measurability}) as the set of probability measures induced on the product (sequence) spaces of the states, measurements, and actions; that is, given an initial state distribution and a policy, one can uniquely define a probability measure on the product space of the states, measurements, and actions. Certain measurability, compactness, and convexity properties of strategic measures for classical stochastic control problems were studied in \cite{dynkin1979controlled,piunovskii1998controlled,feinberg1996measurability,blackwell1976stochastic}. 

In \cite{YukselSaldiSICON17}, strategic measures for decentralized stochastic control problems were introduced and many of their properties were established. For decentralized stochastic control problems, considering the set of strategic measures along with compactification and/or convexification of these sets of measures through introducing private and/or common randomness allow one to place operationally flexible topologies (such as those leading to a standard Borel space, e.g., weak convergence topology, among others) on the set of strategic measures, as we will study in the following. 

\subsection{Measurable policies as a subset of randomized policies and strategic measures}\label{distributionalPolicy}

A common method in control theory is to view a measurable policy as a special case of {\it relaxed} policies where relaxation is often employed by randomization. Such an approach has been ubiquitously adopted in various fields often with different terminology (e.g., relaxed controls (Young topology) in optimal deterministic control \cite{mascolo1989relaxation} \cite{young1937generalized}, distributional strategies in economics \cite{milgrom1985distributional} \cite{mertens2015repeated}, local hidden variables in quantum information theory \cite{BaLiMaPiPoRo05,ClHoToWa04,Wat18}, optimal quantization \cite{YukselOptimizationofChannels} etc.). 

We recall here the following representation result \cite{BorkarRealization} (see also Section \ref{YoungM}). Let $\mathbb{X}, \mathbb{M}$ be Borel spaces. Let the notation ${\cal P}(\mathbb{X})$ denote the set of probability measures on $\mathbb{X}$. Consider the set of probability measures
\begin{eqnarray}\label{extremePointQuan0}
\Theta: = \left\{\zeta \in {\cal P}(\mathbb{X} \times \mathbb{M}): \zeta(dx,dm) = P(dx) \, Q^f(dm|x), Q^f(\cdot | x) = 1_{\{f(x) \in \cdot\}}, f : \mathbb{X} \to \mathbb{M} \right\}, \nonumber
\end{eqnarray}
on $\mathbb{X} \times \mathbb{M}$ having fixed input marginal $P$ on $\mathbb{X}$ and the stochastic kernel from $\mathbb{X}$ to $\mathbb{M}$ is realized by some measurable function $f : \mathbb{X} \to \mathbb{M}$. We equip this set with weak convergence topology. This set is the (Borel measurable) set of the extreme points of the set of probability measures on $\mathbb{X} \times \mathbb{M}$ with a fixed marginal $P$ on $\mathbb{X}$. For compact $\mathbb{M}$, the Borel measurability of $\Theta$ follows from \cite{Choquet} since the set of probability measures on $\mathbb{X} \times \mathbb{M}$ with a fixed marginal $P$ on $\mathbb{X}$ is a convex and compact set in a complete separable metric space, and therefore, the set of its extreme points is Borel measurable.  Moreover, the non-compact case holds by \cite[Lemma 2.3]{BorkarRealization}.
Furthermore, given a fixed marginal $P$ on $\mathbb{X}$, any stochastic kernel $Q$ from $\mathbb{X}$ to $\mathbb{M}$ can be identified by a probability measure $\xi \in {\cal P}(\Theta)$ such that
\begin{eqnarray}\label{convR0}
Q(\cdot|x)  = \int_{\Theta} \xi(dQ^f) \, Q^f(\cdot|x). 
\end{eqnarray}
In particular, a stochastic kernel can thus be viewed as an integral representation over probability measures induced by deterministic policies.

For a team setup, for any DM~$k$, let
\begin{align} 
&\Theta^k \coloneqq \bigg\{\zeta \in {\cal P}(\mathbb{Y}^k \times \mathbb{U}^k): \zeta = P_k \,Q^{\gamma^k}, \nonumber \\
&\phantom{xxxxxxxxxxxxxxxxxxxxxxx} Q^{\gamma^k}(\cdot|y^k) = 1_{\{\gamma^k(y^k) \in \cdot\}}, \gamma^k \in \Gamma^k, P_k(\cdot) = P(y^k \in \cdot) \bigg\}. \nonumber 
\end{align}
For a static team, $P_k$ would be fixed; that is, independent of the policies of the preceding DMs. Therefore, in static case, in view of (\ref{convR0}), any element $\zeta \in {\cal P}(\mathbb{Y}^k \times \mathbb{U}^k)$ with fixed marginal $P_k$ on $\mathbb{Y}^k$ can be expressed as the mixture of $\Theta^k$
\begin{eqnarray}\label{convR}
\zeta (A)  = \int_{\Theta^k} \xi^k(dQ) \, Q(A), \quad  A \in {\cal B}(\mathbb{Y}^k \times \mathbb{U}^k),
\end{eqnarray}
for some $\xi \in {\cal P}(\Theta^k)$. In the sequel, we generalize this idea to the set of strategic measures induced by measurable policies and define various relaxed policies that are obtained as a mixture of measurable policies. Indeed, instead of viewing $N$-tuple of policies as the joint strategy of DMs, we regard the induced probability distribution on the product space of state, measurements, and actions as the joint strategy and name it \textit{strategic measure}. However, we will recognize that the $N$-tuple view above is also of both operational and mathematical interest as will be detailed in Section \ref{TopologyOnPolicies}.

\subsection{Sets of strategic measures for static teams}\label{setsOfStrategicMeasures}

Consider a static team problem defined under Witsenhausen's intrinsic model.  In the following, $B = B^0 \times \prod_{k=1}^N (A^k \times B^k)$ are used to denote the cylindrical Borel sets in $\Omega_0 \times \prod_{k=1}^N (\mathbb{Y}^k \times \mathbb{U}^k)$. 

\smallskip

Let $L_A(\mu)$ be the set of strategic measures induced by all admissible measurable policies with $(\omega_0, {\bf y}) \sim \mu$; that is, $P \in L_A(\mu) \subset {\cal P}\bigg(\Omega_0 \times \prod_{k=1}^N (\mathbb{Y}^k \times \mathbb{U}^k)\bigg)$ if and only if
\begin{align}\label{measurableStrategic}
P(B) = \int_{B^0 \times \prod_{k=1}^N A^k} \mu(d\omega_0, d{\bf y}) \, \prod_{k=1}^N 1_{\{u^k = \gamma^k(y^k) \in B^k\}},  
\end{align}
for all cylindrical $B \in {\cal B}\left(\Omega_0 \times \prod_{k=1}^N (\mathbb{Y}^k \times \mathbb{U}^k)\right)$, where $\gamma^k \in \Gamma^k$ for $k=1,\ldots,N$. Let $L_A(\mu,\underline{\gamma})$ be the strategic measure under a particular strategy $\underline{\gamma} \in {\bf \Gamma}$.

The first relaxation is obtained via individual randomization of policies. Namely, let $L_R(\mu)$ be the set of strategic measures induced by all individually randomized team policies where $\omega_0, {\bf y} \sim \mu$; that is, 
\[L_R(\mu) := \bigg\{P \in {\cal P}\bigg(\Omega_0 \times \prod_{k=1}^N (\mathbb{Y}^k \times \mathbb{U}^k)\bigg): P(B) = \int_{B} \mu(d\omega_0, d{\bf y}) \prod_{k=1}^N \Pi^k(d u^k| y^k)\bigg\},\]
where $\Pi^k$ takes place from the set of stochastic kernels from $\mathbb{Y}^k$ to $\mathbb{U}^k$ for each $k=1,\ldots,N$.

Another relaxation, which is stronger than the former one, is obtained by taking the mixture of the elements of $L_A(\mu)$. To this end, define $\Upsilon = [0,1]^{N}$. We then let 
\begin{eqnarray}\label{LC}
L_{C}(\mu) := \bigg\{P \in {\cal P}\bigg(\Omega_0 \times \prod_{k=1}^N (\mathbb{Y}^k \times \mathbb{U}^k)\bigg): P(B) = \int \eta(dz) L_A(\mu,\underline{\gamma}(z))(B), \, \eta \in{\cal P}(\Upsilon) \bigg\},  \nonumber 
\end{eqnarray}
where $\underline{\gamma}(z)$ denotes a collection of team policies measurably parametrized by $z \in \Upsilon$ so that the map $L_A(\mu,\underline{\gamma}(\cdot)): \Upsilon \to L_A(\mu)$ is Borel measurable as $L_A(\mu)$ is a Borel subset of ${\cal P}\bigg(\Omega_0 \times \prod_{k=1}^N (\mathbb{Y}^k \times \mathbb{U}^k)\bigg)$ under weak convergence topology (as we will see in Theorem \ref{FeinbergStrategicMeasurability3}).

Let $L_{CR,\eta}$ denote the set of strategic measures that are induced by some \textit{fixed} but common independent randomness (with probability measure $\eta$) and arbitrary private independent randomness; that is,
\begin{eqnarray}
L_{CR,\eta}(\mu)&:=& \bigg\{P \in {\cal P}\bigg(\Omega_0 \times \prod_{k=1}^N (\mathbb{Y}^k \times \mathbb{U}^k)\bigg): \nonumber \\
&&  \quad \quad \quad \quad \quad P(B) = \int_{B \times \Upsilon} \eta(dz) \mu(d\omega_0, d{\bf y}) \prod_k \Pi^k(d u^k| y^k,z) \bigg\}, \nonumber 
\end{eqnarray}
where $\Pi^k$ takes place from the set of stochastic kernels from $\mathbb{Y}^k \times \Upsilon$ to $\mathbb{U}^k$ for each $k=1,\ldots,N$. Here, the common randomness $\eta$ is fixed.

Let $L_{CCR}$ denote the set of strategic measures that are induced by some arbitrary but common independent randomness and arbitrary private independent randomness, as in $L_{C}(\mu)$; that is,
\begin{eqnarray}
L_{CCR}(\mu)&:=& \bigg\{P \in {\cal P}\bigg(\Omega_0 \times \prod_{k=1}^N (\mathbb{Y}^k \times \mathbb{U}^k)\bigg): \nonumber \\
&&  \quad \quad \quad \quad \quad P(B) = \int_{B \times \Upsilon} \eta(dz) \mu(d\omega_0, d{\bf y}) \prod_k \Pi^k(d u^k| y^k,z), \, \eta \in{\cal P}(\Upsilon) \bigg\}, \nonumber 
\end{eqnarray}
where $\Pi^k$ takes place from the set of stochastic kernels from $\mathbb{Y}^k \times \Upsilon$ to $\mathbb{U}^k$ for each $k=1,\ldots,N$. Here, the common randomness $\eta$ is arbitrary, unlike $L_{CR,\eta}(\mu)$. The following result, essentially from \cite{YukselSaldiSICON17}, states some structural results about above-defined sets of strategic measures. In particular, it establishes convexity related properties of these sets. 

\begin{theorem}\label{FeinbergStrategic}
Consider a static team problem. Then, we have the following characterizations. \\
$(i)$ $L_R(\mu)$ has the following representation:
\begin{align}\label{representLR}
&L_{R}(\mu) = \bigg\{P \in {\cal P}\bigg(\Omega_0 \times \prod_{k=1}^N (\mathbb{Y}^k \times \mathbb{U}^k)\bigg): P(B) = \int U(dz) L_A(\mu,\underline{\gamma}(z))(B), \nonumber \\
& \phantom{xxxxxxxxxxxxxxxxxxx}   U \in{\cal P}(\Upsilon), U(dv_1,\cdots,dv_N) = \prod_s \eta_k(dv_k), \eta_k \in {\cal P}([0,1]) \bigg\}; \nonumber
\end{align}
that is, $U \in {\cal P}(\Upsilon)$ is constructed by the product of $N$ independent random variables on $[0,1]$.\\
$(ii)$  $L_{C}(\mu)=L_{CCR}(\mu)$ and this is a convex set. The set of extreme points of $L_{C}(\mu)$ is $L_A(\mu)$. Furthermore, $L_R(\mu) \subset L_C(\mu)$. \\
$(iii)$ We have the following equalities:
\[\inf_{\underline{\gamma} \in {\bf \Gamma}} J(\underline{\gamma}) = \inf_{P \in L_A(\mu)} \int P(ds) c(s) = \inf_{P \in L_R(\mu)} \int P(ds) c(s) = \inf_{P \in L_{C}(\mu)} \int P(ds) c(s).\]
In particular, deterministic policies are optimal among the randomized class. In other words, individual and common randomness does not improve the optimal team cost. \\
$(iv)$ The sets $L_R(\mu)$ and $L_{CR,\eta}(\mu)$ are not convex. In particular, the presence of independent or (fixed) common randomness does not convexify the set of strategic measures.\\
$(v)$ $L_R(\mu)$ and $L_C(\mu)$ are not necessarily weakly closed.
\end{theorem}
\sy{
With $\mathbb{X}, \mathbb{Y}$ standard Borel spaces, let us recall the following: we say that a sequence of probability measures (on $\mathbb{X}$) $P_n \to P$ {\it weakly} if $\int P_n(dx) f(x) \to \int P(dx) f(x)$ for all $f \in C_b(\mathbb{X})$. If this holds for every $f$ measurable and bounded, we say convergence holds {\it setwise}. }

To establish (v), we will present two counterexamples below in Theorem~\ref{counterEx} and Theorem~\ref{counterExLC}.

In the economics and game theory literature, information structures are also studied extensively. Stochastic team problems are termed as {\it identical interest games}. In this literature, $L_C(\mu)$ appears in the analysis of Aumann's correlated equilibrium \cite{aumann1987correlated}. Common and independent randomness discussions appear in the analysis of comparison of information structures \cite{LehrerRosenbergShmaya}. For further discussions, including a multi-stage generalization known as communication equilibria, see \cite{forges1986approach}. For a detailed treatment, we refer the reader to \cite[p. 131]{mertens2015repeated}.

\subsection{Sets of strategic measures for dynamic teams in the absence of static reduction}

Note that if the dynamic team setup admits a static reduction (in particular independent static reduction), then one can define strategic measures by considering equivalent static problem and characterize the convexity properties of the set of strategic measures, as done in the previous section. In this section, we suppose that dynamic team does not admit a static reduction. Let $\mu$ be the distribution of $\omega_0$. Recall that in dynamic setup, the distribution of measurements ${\bf y}$ is not fixed as opposed to the static case. In this case, we present the following characterization for strategic measures in dynamic sequential teams. Let, for all $n \in {\cal N}$, 
\[h_n = \{\omega_0,y^1,u^1,\cdots,y^{n-1},u^{n-1},y^n,u^n\},\]
and $p_n(dy^n|h_{n-1}) := P(dy^n | h_{n-1})$ be the transition kernel characterizing the measurements of DM~$n$ according to the intrinsic model. We note that this may be obtained by the relation:
\sy{\begin{eqnarray}
&& p_n(y^n \in \cdot | \, \omega_0,y^1,u^1,\cdots,y^{n-1},u^{n-1}) \nonumber \\
&& \, := P\bigg(\eta^i(\omega,{\bf u}^{[1,i-1]}) \in \cdot  \bigg| \, \omega_0,y^1,u^1,\cdots,y^{n-1},u^{n-1}\bigg) \nonumber \\
&& \, = P\bigg(g^n(\omega_0,\omega_n, u^{1},\cdots,u^{n-1}) \in \cdot  \bigg| \, \omega_0,y^1,u^1,\cdots,y^{n-1},u^{n-1}\bigg). \label{kernelDefn}
\end{eqnarray}
}
Note that once a policy is fixed, $p_n(dy^n|h_{n-1})$ represents the conditional distribution of $y^n$ given the past history $h_{n-1}$. Let $L_A(\mu)$ be the set of strategic measures induced by measurable policies and let $L_R(\mu)$ be the set of strategic measures induced by individually randomized policies for the dynamic team. We have the following characterizations of $L_A(\mu)$ and $L_R(\mu)$ that are quite useful when establishing the closedness of these sets. 

\smallskip

\begin{theorem}[{\cite[Theorem 2.2]{YukselSaldiSICON17}}]\label{StrategicCharacterization} 
Consider a dynamic team problem that does not admit a static reduction. Then, we have the following characterizations. \\
$(i)$ A probability measure $P \in {\cal P}\bigg(\Omega_0 \times \prod_{k=1}^N (\mathbb{Y}^k \times \mathbb{U}^k)\bigg)$ is a strategic measure induced by a measurable policy (that is in $L_A(\mu)$) if and only if, for every $n=1,\ldots,N$, we have
\[\int P(dh_{n-1},dy^n) \, g(h_{n-1},y^{n}) = \int P(dh_{n-1}) \, \bigg(\int_{\mathbb{Y}^n} g(h_{n-1},z) \, p_n(dz | h_{n-1}) \bigg)\]
and
\[\int P(dh_n) \, g(h_{n-1},y^n,u^{n}) = \int P(dh_{n-1},dy^n) \, \bigg(\int_{\mathbb{U}^n} g(h_{n-1},y^{n},a) \, 1_{\{\gamma^n(y^n) \in da\}} \bigg),\]
for all continuous and bounded function $g$ with appropriate arguments, where $P(d\omega_0) = \mu(dw_0)$ and $\gamma^n \in \Gamma^n$. \\
$(ii)$ A probability measure $P \in {\cal P}\bigg(\Omega_0 \times \prod_{k=1}^N (\mathbb{Y}^k \times \mathbb{U}^k)\bigg)$  is a strategic measure induced by a individually randomized policy (that is in $L_R(\mu)$)  if and only if, for every $n = 1,\ldots,N$, we have
\begin{eqnarray}\label{convD1}
\hspace{-20pt} \int P(dh_{n-1},dy^n) \, g(h_{n-1},y^{n}) = \int P(dh_{n-1}) \, \bigg(\int_{\mathbb{Y}^n} g(h_{n-1},z) \, p_n(dz | h_{n-1}) \bigg)  
\end{eqnarray}
and
\begin{eqnarray}\label{convD2}
\hspace{-20pt} \int P(dh_n) \, g(h_{n-1},y^n,u^{n}) = \int P(dh_{n-1},dy^n) \, \bigg(\int_{\mathbb{U}^n} g(h_{n-1},y^{n},a^n) \, \Pi^n(da^n | y^n) \bigg)  
\end{eqnarray}
for all continuous and bounded function $g$ with appropriate arguments, where $P(d\omega_0) = \mu(dw_0)$ and $\Pi^n$ is a stochastic kernel from $\mathbb{Y}^n$ to $\mathbb{U}^n$. 
\end{theorem}

\begin{remark}
A result similar to Theorem \ref{FeinbergStrategic} can also be stated for the dynamic case, in particular with regard to $L_A(\mu)$ being the set of extreme points of the convex hull of $L_R(\mu)$. The reader is referred to \cite[Theorem 2.3]{YukselSaldiSICON17} which essentially establishes this; see also \cite[Theorem 1.c]{Feinberg1} for related discussions.

\sy{A celebrated result in economics theory, known as Kuhn's theorem \cite{kuhn1953extensive}, notes that the convex hull of admissible (i.e. those in $L_A(\mu)$) strategic measures (hence $L_C(\mu)$) is equivalent to $L_R(\mu)$ when the information structure is classical. We can thus state that this does not apply in the absence of classical-ness, as $L_R(\mu)$ would not be convex (if the information structure is not classical, then convexity fails \cite[p.12]{YukselSaldiSICON17}), but the convex hull of admissible policies is, by definition, convex; but the convex hull of $L_R(\mu)$ is $L_C(\mu)$. }

\end{remark}

%

\subsection{Measurability properties of sets of strategic measures}

As noted earlier, the set $L_A(\mu)$ is a Borel subset of ${\cal P}\left(\Omega_0 \times \prod_k (\mathbb{Y}^k \times \mathbb{U}^k)\right)$ under weak convergence topology. The same is true for $L_R(\mu)$, which is stated in the following theorem. This result will be crucial in the analysis to follow.

\smallskip

\begin{theorem}[{\cite[Theorem 2.10]{YukselSaldiSICON17}}]\label{FeinbergStrategicMeasurability3}
Consider a sequential (static or dynamic) team. 
\begin{itemize}
\item[(i)] The set of strategic measures $L_R(\mu)$ is Borel when viewed as a subset of the space of probability measures on $\Omega_0 \times \prod_{k=1}^N (\mathbb{Y}^k \times \mathbb{U}^k)$ under the topology of weak convergence.
\item[(ii)] The set of strategic measures $L_A(\mu)$ is Borel when viewed as a subset of the space of probability measures on $\Omega_0 \times \prod_{k=1}^N (\mathbb{Y}^k \times \mathbb{U}^k)$ under the topology of weak convergence.
\end{itemize}
\end{theorem}


For further properties of the sets of strategic measures, see \cite{YukselSaldiSICON17}.

\section{Relative compactness and closedness of strategic measures, and existence of optimal policies}\label{regularity} 

Existence of optimal policies for static and a class of sequential dynamic teams have been studied recently in \cite{gupta2014existence,YukselSaldiSICON17,SaldiArXiv2017}. More specific setups have been studied in \cite{WuVer11}, \cite{wit68}, \cite{YukselOptimizationofChannels} and \cite{YukselBasarBook}. Existence of optimal team policies has been established in \cite{charalambous2017centralizedI} for a class of continuous-time decentralized stochastic control problems. For a class of teams which are convex, one can reduce the search space to a smaller parametric class of policies, such as linear policies for quasi-classical linear quadratic Gaussian problems \cite{rad62,KraMar82,HoChu}. 

The following theorem states a general existence result for static teams and for dynamic teams admitting static reduction. Its proof depends on Weierstrass Extreme Value Theorem.  

\begin{theorem}\label{existenceT}
Consider a static team or the static reduction of a dynamic team with $c$ denoting the cost function. Let $c$ be lower semi-continuous in ${\bf u}$ for every fixed $\omega_0, {\bf y}$ and $L_R(\mu)$ or $L_C(\mu)$ be a compact set under weak convergence topology. Then, there exists an optimal team policy. This policy can be chosen deterministic and hence induces a strategic measure in $L_A(\mu)$.  
\end{theorem}

\begin{remark}
Since the cost function $c_s$ in independent static reduction of a dynamic team also depends on the measurements ${\bf y}$, we include ${\bf y}$ as an argument to the cost function $c$ in the previous theorem.
\end{remark}

By Theorem~\ref{existenceT}, to prove the existence of optimal team policy, it is sufficient to establish the lower semi-continuity of the cost function $c$ and compactness of $L_R(\mu)$ or $L_C(\mu)$ or any other subset of $L_C(\mu)$ that is sufficient for optimality. However, we recall that unless certain conditions are imposed, the conditional independence property dictated by IS is not preserved under weak convergence (indeed this is also true even for setwise convergence which is a stronger convergence notion). Hence, $L_R(\mu)$ and $L_C(\mu)$ are in general not closed and so compact. The following theorem establishes that $L_R(\mu)$ is not closed under the weak convergence topology.
 
\begin{theorem}[{\cite[Theorem 2.7]{YukselSaldiSICON17}}] \label{counterEx}
Consider a sequence of probability measures $P_n \in {\cal P}(\mathbb{U}^1 \times \mathbb{Y} \times \mathbb{U}^2)$ so that for all $n$: 
\[P_n(du^1| y, u^2) = P_n(du^1|y).\]
If $P_n \to P$ setwise (and thus also weakly), it is not necessarily the case that 
\[P(du^1| y, u^2) = P(du^1|y).\]
That is, conditional independence of $u^1$ and $u^2$ given $y$ is not preserved under setwise convergence. In particular, $L_R(\mu)$ is not (weakly or setwise) closed.
\end{theorem}

Before we present the proof of Theorem~\ref{counterEx}, let us note that, for each $n$, we can view $P_n$ as a strategic measure of the following static team problem. In this team problem, $\Omega$ is a degenerate space and we have two DMs that are sharing the same measurement $y$. Therefore, Theorem~\ref{counterEx} states that the set of strategic measures $L_R(\mu)$ is not closed under setwise (and so weak convergence) topology.   

\begin{proof} It suffices to provide a counterexample. We build on an example from \cite{YukselOptimizationofChannels} (used in a different context) in the following. Let $\mathbb{Y}=[0,1]$, $\mathbb{U}^1=\mathbb{U}^2=\{0,1\}$, and $y \sim m$ where $m$ the Lebesgue measure (uniform
distribution) on $[0,1]$. Let
\begin{eqnarray}\label{star2}
  L_{nk}= \left[\frac{2k-2}{2n},\frac{2k-1}{2n}\right), \quad
    R_{nk}= \left[\frac{2k-1}{2n},\frac{2k}{2n}\right)
\end{eqnarray}
and define the {\it square wave} function
\[
  h_n(t) = \sum_{k=1}^n\bigl(1_{\{t \in L_{nk}\}} - 1_{\{t \in R_{nk}\}} \bigr).
\]
Define further $f_n(t)=h_n(t)+1$. 

 \begin{figure}
      \centering
      \begin{tikzpicture}
        \begin{groupplot}[group style={group size=3 by 1},xmin=0,xmax=1,ymin=-0.2,ymax=2.2,height=4.5cm,width=4.5cm,ytick={0,1,2},xtick={0,0.5,1},no markers]
          \nextgroupplot[title={$f_1$}];
          \addplot[domain=0:0.5,blue] {2};
          \addplot[domain=0.5:1,blue] {0};
          \draw[dotted] (axis cs:0.5,2) -- (axis cs:0.5,0);
          \nextgroupplot[title={$f_2$}];
          \addplot[domain=0:0.25,blue] {2};
          \addplot[domain=0.25:0.5,blue] {0};
          \addplot[domain=0.5:0.75,blue] {2};
          \addplot[domain=0.75:1,blue] {0};
          \draw[dotted] (axis cs:0.25,2) -- (axis cs:0.25,0);
          \draw[dotted] (axis cs:0.5,2) -- (axis cs:0.5,0);
          \draw[dotted] (axis cs:0.75,2) -- (axis cs:0.75,0);
          \nextgroupplot[title={$f_3$}];
          \addplot[domain=0:1/6,blue] {2};
          \addplot[domain=1/6:2/6,blue] {0};
          \addplot[domain=2/6:3/6,blue] {2};
          \addplot[domain=3/6:4/6,blue] {0};
          \addplot[domain=4/6:5/6,blue] {2};
          \addplot[domain=5/6:1,blue] {0};
          \draw[dotted] (axis cs:1/6,2) -- (axis cs:1/6,0);
          \draw[dotted] (axis cs:2/6,2) -- (axis cs:2/6,0);
          \draw[dotted] (axis cs:3/6,2) -- (axis cs:3/6,0);
          \draw[dotted] (axis cs:4/6,2) -- (axis cs:4/6,0);
          \draw[dotted] (axis cs:5/6,2) -- (axis cs:5/6,0);
        \end{groupplot}
      \end{tikzpicture}
      \caption{Plots of probability density functions $f_1$, $f_2$, and $f_3$.}\label{fig:squarewave}
    \end{figure}
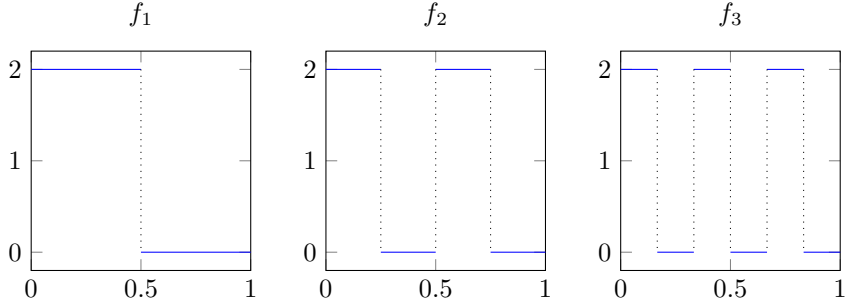

Let
  $B_{n,1} = \bigcup_{k=1}^n L_{nk}$ and $B_{n,2}= [0,1]\setminus
  B_{n,1}$. Define $\{Q_n\}$ as the sequence of $2$-cell quantizers
  given by
\[
Q_n(1|y)=1_{\{y \in B_{n,1}\}}, \quad Q_n(0|y)=1_{\{y\in  B_{n,2}\}} .
\]
Let \[P_n(u^1 = a | y) = P_n(u^2 = a | y) = Q_n(a|y),\]
where $a \in \{0,1\}$. Define $P \in {\cal P}(\mathbb{U}^1 \times \mathbb{Y} \times \mathbb{U}^2)$ as $P(a,A,b) = 1_{\{a=b\}} \frac{1}{2}m(A)$, where $a,b \in \{0,1\}$ and  $A\in \mathcal{B}([0,1])$. 

By the proof of the Riemann-Lebesgue lemma (\cite{WhZy77}, Thm.\ 12.21), observe that for all $A\in \mathcal{B}([0,1])$,
\[
\lim_{n \to \infty} \int_{A} Q_n(1|y)m(dy) = \lim_{n\to \infty} \int_{A} \frac{1}{2}
f_n(t) \, m(dt) = \frac{1}{2}m(A),
\]
and thus for all $A\in \mathcal{B}([0,1])$
\begin{eqnarray}
&&\lim_{n \to \infty} P_n(u^1=1,y \in A,u^2=1) \nonumber \\
&&=\lim_{n \to \infty} \int_{A}   P_n(u^1=1|y)P_n(u^2=1|y)m(dy)  \nonumber \\
&&= \lim_{n \to \infty} \int_{A}   P_n(u^1=1|y)m(dy) \nonumber \\
&& = \frac{1}{2}m(A) \nonumber \\
&&= P(1,A,1)
\end{eqnarray}
A similar property applies for $(u^1,u^2)=(0,0), (0,1)$ and $(1,0)$ so that 
\[\lim_{n \to \infty} P_n(u^1=a,y \in A,u^2=b) \to P(a,A,b) = 1_{\{a=b\}} \frac{1}{2}m(A) \]
Thus, $P_n \to P$ setwise. But even though $P_n$ satisfies the conditional independence property that $P_n(u^1=1|y,u^2) = Q_n(1|y)$, $P$ does not satisfy the conditional independence property of $u^1$ and $u^2$ given $y$: Under $P$, $u^1$ and $y$ are independent but $u^1=u^2$ and thus $P(u^1=a|y,u^2=b) = 1_{\{a=b\}} \neq {1 \over 2} = P(u^1=a|y)$. Thus, setwise (and hence weak) convergence does not preserve the conditional independence property. 
\end{proof}

We can also show via counterexample that $L_C(\mu)$ is not closed either with respect to the weak convergence topology.

\begin{theorem} \label{counterExLC}
$L_C(\mu)$ is not weakly or setwise closed.
\end{theorem}

\begin{proof}
Consider the same example as in the proof of Theorem \ref{counterEx}, modified as follows.
Let $\mathbb{Y}=[0,1] \cup (1,2]$, $\mathbb{U}^1=\mathbb{U}^2=\{0,1,2,3\} $. Let $x$ be a $\{0,1\}$-valued uniformly distributed Bernoulli random variable and let $y \sim 1_{\{x=0\}} v +  1_{\{x=1\}}(v+1)$ where $v \sim m$ with $m$ being the Lebesgue measure (uniform
distribution) on $[0,1]$. 
Let $L_{nk}$ and $R_{nk}$ be as given in (\ref{star2})
and define the {\it square wave} function, as before,
\[
  h_n(t) = \sum_{k=1}^n\bigl(1_{\{t \in L_{nk}\}} - 1_{\{t \in R_{nk}\}} \bigr).
\]
Define also $f_n(t)=h_n(t)+1$. Let $B_{n,1} = \bigcup_{k=1}^n L_{nk}$ and $B_{n,2}= [0,1]\setminus
  B_{n,1}$. Define $\{Q_n\}$ as the sequence of $4$-cell quantizers given by
\[
Q_n(1|y)=1_{\{y \in B_{n,1}\}}, \quad Q_n(0|y)=1_{\{y\in  B_{n,2}\}} , \qquad \qquad y \in [0,1]
\]
and
\[
Q_n(3|y)=1_{\{y \in 1+ B_{n,1}\}}, \quad Q_n(2|y)=1_{\{y\in  1+ B_{n,2}\}} , \qquad \qquad y \in (1,2].
\]

\noindent Let \[P_n(u^1 = a | y) = P_n(u^2 = a | y) = Q_n(a|y),\]
where $a \in \{0,1,2,3\}$. Define $P \in {\cal P}(\mathbb{U}^1 \times \mathbb{Y} \times \mathbb{U}^2)$ as 
\[P(a,A,b) = \frac{1}{4} \, 1_{\{a=b, a \in \{0,1\} \}} \, m(A \cap [0,1]) + \frac{1}{4} 1_{\{a=b, a \in \{2,3\} \}} \, m(A \cap (1,2]) \]
 where $a,b \in \{0,1,2,3\}$ and  $A\in \mathcal{B}([0,2])$. 
 
We can show that $P_n \to P$ setwise (and so weakly). In the limit, it is the case that $u^1=u^2$ almost surely. However, notice that $P \notin L_C(\mu)$: $P$ is not a mixture of conditionally independent random variables given $y$, with the mixture being independent of $y$. Here, the mixture representation of $P$ in terms of the two independent $\{0,1\}$-valued and $\{2,3\}$-valued random variables is so that the mixing random variable is correlated with $y$ (and gives information on $y$).
\end{proof}

We refer the reader to \cite{aldous1981weak} \cite{hellwig1996sequential} \cite[Theorem 1.1]{backhoff2019all} \cite{YukselOptimizationofChannels}, \cite{barbie2014topology} for further related results on such intricacies on conditional independence properties. A sufficient condition for compactness of $L_R(\mu)$ under the weak convergence topology was reported in \cite{gupta2014existence}. We re-state this result in a brief and different form below for reader's convenience.

\begin{theorem}[\cite{gupta2014existence}]\label{SuffCon1} 
Consider a static team where the action sets $\mathbb{U}^i$ are compact. Furthermore, the measurements 
\[P(d{\bf y}|\omega_0) = \prod_{i=1}^N g^i(dy^i|\omega_0)\]
satisfy $g^i(dy^i|\omega_0) = \xi^i(y^i,\omega_0) \, \nu^i(dy^i)$ for some positive measure $\nu^i$ and continuous $\xi^i$ such that, for every $\epsilon >0$, there exists $\delta > 0$ so that for $\rho_i(a,b)<\delta$ (where $\rho_i$ is a metric on $\mathbb{Y}^i$)
\[|\xi^i(b,\omega_0) - \xi^i(a,\omega_0) | \leq \epsilon \, h^i(a,\omega_0),\]
with $\sup_{\omega_0} \int h^i(a,\omega_0) \, \nu^i(dy^i) < \infty$. Hence, static team admits independent static reduction. Then, the set $L_R(\mu)$ is compact under weak convergence topology. Therefore, if $c(\omega_0,{\bf u})$ is lower semi-continuous in ${\bf u}$ for any $\omega_0$, then there exists an optimal team policy (which is deterministic and hence in $L_A(\mu)$).
\end{theorem}

The results in \cite{gupta2014existence} also apply to static reduction of sequential dynamic teams, and a class of teams with non-compact action spaces that however has moment-type cost function leading to a tightness condition on the set of strategic measures with a finite cost. In particular, the existence result applies to the celebrated counterexample of Witsenhausen \cite{wit68}. 

However, the result above can be significantly relaxed where the relaxation is not only in the compactness condition, this can be modified by the usual tightness conditions which will also be presented below. The generalization is with respect to the topology considered: in the following, we present the most general conditions to our knowledge on existence of optimal policies. In the following, we slightly strengthen \cite[Theorem 5.2]{YukselWitsenStandardArXiv} to allow for lower semi-continuity of the cost function only in the actions.

\begin{theorem}[{\cite[Theorem 5.2]{YukselWitsenStandardArXiv}}]\label{existenceRelaxed3}
Consider a static or a dynamic team that admits an independent static reduction. Let $c$ be lower semi-continuous in ${\bf u}$ for any $\omega_0, {\bf y}$. Suppose further that $\mathbb{U}^i$ is $\sigma-$compact (that is, $\mathbb{U}^i=\cup_n K_n$ for a countable collection of increasing compact sets $K_n$) and, without any loss, the control laws can be restricted to those with $E[\phi^i(u^i)]\leq M$ for some lower semi-continuous $\phi^i: \mathbb{U}^i \to \mathbb{R}_+$ which satisfies $\lim_{n \to \infty} \inf_{u^i \notin K_n} \phi^i(u^i) = \infty$.
Then, an optimal team policy exists. 
\end{theorem}

\begin{proof} For each $m$, consider the strategic measure: 
\[P(d\omega_0) \prod_{k=1}^N(Q^k\gamma^k)_m(dy^k,du^k).\]
Suppose that $(Q^k\gamma^k)_m(dy^k,du^k)$, converges to $(Q^k\gamma^k)(dy^k,du^k)$ weakly. By \cite[p. 57]{Par67}, the joint product measure $P(d\omega_0) \prod_{k=1}^N(Q^k\gamma^k)_m(dy^k,du^k)$
will converge weakly to $P(d\omega_0) \prod_{k=1}^N(Q^k\gamma^k)(dy^k,du^k)$; see also \cite[Section 5]{serfozo1982convergence}.
Recall the $w$-$s$ topology introduced by Sch\"al \cite{schal1975dynamic}: The $w$-$s$ topology on the set of probability measures ${\cal P}(\mathbb{X} \times \mathbb{U})$ is the coarsest topology under which $\int f(x,u) \nu(dx,du): {\cal P}(\mathbb{X} \times \mathbb{U}) \to \mathbb{R}$ is continuous for every measurable and bounded $f$ which is continuous in $u$ for every $x$ (but unlike weak topology, $f$ does not need to be continuous in $x$). Since the marginals on $\prod_k \mathbb{Y}^k$ is fixed, \cite[Theorem 3.10]{schal1975dynamic} (or  \cite[Theorem 2.5]{balder2001}) establishes that the set of all probability measures with a fixed marginal on $\prod_k \mathbb{Y}^k$ is (sequentially) compact under the $w$-$s$ topology. If the function $c_s(\omega_0,{\bf y},\cdot)$ were continuous, this in turn would ensure that the function (by a truncation and then a limiting argument)
\[\int P(d\omega_0) \prod_{k=1}^N(Q^k\gamma^k)(dy^k,du^k) \, c_s(\omega_0,{\bf y},{\bf u}),\]
is lower semi-continuous under the $w$-$s$ topology. Now, by approximating the lower semi-continuous function $c_s(\omega_0,{\bf y},\cdot)$ from below by continuous functions, applying the argument above, and taking the limit, we conclude that the lower semi-continuity also applies when $c_s(\omega_0,{\bf y},\cdot)$ is only lower semi-continuous.

Finally, since the set of admissible strategic measures is sequentially compact under the $w$-$s$ topology,  existence of an optimal team policy follows. The proof for dynamic case follows analogously.
\end{proof}

\smallskip

\begin{remark}\label{detRemark}
Building on \cite[Theorems 2.3 and 2.5]{YukselSaldiSICON17} and \cite[p. 1691]{gupta2014existence} (due to Blackwell's theorem on irrelevant information  \cite{Blackwell2,Blackwell3}, \cite[p. 457]{YukselBasarBook}), an optimal policy, when exists, can be assumed to be {\it deterministic}.
\end{remark}

\smallskip

Note that Theorem \ref{existenceRelaxed3} provide weaker conditions when compared with Theorem \ref{SuffCon1}. So far, we present existence results for static or dynamic teams that admit independent static reduction. In the following, we present existence results for teams that do not admit independent static reduction. 

\begin{theorem}[{\cite[Theorem 2.9]{YukselSaldiSICON17}}]\label{SuffCon2}
Consider a sequential team with a classical information structure with the further property that $\sigma(\omega_0) \subset \sigma(y^1)$ (under every policy, $y^1$ contains $\omega_0$). Suppose further that $\prod_{k=1}^N \mathbb{U}^k$ is compact. If $c$ is lower semi-continuous and each of the kernels $p_n$ (defined in (\ref{kernelDefn})) is weakly continuous so that
\begin{eqnarray}\label{weakConKer}
\int f(y^n) \, p_n(dy^n | \omega_0,y^{1},\ldots,y^{n-1},u^{1},\cdots,u^{n-1})
\end{eqnarray} 
is continuous in $\omega_0,y^{1},\cdots,y^{n-1},u^{1},\cdots,u^{n-1}$ for every continuous and bounded $f$, then there exists an optimal team policy which is deterministic.
\end{theorem}

A further existence result along similar lines, for a class of static teams, is presented next.

\begin{theorem}[{\cite[Theorem 5.6]{YukselWitsenStandardArXiv}}]\label{SuffCon2'''}
Consider a static team with a classical information structure (that is, with an expanding information structure so that $\sigma(y^n) \subset \sigma(y^{n+1}), n \geq 1$). Suppose further that $\prod_{k=1}^N (\mathbb{Y}^k \times \mathbb{U}^k)$ is compact. If 
\[\tilde{c}(y^1,\cdots,y^N,u^1,\cdots,u^N):=E[c(\omega_0,{\bf u}) | {\bf y}, {\bf u}]\]
is lower semi-continuous in ${\bf u}$ for every ${\bf y}$, then there exists an optimal team policy which is deterministic.
\end{theorem}

\begin{proof} Different from Theorem \ref{SuffCon2}, we eliminate the use of $\omega_0$, and study the properties of the set of strategic measures. Different from \cite[Theorem 3.5]{YukselSaldiSICON17}, we relax weak continuity. Once again, instead of the weak topology, we will use the $w$-$s$ topology \cite{Schal}.

Note that when $\prod_{k=1}^N \mathbb{U}^k$ is compact, the set of all probability measures on $\prod_{k=1}^N \mathbb{Y}^k \times \mathbb{U}^k$ forms a weakly compact set. Since the marginals on $\prod_{k=1} \mathbb{Y}^k$ is fixed, \cite[Theorem 3.10]{Schal} (or  \cite[Theorem 2.5]{balder2001}) establishes that the set of all probability measures with a fixed marginal on $\prod_{k=1} \mathbb{Y}^k$ is relatively compact under the $w$-$s$ topology. Therefore, it suffices to ensure the closedness of the set of strategic measures, which leads to the sequential compactness of the set under this topology. To facilitate such a compactness condition, we first {\it expand} the information structure so that DM $k$ has access to all the previous actions $u^{1},\cdots,u^{k-1}$ as well. Later on, we will see that this expansion is redundant. With this expansion, any $w$-$s$ converging sequence of strategic measures will continue satisfying (\ref{convD2}) in the limit due to the fact that there is no conditional independence property in the sequence since all the information is available at DM $k$. That is, $P_n(du^n | y^n,y_{[0,n-1]},u_{[0,n-1]})$ satisfies the conditional independence property trivially as all the information is available. On the other hand, for each element in the sequence of conditional probability measures, the conditional probability for the measurements writes as
$P(dy^n | y_{[0,n-1]},u_{[0,n-1]}) = P(dy^n | y^{n-1})$.  We wish to show that this also holds for the $w$-$s$ limit measure. Now, we have that for every $n$, $y^n \leftrightarrow y^{n-1} \leftrightarrow h_{n-1}$ forms a Markov chain. By considering the convergence properties only on continuous functions and bounded $f$, as in (\ref{convD1}), with $P_m \to P$ weakly, we have that
\begin{eqnarray}
&&\int P(dy^n | y^{n-1}) P_m(dh_{n-1}) f(y^n,h_{n-1}) = \int \bigg( P(dy^n | y^{n-1})  f(y^n,h_{n-1}) \bigg) P_m(dh_{n-1}) \nonumber \\
&& \,\to  \int \bigg( P(dy^n | y^{n-1})  f(y^n,h_{n-1}) \bigg) P(dh_{n-1}) =  \int P(dy^n | y^{n-1}) P(dh_{n-1}) f(y^n,h_{n-1}) \nonumber
\end{eqnarray}
Here, $\bigg( P(dy^n | y^{n-1})   f(y^n,h_{n-1}) \bigg)$ is not continuous in $y^{n-1}$, but it is in $h_{n-2}$ and $u^{n-1}$ by an application of the dominated convergence theorem, and since \[P_m(dy^1,\ldots,dy^{n-1},du^{1},\ldots,du^{n-1}) \to P(dy^1,\ldots,dy^{n-1},du^{1},\ldots,du^{n-1}),\]
 in the $w$-$s$ sense (setwise in the measurement variable coordinates), convergence holds. Thus, (\ref{convD1}) is also preserved. Hence, for any $w$-$s$ converging sequence of strategic measures satisfying (\ref{convD1})-(\ref{convD2}) so does the limit since the team is static and with perfect-recall. By \cite[Theorem 3.7]{Schal}, and the generalization of Portmanteau theorem for the $w$-$s$ topology, the lower semi-continuity of $\int \mu(d{\bf y}, d{\bf u}) \tilde{c}({\bf y},{\bf u})$ over the set of strategic measures leads to the existence of an optimal strategic measure. As a result, the existence follows from similar steps to that of Theorem \ref{existenceT}. Now, we know that an optimal policy will be deterministic (see Remark \ref{detRemark}). Thus, a deterministic policy may not make use of randomization, therefore DM $k$ having access to $\{y^{k},y^{k-1},y^{k-2}, \cdots\}$ is informationally equivalent to him having access to $\{y^{k},(y^{k-1},u^{k-1}),(y^{k-2},u^{k-2})\}$ for an optimal policy. Thus, an optimal team policy exists.
\end{proof}

We can report that Theorem \ref{existenceRelaxed3} (for static teams or dynamic teams with an independent static reduction) and Theorems \ref{SuffCon2} and \ref{SuffCon2'''} (for sequential teams that do not allow an independent static reduction) are the most general existence results, to our knowledge, for sequential team problems considered in this paper. These results complement each other and cover a very large class of decentralized stochastic control problems. 

In the following, we will discuss lower and upper bounds on the optimal costs, numerical programs, and convex relaxations.

\section{Convex Relaxations}\label{convexRelax}

In this section, we introduce convex relaxations for static teams or dynamic teams that admit independent static reduction through relaxing conditional independence among actions. Therefore, in this section, the joint distribution of $\omega_0,{\bf y}$ is in product form $\mu(d\omega_0,d{\bf y}) = P(d\omega_0) \, \prod_{k=1}^N Q^k(dy^k)$. These relaxations can be classified in increasing order as \textit{quantum-correlated relaxation} and \textit{non-signaling relaxation}. These classes are new to team decision theory and we believe that these new classes provide novel perspectives and results to team decision theory and decentralized stochastic control in the future. 

We saw earlier in the article that individual or common (independent) randomness do not improve optimal team cost, whereas we will see that quantum-correlated and non-signaling relaxations in general improve optimal cost. Moreover, the optimization problem associated with non-signaling case can be written as a linear program. 

A related hierarchy of policies was introduced in \cite{AuFeRaScWi16} to study games with local information. In \cite{AuFeRaScWi16}, advantages of quantum-correlated and non-signaling equilibria over classical ones were discussed and it was established that quantum-correlated and non-signaling equilibria are socially more beneficial. However, our aim is not to show benefits of quantum-correlated and non-signaling policies over classical ones, instead our main motivation to obtain an approximation to the classical team problem by using these new classes of policies. 

\sy{The non-signaling type relaxation for team problems was studied in \cite{YukselSaldiSICON17} with a different name, building on \cite{anantharam2007common} where it was proved that the set of strategic measures corresponding to the extreme points of non-signaling policies is a strict superset of the set of strategic measures corresponding to measurable policies. This set arises in applications in information theory in the context of converse theorems in multi-terminal source coding (see e.g. Berger-Tung inner-outer bounds \cite{wagner2008improved}). In quantum information theory literature \cite{brunner2014bell}, this correlation structure has evidently also been studied extensively, as we will present later in this section.}

\subsection{Non-signaling Relaxation}\label{nonsignal}

A joint conditional distribution of actions given observations is \textit{non-signaling} if, for any $i \in {\cal N}$, the marginal distribution of action of DM~$i$ given its observation does not give any information about the observations of other agents \cite{MaAcGi06,BaLiMaPiPoRo05}. Non-signaling has been investigated in quantum mechanics due to its close connection to foundations of quantum mechanics and relativity \cite{PoRo94}.  

Let {\bf $L_{NS}(\mu)$} denote the set of non-signaling strategic measures. Formally, non-signaling strategic measures are defined as follows. An element $P \in L_{NS}(\mu)$ if $P(d\omega_0, d{\bf y}) = \mu(d\omega_0, d{\bf y})$, and for any subset $\{k_1,\ldots,k_M\}$ of $\{1,\ldots,N\}$, the actions of the agents in $\{k_1,\ldots,k_M\}$ given their measurements are independent of measurements of agents in $\{1,\ldots,N\}\setminus\{k_1,\ldots,k_M\}$ and $\omega_0$; that is, for any $\{k_1,\ldots,k_M\} \subset \{1,\ldots,N\}$, we have 
\begin{align}
P(du^{k_1},\ldots,du^{k_M}\, | \, \omega_0,{\bf y}) = P(du^{k_1},\ldots,du^{k_M} \, | \, y^{k_1},\ldots,y^{k_M}). \label{nonsignaling}
\end{align}

We also define a more relaxed version of $L_{NS}(\mu)$, denoted by {\bf $L_M(\mu)$} (which we will term as local-Markov correlations), as follows: An element $P \in L_{M}(\mu)$ if $P(d\omega_0, d{\bf y}) = \mu(d\omega_0, d{\bf y})$, and for any $k \in \{1,\ldots,N\}$, the action $u^k$ of agent $k$ given the measurement $y^k$ is independent of measurements of agents in $\{1,\ldots,N\}\setminus\{k\}$ and $\omega_0$; that is, we have 
\begin{align}
P(du^{k} | \, \omega_0,{\bf y}) = P(du^{k}  | y^k).\label{LM}
\end{align}


At first sight, it is tempting to claim that $L_C(\mu)$ is the same as the set of non-signaling strategic measures $L_{NS}(\mu)$ or $L_M(\mu)$. Indeed, in \cite{anantharam2007common} this question was raised: a counterexample to the claim that $L_C(\mu)$ is equivalent to the set $L_{NS}(\mu)$ was given to establish that the set of extreme points of non-signaling policies is not $L_A(\mu)$, which would imply that non-signaling policies are more general than randomized ones as $L_A(\mu)$ is the set of extreme points of $L_C(\mu)$ (see Theorem~\ref{FeinbergStrategic}). In the quantum information literature, this result was known long ago, as we will review later.

It turns out that the non-signaling condition (\ref{nonsignaling}) can be derived from fewer linear constraints, which will be described below. These constraints indeed enable us to write the optimization problem associated with non-signaling policies as a linear program. 

\begin{lemma}[{\cite[Section II-A]{MaAcGi06}}]\label{equiv}
An element $P$ with $P(d\omega_0, d{\bf y}) = \mu(d\omega_0, d{\bf y})$ is a non-signaling strategic measure if and only if it satisfies the following condition: 

(\textbf{N}) For each $k \in \{1,\ldots,N\}$, the marginal distribution of actions excluding  $u^k$ is independent of the observation $y^k$ and $\omega_0$:
\begin{align}
P(d{\bf u}^{-k}|\omega_0,{\bf y}) = P(d{\bf u}^{-k}|{\bf y}^{-k}).\label{cons}
\end{align}
\end{lemma}

\begin{proof}
Note that each constraint (\ref{cons}) is linear in $P$. The equivalence of (\ref{nonsignaling}) and (\textbf{N}) can be established as follows. First, it is immediate that (\ref{nonsignaling}) implies (\textbf{N}). Conversely, let $P$ satisfies (\textbf{N}). Fix any subset $\{k_1,\ldots,k_M\}$ of $\{1,\ldots,N\}$. Let $\{l_1,\ldots,l_T\} \coloneqq \{1,\ldots,N\}\setminus\{k_1,\ldots,k_M\}$. Then, using the condition (\textbf{N}), we can prove that 
\begin{align}
P(du^{k_1},\ldots,du^{k_M}\, | \, \omega_0,{\bf y})
= P(du^{k_1},\ldots,du^{k_M}\, | \, \hat{\omega}_0,\hat{y}^{l_1},\ldots,\hat{y}^{l_T},y^{k_1},\ldots,y^{k_M}), \nonumber 
\end{align} 
for all $\omega_0,y^{l_1},\ldots,y^{l_T}$ and $\hat{\omega}_0,\hat{y}^{l_1},\ldots,\hat{y}^{l_T}$, which implies that 
\begin{align}
P(du^{k_1},\ldots,du^{k_M}\, | \, \omega_0,{\bf y}) = P(du^{k_1},\ldots,du^{k_M} \, | \, y^{k_1},\ldots,y^{k_M}). \nonumber 
\end{align} 
Indeed, we have
\begin{align}
P(du^{k_1},\ldots,du^{k_M}\, | \, \omega_0,{\bf y})
&= \int_{\mathbb{U}^{l_2}\times\ldots\times\mathbb{U}^{l_T}} P(d{\bf u}^{-l_1} \, | \, \omega_0,{\bf y})  \nonumber \\
&= \int_{\mathbb{U}^{l_2}\times\ldots\times\mathbb{U}^{l_T}} P(d{\bf u}^{-l_1} \, | \, \hat{\omega}_0,\hat{y}^{l_1}, {\bf y}^{-l_1})  \,\, \text{(by (\textbf{N}))} \nonumber \\
&= P(du^{k_1},\ldots,du^{k_M}\, | \, \hat{\omega}_0,\hat{y}^{l_1}, {\bf y}^{-l_1}) \nonumber \\
&= \int_{\mathbb{U}^{l_1}\times\mathbb{U}^{l_3}\times\ldots\times\mathbb{U}^{l_T}} P(d{\bf u}^{-l_2} \, | \, \hat{\omega}_0,\hat{y}^{l_1}, {\bf y}^{-l_1})  \nonumber \\
&= \int_{\mathbb{U}^{l_1}\times\mathbb{U}^{l_3}\times\ldots\times\mathbb{U}^{l_T}} P(d{\bf u}^{-l_2} \, | \, \hat{\omega}_0, \hat{y}^{l_1},\hat{y}^{l_2}, {\bf y}^{-l_1,-l_2})  \,\, \text{(by (\textbf{N}))} \nonumber \\
&= P(du^{k_1},\ldots,du^{k_M} \, | \, \hat{\omega}_0, \hat{y}^{l_1},\hat{y}^{l_2}, {\bf y}^{-l_1,-l_2}) \nonumber \\
&\phantom{x}\vdots \nonumber \\
&= P(du^{k_1},\ldots,du^{k_M}\, | \, \hat{\omega}_0,\hat{y}^{l_1},\ldots,\hat{y}^{l_T}, y^{k_1},\ldots,y^{k_M}). \nonumber 
\end{align}
This completes the proof. 
\end{proof}

The following result establishes convexity and topological properties of the sets $L_{NS}(\mu)$ and $L_M(\mu)$. 
\smallskip

\begin{theorem}\label{nonsignal-top}
\begin{itemize}
\item[(i)] $L_M(\mu)$ is a convex set and it is closed under the weak convergence topology. 
\item[(ii)] $L_{NS}(\mu)$ is a convex set and it is closed under the weak convergence topology.
\end{itemize}
\end{theorem}

\smallskip

\begin{proof} 
Convexity of both sets follow since both $L_M(\mu)$ and $L_{NS}(\mu)$ are intersections of convex sets.

Let $P_n \in L_{M}(\mu)$ converge to $P$ weakly. Then, for any continuous and bounded $g$:
\[ \int P_n(dy^i,d{\bf y}^{-i}, d\omega_0,u^i) g(y^i,{\bf y}^{-i},\omega_0,u^i) \to \int P(dy^i,d{\bf y}^{-i}, d\omega_0,u^i) g(y^i,{\bf y}^{-i},\omega_0,u^i)\]
Since $P_n \in L_{M}(\mu)$, we have that
\begin{align}
&\int P_n(dy^i,d{\bf y}^{-i}, d\omega_0,u^i) g(y^i,{\bf y}^{-i},\omega_0,u^i) \nonumber \\
&\phantom{xxxxxxxxxxxxxx}= \int \int P_n(d{\bf y}^{-i},d\omega_0 | y^i) P_n(dy^i, du^i) g(y^i,{\bf y}^{-i},\omega_0,u^i). \nonumber 
\end{align}
But since the information structure is static, we write
\begin{align} 
&\int P_n(dy^i,d{\bf y}^{-i}, d\omega_0,u^i) g(y^i,{\bf y}^{-i},\omega_0,u^i) \nonumber \\
&\phantom{xxxxxxxxxxxxxx}= \int \int \mu(d{\bf y}^{-i},d\omega_0 | y^i) P_n(dy^i, du^i) g(y^i,{\bf y}^{-i},\omega_0,u^i). \nonumber
\end{align}
Now, 
\[h(y^i,u^i) := \bigg( \int \mu(d{\bf y}^{-i},d\omega_0 | y^i) g(y^i,{\bf y}^{-i},\omega_0,u^i) \bigg) \]
is continuous in $u^i$ as $g$ is continuous by dominated convergence theorem. Since the marginal on $y^i$ is fixed, as noted earlier, \cite[Theorem 3.10]{Schal} (or  \cite[Theorem 2.5]{balder2001}) establishes that the set of all tight probability measures with a fixed marginal on $\prod_k \mathbb{Y}^k$ is relatively compact under the $w$-$s$ topology. This implies that $P_n$ converges to $P$ in $w$-$s$ topology as the sequence $\{P_n\}$ is tight. Hence,
\[ \int P_n(dy^i,d{\bf y}^{-i}, d\omega_0,u^i) g(y^i,{\bf y}^{-i},\omega_0,u^i) \to \int h(y^i,u^i) P(dy^i,du^i)   \]
\[ =  \int \int \mu(d{\bf y}^{-i},d\omega_0 | y^i) P(dy^i, du^i) g(y^i,{\bf y}^{-i},\omega_0,u^i)\]
and hence $P \in L_{M}(\mu)$ also.

By viewing each subset $\{k_1,\ldots,k_m\}$ of $\{1,\ldots, N\}$ as a single decision maker with a collective action ${\bf u}_m:=\{u^{k_1},\ldots,u^{k_m}\}$ and ${\bf y}_m:=\{y^{k_1},\ldots,y^{k_m}\}$, applying the same analysis above shows that each set is closed, and since intersection of closed sets is closed, the closedness of $L_{NS}(\mu)$ follows. As noted above using the condition (\textbf{N}), without any loss, we can only consider the sets of cardinality $N-1$.
\end{proof}

Since the constraints $P(d{\bf u}^{-i} | y^i, {\bf y}^{-i}, \omega_0) = P(d{\bf u}^{-i}|{\bf y}^{-i}), \, i \in {\cal N}$, in $L_{NS}(\mu)$ for $P$ are linear, the optimization problem associated with non-signaling strategic measures can be written as a linear program over appropriate vector spaces. One such formulation will be given in Section~\ref{dual}. 

\subsection{Quantum correlations}\label{quantum}

To introduce quantum-correlated strategic measures, we briefly introduce the mathematical formalism necessary to discuss quantum operations. We refer the reader to books \cite{NiCh02,Wat18,Hol11} for basics of quantum information and computation. In this section, we suppose that all Hilbert spaces are complex and separable.

\def\Tr{\mathop{\rm Tr}}
\def\Pos{\mathop{\rm Pos}}
\def\Id{\mathop{\rm Id}}
\def\D{{\mathcal D}}
\def\H{{\mathcal H}}
\def\L{{\mathcal L}}

Quantum physical systems are described by complex Hilbert spaces. For a Hilbert space $\H$ with an inner product $\langle \cdot , \cdot \rangle$, let $\L(\H)$ and $\D(\H)$ denote the set of Hermitian and positive Hermitian operators with unit trace. A state $\rho$ of a quantum system, living in $\H$, is an element of $\D(\H)$; that is, it is a positive Hermitian operator with unit trace. A quantum state $\rho$ is said to be \textit{pure} if it has rank equal to $1$. Equivalently, $\rho$ is pure state if there exists a unit vector $u \in \H$ such that 
$$
\rho = u \, u^*,
$$
where $u^*$ is the complex conjugate of the vector $u$. By spectral decomposition, every quantum state can be written as a mixture of pure states. Hence, the extreme points of $\D(\H)$ coincides with the set of pure states.

A measurement on a quantum system is given by Positive Operator Valued Measure (POVM) \cite{Hol11}; that is, given any measurable outcome space $(\mathbb{X},{\cal X})$, a POVM is a mapping $M: {\cal X} \rightarrow \L(\H)$ from the $\sigma$-algebra ${\cal X}$ of $\mathbb{X}$ to the Hermitian operators $\L(\H)$ on $\H$ such that
\begin{itemize}
\item[(1)] $M(\emptyset) = 0, \,\, M(\mathbb{X}) = \Id$;
\item[(2)] $M(B) \geq 0 \,\, \forall B \in {\cal X}$;
\item[(3)] $M(\bigcup_j B_j) = \sum_{j} M(B_j)$ for any disjoint collection $\{B_j\},$ where the series is weakly convergent. 
\end{itemize} 
Here, weak convergence means the following: 
$$
\lim_{n \rightarrow \infty} \left \langle \psi, \sum_{j=1}^n M(B_j) \, \psi \right \rangle = \left \langle \psi, \sum_{j=1}^{\infty} M(B_j) \, \psi \right \rangle,
$$ 
for all $\psi \in {\cal H}$.

If the measurable outcome space $(\mathbb{X},{\cal X})$ is finite in POVM $M: {\cal X} \rightarrow {\cal L}(\H)$, it is sufficient to define POVM $M: \mathbb{X} \rightarrow {\cal L}(\H)$ on $\mathbb{X}$ and extend the definition to ${\cal X}$ as follows: for any $A \in {\cal X}$, define $M(A) = \sum_{x \in A} M(x)$. For instance, if $\mathbb{X} =\{0,1\}$, then we will only need to define $M(0)$ and $M(1)$. In finite case, $M: \mathbb{X} \rightarrow {\cal L}(\H)$ should satisfy the following conditions: $M(x) \geq 0 \,\, \forall x \in \mathbb{X}$ and 
$\sum_{x} M(x) = \Id$.

When one applies the measurement $M: {\cal X} \rightarrow \L(\H)$ to the quantum system in the state $\rho$, the probability of obtaining the event $B \in {\cal X}$ is given by 
$$P(B \, | \, \rho) =  \Tr\bigl\{M(B) \, \rho \bigr\}.$$
Note that since $\Tr\{\rho\}=1$ and $M(B) \geq 0$ for all $B \in {\cal X}$, we have $P(B|\rho) \geq 0$ for all $B \in {\cal X}$ and $P(\mathbb{X}|\rho) = 1$. Moreover, for any disjoint collection $\{B_j\}$ of events, by (3) we have 
$$
P(B \, | \, \rho) = \sum_j P(B_j \, | \, \rho),
$$
where $B \coloneqq \bigcup_j B_j$. Hence, $P(\cdot|\rho)$ is indeed a probability measure on $\mathbb{X}$. We note that if $\rho$ is a pure state given by a unit vector $u$, then
$$P(B \, | \, \rho) =  \langle u, M(A) \, u \rangle.$$

In quantum physics, a compound of $N$ quantum systems with the underlying Hilbert spaces $\{\H_i, i=1,\ldots,N\}$  is represented by the tensor product $\H_1 \otimes \H_2 \otimes \ldots \otimes \H_N$ of the individual Hilbert spaces. Therefore, any state $\rho$ (called the \textit{compound state}) of this compound quantum system is an element of $\D(\H_1 \otimes \H_2 \ldots \otimes \H_N)$. If the state $\rho$ is in product form 
$$
\rho = \rho_1 \otimes \ldots \otimes \rho_N,
$$
then it is called \textit{product state}. In this case, the individual quantum systems is said to be independent. If $\rho$ is a mixture of product states, it is called \textit{separable}. If the state is not separable, then it is called \textit{entangled}.

{\bf Quantum-correlated strategic measures.} With these definitions, we can now define quantum-correlated strategic measures. An element $P$ 
is a \textit{quantum-correlated} strategic measure if agents have access to a part of a quantum compound state $\rho \in \D(\H_1 \otimes \H_2 \ldots \otimes \H_N)$, where $\{\H_i, i=1,\ldots,N\}$ is a collection of arbitrary Hilbert spaces, and, for each $i=1,\ldots,N$, DM~$i$ makes measurements $M^{\sy{i},y^i}: \mathbb{U}^i \rightarrow \L(\H_i)$ on $i^{th}$ part of the state $\rho$ depending on its observations $y^i$ to generate its action $u^i$ as the output of the measurement (recall that ${\cal U}^i$ is the Borel $\sigma$-algebra on the action space $\mathbb{U}^i$ of DM~$i$); that is, the conditional distribution $P(d{\bf u}|\omega_0,{\bf y})$ is of the following form:
\begin{align}
P(d{\bf u}|\omega_0,{\bf y}) = \Tr\left\{\left( M^{1,y^1}(du^1) \otimes \ldots \otimes M^{N,y^N}(du^N) \right) \rho \right\}. \nonumber
\end{align}

Let the state $\rho$ be separable; that is,
$$
\rho = \int \nu(dz) \, \rho_1^z \otimes \ldots \otimes \rho_N^z. 
$$ 
Suppose that on each individual quantum system, we are performing the following measurement $M^{i,y^i}: \mathbb{U}^i \rightarrow \L(\H_i)$. Then, the corresponding joint probability of obtaining the event $du^1\times\ldots\times du^N$ given observations $(y^1,\ldots,y^N)$ and state $\omega_0$ is the following:
\begin{align}
P(du^1,\ldots,du^N|\omega_0,{\bf y}) &= \int \nu(dz) \, \Tr\{M^{1,y^1}(du^1) \rho_1^z\} \, \ldots \, \Tr\{M^{N,y^N}(du^N) \rho_N^z\} \nonumber \\
&= \int \nu(dz) \, P(du^1|y^1,z) \, \ldots \, P(du^N|y^N,z), \nonumber 
\end{align}
where $P(du^i|y^i,z) = \Tr\{M^{i,y^i}(du^i) \rho_i^z\}$. Therefore, $P \in L_C(\nu)$. Hence, separable states can only generate correlations that can be classically generated; that is, if individual quantum systems share a common randomness $\nu(dz)$, then they can also realize $P(du^1,\ldots,du^N|\omega_0,{\bf y})$  independently using classical devices given the common randomness. Hence, to generate correlations that cannot be realized classically, one must use entangled states.

Let $L_Q(\mu)$ denote the set of quantum-correlated strategic measures. In the next section, we will show that $L_C(\mu)$ is a strict subset of $L_Q(\mu)$. In  view of discussion in the previous paragraph, the strategic measure $P \in L_Q(\mu) \setminus L_C(\mu)$ if the shared state $\rho$ is entangled.

As an example, consider the two-agent team problem with $\mathbb{Y}^1 = \mathbb{Y}^2 = \mathbb{U}^1 =\mathbb{U}^2 = \{0,1\}$. Suppose that $\Omega_0$ is degenerate. Let $\H_1 = \H_2 = \mathbb{C}^2$. Consider the following states on $\mathbb{C}^2 \times \mathbb{C}^2$:
\begin{align}
\sigma &= \frac{1}{2} E_{1,1} \otimes E_{1,1} + \frac{1}{2} E_{2,2} \otimes E_{2,2} \nonumber \\
\tau &= \frac{1}{2} E_{1,1} \otimes E_{1,1} + \frac{1}{2} E_{1,2} \otimes E_{1,2} + \frac{1}{2} E_{2,1} \otimes E_{2,1} + \frac{1}{2} E_{2,2} \otimes E_{2,2}, \nonumber 
\end{align}
where $E_{a,b}$ is a matrix on $\mathbb{C}^2$ such that $E_{a,b}(i,j) = 1_{\{(i,j)=(a,b)\}}$, for $a,b \in \{1,2\}$. Note that $E_{1,1}, E_{2,2} \in \D(\mathbb{C}^2)$. Hence, $\rho$ is a mixture of product states, and so, it is separable. However, note that $E_{1,2}, E_{2,1} \notin \D(\mathbb{C}^2)$. Therefore, $\tau$ is an entangled state (indeed it is maximally entangled state). For any angle $\theta$, let us define 
$$
\Pi_{\theta} = \begin{pmatrix}
\cos^2(\theta) & \cos(\theta) \sin(\theta) \\
\cos(\theta)\sin(\theta) & \sin^2(\theta)
\end{pmatrix} = \begin{pmatrix} \cos(\theta) \, \sin(\theta) \end{pmatrix} \, \begin{pmatrix} \cos(\theta) \\ \sin(\theta) \end{pmatrix}.
$$
Using $\Pi_{\theta}$, we now define two collection of measurements $\{M^{1,y^1}: \mathbb{U}^1 \rightarrow {\cal L}(\mathbb{C}^2); y^1 \in \mathbb{Y}^1\}$, $\{M^{2,y^2}: \mathbb{U}^2 \rightarrow {\cal L}(\mathbb{C}^2); y^2 \in \mathbb{Y}^2\}$ for DM~$1$ and DM~$2$, respectively, as follows: 
$$
M^{1,0}(0)=\Pi_0, M^{1,0}(1)=\Pi_{\pi/2}; \,\, M^{1,1}(0)=\Pi_{\pi/4}, M^{1,1}(1)=\Pi_{3\pi/4}
$$ 
and 
$$
M^{2,0}(0)=\Pi_{\pi/8}, M^{2,0}(1)=\Pi_{5\pi/8}; \,\, M^{2,1}(0)=\Pi_{7\pi/8}, M^{2,1}(1)=\Pi_{3\pi/8}.
$$ 
Namely, if the observation of DM~$1$ is $y^1$, then DM~$1$ applies the measurement $M^{1,y^1}: u^1  \mapsto M^{1,y^1}(u^1)$ to the quantum state. Similarly, if the observation of DM~$2$ is $y^2$, then DM~$2$ applies the measurement $M^{2,y^2}: u^2  \mapsto M^{2,y^2}(u^2)$ to the quantum state.

Now, if we apply these measurements to the quantum states $\sigma$ and $\tau$, we obtain two quantum-correlated strategic measures $P_{\sigma}$ and $P_{\tau}$, respectively. The conditional distributions of these quantum-correlated strategic measures on actions given observations can be computed as follows:
\small
\begin{align}
&P_{\sigma}(u^1,u^2|y^1,y^2) \nonumber \\
&= \frac{1}{2} \left( \Tr\{M^{1,y^1}(u^1) E_{1,1}\} \Tr\{M^{2,y^2}(u^2) E_{1,1}\} + \Tr\{M^{1,y^1}(u^1) E_{2,2}\} \Tr\{M^{2,y^2}(u^2) E_{2,2}\} \right) \nonumber \\
&P_{\tau}(u^1,u^2|y^1,y^2) \nonumber \\
&= \frac{1}{2} \left( \Tr\{M^{1,y^1}(u^1) E_{1,1}\} \Tr\{M^{2,y^2}(u^2) E_{1,1}\} + \Tr\{M^{1,y^1}(u^1) E_{1,2}\} \Tr\{M^{2,y^2}(u^2) E_{1,2}\} \right)\nonumber \\
&\phantom{xx}+ \frac{1}{2} \left( \Tr\{M^{1,y^1}(u^1) E_{2,1}\} \Tr\{M^{2,y^2}(u^2) E_{2,1}\} + \Tr\{M^{1,y^1}(u^1) E_{2,2}\} \Tr\{M^{2,y^2}(u^2) E_{2,2}\} \right). \nonumber 
\end{align}
\normalsize
Since $E_{1,1}, E_{2,2} \in \D(\mathbb{C}^2)$, the conditional distribution $P_{\sigma}(u^1,u^2|y^1,y^2)$ is a mixture of the following independent conditional probability measures 
$$P_{1,\sigma}(u^1|y^1) P_{1,\sigma}(u^2|y^2) \coloneqq \Tr\{M^{1,y^1}(u^1) E_{1,1}\} \Tr\{M^{2,y^2}(u^2) E_{1,1}\}$$ and $$P_{2,\sigma}(u^1|y^1) P_{2,\sigma}(u^2|y^2) \coloneqq \Tr\{M^{1,y^1}(u^1) E_{2,2}\} \Tr\{M^{2,y^2}(u^2) E_{2,2}\}.$$ 
Hence $P_{\sigma}$ can be classically realized via common randomness; that is, $P_{\sigma} \in L_C(\mu)$. However, this is not true for $P_{\tau}$ since $E_{1,2}, E_{2,1} \notin \D(\mathbb{C}^2)$, and so, 
$$\Tr\{M^{1,y^1}(u^1) E_{1,2}\} \Tr\{M^{2,y^2}(u^2) E_{1,2}\}$$ and $$\Tr\{M^{1,y^1}(u^1) E_{2,1}\} \Tr\{M^{2,y^2}(u^2) E_{2,1}\}$$
are not probability measures. Therefore, $P_{\tau}$ cannot be realized classically; that is, $P_{\tau} \in L_Q(\mu) \setminus L_C(\mu)$.

Here, one can view the following collection as the joint strategy of the decision makers:
$$
\left\{\H_i \, (i \in {\cal N}), \rho \in  \D(\H_1 \otimes \H_2 \ldots \otimes \H_N), \{M^{i,y^i}: {\cal U}^i \rightarrow \L(\H_i); y^i \in \mathbb{Y}^i\} \, (i \in {\cal N}) \right \}. 
$$
Namely, at the beginning of the problem, decision makers agree on the Hilbert spaces $\H_i  \, (i \in {\cal N})$ and share a compound state $\rho \in  \D(\H_1 \otimes \H_2 \ldots \otimes \H_N)$ (i.e., entangled state) to be used in the measurements. Then, they choose decentralized measurements $\{M^{i,y^i}: {\cal U}^i \rightarrow \L(\H_i); y^i \in \mathbb{Y}^i\} \, (i \in {\cal N})$ without communicating with each other.

Given the Hilbert spaces and the compound state, it is possible to view decentralized measurements as the decentralized strategies of the decision makers. However, when designing optimal quantum-correlated strategic measures, decision makers must also find the optimal Hilbert spaces and the corresponding optimal compound state. Note that maximally entangled state is in general  not the optimal choice \cite{BrGiSc05}. Therefore, the optimization variable in the decision problem corresponding to the quantum-correlated strategic measures is in general quite large even if the measurement and the action spaces are finite. However, note that \cite{ClHoToWa04,KeReTo10} (see also Example~\ref{XORteam}) give evidence that quantum-correlated teams with finite measurement and action spaces can be computationally tractable as opposed to their classical randomized counterparts even though the optimization space is quite large compared to the classical case. Namely, in these papers, optimization problems associated with quantum-correlated strategic measures can be written as (or can be approximated by) semi-definite programs whose sizes scale with the cardinality of the measurement and action spaces. As a result, they can be solved exactly or approximately in polynomial time. In particular, \cite{ClHoToWa04} computes the optimal value of XOR team and \cite{KeReTo10} approximates the optimal value of the unique teams via semi-definite programs. 

However, there are other instances of finite team problems \cite{KeKoMaToVi11,Vid16,NaVi18,ji2020mip}, where exact or approximate computation of the optimal cost of quantum-correlated teams is NP-hard, and therefore, cannot be cast as a semi-definite program. One of the reasons for such negative complexity results might be the fact that the optimal value of quantum-correlated teams can be attained via quantum systems living in infinite-dimensional Hilbert spaces. Indeed, there are some evidences in the literature that supports this observation (see, e.g., \cite{PaVe10}). Therefore, in order to obtain a computationally tractable problem, we may need to put a constraint on the dimensions of the Hilbert spaces where quantum systems live in.

To that end, let $L_{Q(d)}(\mu) \subset L_Q(\mu)$ denote the set of quantum-correlated strategic measures where each agent has access to a part of quantum compound state $\rho \in \D(\H_1 \otimes \H_2 \ldots \otimes \H_N)$ with $\dim(\H_i) \leq d$ for all $i \in {\cal N}$. 

In the following, we state convexity and topological properties of the sets $L_Q(\mu)$ and $L_{Q(d)}(\mu)$. 

\smallskip

\begin{theorem}\label{quantum-convex-closed}
\begin{itemize}
\item[(i)] $L_Q(\mu)$ is a convex set but it is not closed under the weak convergence topology. 
\item[(ii)] $L_{Q(d)}(\mu)$ is not convex if $L_C(\mu) \not \subset L_{Q(d)}(\mu)$. In particular, $L_{Q(1)}(\mu)$ is not convex and $L_{Q(2)}(\mu)$ is not convex for certain models. 
\item[(iii)]For each $d$, $L_{Q(d)}(\mu)$ is closed under the weak convergence topology if the observation spaces are finite and action spaces are compact. 
\item[(iv)] $L_Q(\mu)$ is a strict super-set of $\bigcup_d L_{Q(d)}(\mu)$. 
\end{itemize}
\end{theorem}

\begin{proof}
Proof of (i): The convexity of $L_Q(\mu)$ can be proved easily and so we omit the details. Non-closedness of $L_Q(\mu)$ has been a longstanding open problem and is proved  recently in \cite{slofstra2019set} who showed that the weak closure of the set is strictly larger than the set itself by explicit construction.

Proof of (ii): By \cite[Corollary 2]{MaWo15}, $L_{Q(d)}(\mu)$ is not convex if $L_C(\mu) \not \subset L_{Q(d)}(\mu)$. In particular, $L_{Q(1)}(\mu)$ is not convex by \cite[Proposition 1]{MaWo15} and $L_{Q(2)}(\mu)$ is not convex for certain models by \cite[Proposition 2]{MaWo15}.

Proof of (iii):  Since both the set of states and the set of POVMs of fixed (local) dimension are compact in finite measurement space and compact action space case and the trace is a continuous map, $L_{Q(d)}(\mu)$ is closed (see \cite[Appendix B]{PhysRevA.97.022104} \syr{where the finite action case is considered but the same argument applies for the compact case also}). 

Proof of (iv): The fact that $L_Q(\mu)$ is a strict super-set of $\bigcup_d L_{Q(d)}(\mu)$ is proved in \cite{CoSt20}; that is, authors show that there exists a joint strategy which is not attainable by quantum systems of any arbitrary finite dimension, but is attained
exclusively by infinite-dimensional quantum systems
\end{proof}


The following result is a corollary of Theorem~\ref{quantum-convex-closed}-(iii).

\begin{corollary}
Consider a static team decision problem with finite measurement spaces (and hence independent measurements reduction) and compact action spaces. Then, the following team optimal problem $\inf_{P \in L_{Q(d)}(\mu)} J(P)$ has an optimal team policy.
\end{corollary}

To show the effectiveness (both computationally and optimally) of quantum-correlated strategic measures over classical randomized ones, we consider the XOR team problem and a variation below.  

\begin{example}\label{XORteam}
In XOR team, we have two agents with binary action spaces $\{0,1\}$; that is, $\mathbb{U}^1=\mathbb{U}^2=\{0,1\}$. Measurements are generated independently and uniformly over some finite sets $\mathbb{Y}^1$ and $\mathbb{Y}^2$ via probability measure $\mu(y^1,y^2)$. Hence, there is no $\omega_0$ variable in the problem, and so the problem is automatically static. The reward\footnote{All results in this paper are also true for maximization of a reward function.} function is defined as
\begin{align}
r(y^1,y^2,u^1,u^2) =
\begin{cases}
1, & \text{if $u^1 \oplus u^2 = h(y^1,y^2)$} \\
-1, & \text{otherwise},
\end{cases}\nonumber
\end{align}
where $h: \mathbb{Y}^1 \times \mathbb{Y}^2 \rightarrow \{0,1\}$ is some arbitrary binary-valued function. This team problem with quantum-correlated policies can be written as a semi-definite program due to Tsirelson's Theorem \cite[Theorem 6.62]{Wat18}.

Indeed, let us define 
$$g(y^1,g^2) = \mu(y^1,y^2) \, (-1)^{h(y^1,y^2)}.$$ 
Given any $\rho \in \D(\H_1 \otimes \H_2)$ for some finite-dimensional Hilbert spaces $\H_1,\H_2$ (finite dimensional Hilbert spaces are sufficient for this problem) and given any two collection of POVMs $\{M^{1,y^1}, y^1 \in \mathbb{Y}^1\}$, $\{M^{2,y^2}, y^2 \in \mathbb{Y}^2\}$, the corresponding strategic measure $P$ is $P(y^1,y^2) = \mu(y^1,y^2)$ and 
\begin{align}
&P(u^1,u^2|y^1,y^2) = \Tr\left\{\left( M^{1,y^1}(u^1) \otimes  M^{2,y^2}(u^2) \right) \rho \right\} \nonumber \\
\intertext{and its expected reward function can be written as}
&J(P) = \sum_{y^1,y^2}  g(y^1,y^2) \, \Tr\bigl\{(M^{1,y^1}(0) - M^{1,y^1}(1))  \otimes  (M^{2,y^2}(0) - M^{2,y^2}(1)) \rho \bigr\}. \nonumber 
\end{align}
Note that an operator $H$ is Hermitian with operator norm $\|H\| \leq 1$ if and only if it can be written as $H = M(0) - M(1)$, where $M(i)$ ($i=0,1$) are positive semi-definite operators with $M(0) + M(1) = \Id$. Therefore, for any pair $(y^1,y^2)$, $H^{y^1} = M^{1,y^1}(0) - M^{1,y^1}(1)$ and $H^{y^2} = M^{2,y^2}(0) - M^{2,y^2}(1)$ are Hermitian with operator norms less than $1$. Conversely, for any pair $(y^1,y^2)$, any Hermitian operators $H^{y^1},H^{y^2}$ with operator norms less than $1$ can be decomposed as above. 

Let $X$ be a $|\mathbb{Y}^1| \times |\mathbb{Y}^2|$ real matrix. Tsirelson's Theorem states that the following assertions are equivalent: 
\begin{itemize}
\item[1.] There exist Hilbert spaces $\H_1$ and $\H_2$, a state $\rho \in \D(\H_1 \otimes \H_2)$, and two collections of Hermitian operators
$$\{H^{y^1}, y^1 \in \mathbb{Y}^1\} \, \text{ and } \, \{H^{y^2}, y^2 \in \mathbb{Y}^2\}$$
whose operator norms are less than $1$, and 
$$X(y^1,y^2) = \Tr\left\{\left( H^{y^1} \otimes H^{y^2} \right) \rho\right\},$$
for all $y^1 \in \mathbb{Y}^1$, $y^2 \in \mathbb{Y}^2$.
\item[2.] There exists two collections $\{u_{y^1}, y^1 \in \mathbb{Y}^1\}$, $\{v_{y^2}, y^2 \in \mathbb{Y}^2 \} \subset \mathbb{R}^{|\mathbb{Y}^1| \times |\mathbb{Y}^2|}$ of unit vectors such that 
$$X(y^1,y^2) = u_{y^1}^T v_{y^2}$$ 
for all $y^1 \in \mathbb{Y}^1$, $y^2 \in \mathbb{Y}^2$.
\end{itemize}

Therefore, Tsirelson's Theorem and the above fact about Hermitian operators imply that 
\begin{align}
\sup_{P \in L_Q(\mu)} J(P) = \sup_{\substack{u_{y^1},v_{y^2} \in R^{|\mathbb{Y}^1| \times |\mathbb{Y}^2|} \\ y^1 \in \mathbb{Y}^1, y^2 \in \mathbb{Y}^2}} \sum_{y^1,y^2}   g(y^1,y^2) \, u_{y^1}^T v_{y^2}, \nonumber 
\end{align}
subject to 
\begin{align}
u_{y^1}^T u_{y^1} = 1 \, \text{ and } \, v_{y^2}^T v_{y^2} = 1, \, \text{ for } \, y^1 \in \mathbb{Y}^1,y^2\in \mathbb{Y}^2.\nonumber
\end{align}
This optimization problem is indeed a semi-definite program. Therefore, the optimal value of the XOR team with quantum-correlated policies can be computed in polynomial time as opposed to its classical counterpart. 
\end{example}

\begin{example}\label{CHSH}
A special case for XOR team is the celebrated CHSH (Clauser-Horne-Shimony-Holt) team \cite{Wat18}. In CHSH team, we have binary observation and action spaces $\{0,1\}$ and the reward function is defined as
\begin{align}
r(y^1,y^2,u^1,u^2) =
\begin{cases}
1, & \text{if $u^1 \oplus u^2 = y^1 \cdot y^2$} \\
-1, & \text{otherwise}.
\end{cases}\nonumber
\end{align}
For this problem, the optimal value of randomized policies is $0.5$ \cite[Section 6.3.2]{Wat18}, which can be verified by an inspection of 16 measurable strategic measuresas they are sufficient for optimality. However, quantum-correlated policies can achieve the maximum reward of $2 \sqrt{2}/4$, that is obtained by solving the corresponding semi-definite program. Therefore, for CHSH team, we have 
$$\sup_{P \in L_Q(\mu)} J(P) > J_*,$$ 
that is, quantum-correlated policies improve the optimal value $J^*$ of the original team as opposed to randomized policies. 

Indeed, Let $\H_1 = \H_2 = \mathbb{C}^{2}$ and let $\rho \in \D(\H_1 \otimes \H_2)$ be the maximally entangled state 
$$
\rho = \frac{1}{2} \, \sum_{a,b \in \{1,2\}} E_{a,b} \otimes E_{a,b}.
$$
Recall the two collection of POVMs $\{M^{1,y^1}, y^1 \in \mathbb{Y}^1\}$, $\{M^{2,y^2}, y^2 \in \mathbb{Y}^2\}$:
$$
M^{1,0}(0)=\Pi_0, M^{1,0}(1)=\Pi_{\pi/2}; \,\, M^{1,1}(0)=\Pi_{\pi/4}, M^{1,1}(1)=\Pi_{3\pi/4}
$$ 
and 
$$
M^{2,0}(0)=\Pi_{\pi/8}, M^{2,0}(1)=\Pi_{5\pi/8}; \,\, M^{2,1}(0)=\Pi_{7\pi/8}, M^{2,1}(1)=\Pi_{3\pi/8},
$$ 
where 
$$
\Pi_{\theta} = \begin{pmatrix}
\cos^2(\theta) & \cos(\theta) \sin(\theta) \\
\cos(\theta)\sin(\theta) & \sin^2(\theta)
\end{pmatrix}.
$$
Then, the corresponding quantum-correlated strategic measure $P$ has the following conditional distribution $P(u^1,u^2|y^1,y^2)$:
 
\begin{align}
P(\cdot,\cdot|0,0) &= \left(\frac{2+\sqrt{2}}{8}, \frac{2-\sqrt{2}}{8}, \frac{2+\sqrt{2}}{8}, \frac{2-\sqrt{2}}{8}\right), \nonumber \\
P(\cdot,\cdot|0,1) &= \left(\frac{2-\sqrt{2}}{8}, \frac{2+\sqrt{2}}{8}, \frac{2-\sqrt{2}}{8}, \frac{2+\sqrt{2}}{8}\right) \nonumber \\
P(\cdot,\cdot|1,0) &= \left(\frac{2+\sqrt{2}}{8}, \frac{2-\sqrt{2}}{8}, \frac{2-\sqrt{2}}{8}, \frac{2+\sqrt{2}}{8}\right) \nonumber \\
P(\cdot,\cdot|1,1) &= \left(\frac{2-\sqrt{2}}{8}, \frac{2+\sqrt{2}}{8}, \frac{2+\sqrt{2}}{8}, \frac{2-\sqrt{2}}{8}\right) \nonumber \\
\end{align}

Now, it is straightforward to compute that $J(P)=2 \sqrt{2}/4$. Hence, $P$ is the optimal quantum-correlated strategic measure. 
\end{example}

\subsection{Relations between convex relaxations}\label{relations}

Up to this point we have introduced a number of convex strategic measures, we now present a comparison. First, we note that the joint strategies in the set $L_C(\mu)$ is indeed corresponds to the \textit{local hidden variable} correlations in quantum mechanics. In a famous paper \cite{PhysRev.47.777}, Einstein, Podolsky, and Rosen discussed that due to the probabilistic nature of the quantum mechanics, this theory could not be a complete theory and should be supplemented by a local hidden variable to describe the probabilistic nature, where local hidden variable describes the uncertainties in experimental setup. In other words, they claimed that $L_C(\mu) = L_Q(\mu)$. However, John Bell in his famous paper \cite{PhysicsPhysiqueFizika.1.195} constructed a quantum mechanical system whose statistical behavior cannot be explained via local hidden variable theory; that is, $L_Q(\mu)$ is a strict super-set of $L_C(\mu)$. This result implies that quantum mechanical systems are intrinsically probabilistic and cannot be explained by the lack of precision in experimental measurements. Experimental verification of Bell's prediction has been established by Alain Aspect in \cite{PhysRevD.14.1944}. Additionally, in \cite{PoRo94}, Popescu and Rohrlich also discussed that correlations (strategic measures in our setting) achieved by quantum mechanical systems can be as large as $L_{NS}(\mu)$ which is larger than $L_Q(\mu)$. However, no experimental verification is provided for this claim yet. 

We also note that the lower bounds on the optimal team costs achieved by classical strategies $L_C(\mu)$ is known to be \textit{Bell inequalities} in quantum information theory literature. Bell inequality violations correspond to quantum-correlated strategies that improve the classical optimal team cost \cite{ClHoToWa04}. Identifying Bell inequality violations is a very active research area in quantum information theory. One such violation is established in Example~\ref{CHSH} for CHSH team. Since the considered sets are convex, by Hahn-Banach Separation Theorem \cite{lue69}, identifying Bell inequality violations is equivalent to finding a {\it separating} function in an inner-product form. Accordingly identifying Bell inequality violations would be implied by finding a cost function for which the optimal quantum-correlated team cost is strictly smaller than the optimal classical team cost; that is;
$$
\inf_{P \in L_Q(\mu)} \int P(ds) \, c(s) < \inf_{P \in L_C(\mu)} \int P(ds) \, c(s). 
$$
Since both $L_Q(\mu)$ and $L_C(\mu)$ are convex sets, above inequality states that the hyperplane $\{P: \int P(ds) \, c(s) = J^*\}$ separates the convex sets $L_Q(\mu)$ and $L_C(\mu)$. It is in general very hard to find all Bell inequality violations (i.e., separating hyperplanes) for a certain setup.

Before we establish the relation between convex relaxations, we introduce a further information structure, denoted by $L_{CJ}$, which allows for the common random variable $z$ to be correlated with $\omega_0, {\bf y}$; in this case, the information structure is essentially centralized, since the dependence between $z$ and the exogenous variables are arbitrary: 

\begin{align}\label{LJ}
&L_{CJ}(\mu) \coloneqq \bigg\{P \in {\cal P}\bigg(\Omega_0 \times \prod_{k=1}^N (\mathbb{Y}^k \times \mathbb{U}^k)\bigg): \nonumber \\
&\phantom{xxxxxxxxxxxx}P(B) = \int_{[0,1]^N \times B^0 \times \prod_{k=1}^N A^k} \mu(dz,d\omega_0, d{\bf y}) \prod_k 1_{\{u^k = \gamma^k(y^k,z) \in B^k\}},  \nonumber \\
&\phantom{xxxxxxxxxxxxxxxxxxxxxxxxxxxxxxxx} \gamma^k(\cdot,z) \in \Gamma^k, B \in {\cal B}(\Omega_0 \times \prod_k (\mathbb{Y}^k \times \mathbb{U}^k)) \bigg\}, \nonumber 
\end{align}
where the marginal of $\omega_0, {\bf y}$ is fixed as $\mu$.

\smallskip

\begin{theorem}\label{relation-thm}
We have the following relation between convex relaxations.
\begin{itemize}
\item[(i)] $L_C(\mu)  \subset L_{Q}(\mu) \subset L_{NS}(\mu) \subset L_{M}(\mu) \subset L_{CJ}(\mu).$
\item[(ii)] The sets $L_C(\mu), L_{Q}(\mu), L_{NS}(\mu), L_M(\mu), L_{CJ}(\mu)$ are convex.
\item[(iii)] The inclusions among the convex sets above are strict.
\item[(iv)] There exist problems for which
\[ \inf_{P \in L_Q(\mu)} \int P(ds) \, c(s) < \inf_{P \in L_{C}(\mu)} \int P(ds) \, c(s)\]
\[ \inf_{P \in L_{NS}(\mu)} \int P(ds) \, c(s) < \inf_{P \in L_{Q}(\mu)} \int P(ds) \, c(s)\]
and in particular
\[ \inf_{P \in L_{NS}(\mu)} \int P(ds) \, c(s) < \inf_{P \in L_{C}(\mu)} \int P(ds) \, c(s).\]
\end{itemize}
\end{theorem}

\noindent One implication of the above is that, the convex programs:
\[\inf_{P \in L_Q(\mu)} \int P(ds) \, c(s) \,\, \text{and} \,\,
\inf_{P \in L_{NS}(\mu)} \int P(ds) \, c(s) \]
or $\inf_{P \in L_{M}(\mu)} \int P(ds) \, c(s)$ provide lower bounds to the original problem. In particular, since the optimization problem associated with non-signaling strategic measures can be written as a linear program, it can be solved in a polynomial time if the measurement and action spaces are finite. 

\begin{proof}
Proof of (ii): We have already proved the convexity of $L_C(\mu)$, $L_{Q}(\mu)$, $L_{NS}(\mu)$, and $L_{M}(\mu)$. 

Proof of (i): Let $P \in L_R(\mu)$; that is,
$$P(d{\bf u}|\omega_0,{\bf y}) = \prod_{i=1}^N \Pi^i(du^i|y^i).$$
Let $\H$ be a Hilbert space with dimension $1$. Let $e_w \in \H$ with unit norm. We define 
\begin{align}
\rho &= e_w e_w^* \otimes \ldots \otimes e_w e_w^* \nonumber \\
M^{y^i}(A) &= \Pi^i(A|y^i) \, e_w e_w^* \nonumber 
\end{align}
for all $y^i \in \mathbb{Y}^i$, $A \in {\cal B}(\mathbb{U}^i)$, and $i=1,\ldots,N$. Then, we have
$$P(d{\bf u}|\omega_0,{\bf y}) = \Tr\bigl\{\left(M^{y^1}(du^1) \otimes \ldots \otimes M^{y^N}(du^N)\right) \rho \bigr\}.$$
This implies that $L_R(\mu) \subset L_Q(\mu)$. Since the convex hull of $L_R(\mu)$ is $L_C(\mu)$ and $L_Q(\mu)$ is convex, we have $L_C(\mu) \subset L_Q(\mu)$. 

Let $P \in L_Q(\mu)$; that is, agents have access to a part of a compound quantum state $\rho \in \D(\H_1 \otimes \H_2 \otimes \ldots \otimes \H_N)$, where $\{\H_i, i=1,\ldots,N\}$ is a collection of arbitrary Hilbert spaces, and, for each $i=1,\ldots,N$, DM~$i$ makes measurements $M^{y^i}: {\cal U}^i \rightarrow \L(\H_i)$ on the $i^{th}$ part of the state $\rho$ depending on its observations $y^i$ to generate its action as the output of the following measurement:
\begin{align}
P(d{\bf u}|\omega_0,{\bf y}) = \Tr\bigl\{\left( M^{y^1}(du^1) \otimes \ldots \otimes M^{y^N}(du^N)\right) \rho \bigr\}. \nonumber
\end{align}

We prove that $P$ satisfies the non-signaling condition. Fix any $k \in \{1,\ldots,N\}$. Then, we have 
\begin{align}
\int_{\mathbb{U}^k} P(d{\bf u} \,&| \,\omega_0, y^1,\ldots,y^k,\ldots,y^N) \nonumber \\
&= \int_{\mathbb{U}^k} \Tr\bigl\{\left( M^{y^1}(du^1) \otimes \ldots \otimes M^{y^k}(du^k) \otimes \ldots \otimes M^{y^N}(du^N) \right) \rho \bigr\} \nonumber \\
&= \Tr\biggl\{\biggl(\int_{\mathbb{U}^k} M^{y^1}(du^1) \otimes \ldots \otimes M^{y^k}(du^k) \otimes \ldots \otimes M^{y^N}(du^N)\biggr)\rho \biggr\}    \nonumber \\
&= \Tr\biggl\{\left( M^{y^1}(du^1) \otimes \ldots \otimes \Id \otimes \ldots \otimes M^{y^N}(du^N)\right)\rho \biggr\}    \nonumber \\
&= \int_{\mathbb{U}^k} \Tr\bigl\{\left(M^{y^1}(du^1) \otimes \ldots \otimes M^{\hat{y}^k}(du^k) \otimes \ldots \otimes M^{y^N}(du^N)\right) \rho \bigr\} \nonumber \\
&= \int_{\mathbb{U}^k} P(d{\bf u} \,| \, \omega_0,y^1,\ldots,\hat{y}^k,\ldots,y^N) \nonumber
\end{align}
for all $y^1,\ldots,y^{k-1},y^k,\hat{y}^k,y^{k+1},\ldots,y^N$ and $\omega_0$. Therefore, $P$ is in $L_{NS}(\mu)$.

Proof of (iii) and (iv): We prove (iii) and (iv) by providing an example. This example is CHSH team. Recall that in CHSH team, we have two agents with binary measurement and action spaces $\{0,1\}$. Observations are generated independently and uniformly. Hence, there is no state variable in the problem, and so the problem is automatically static. The reward function is defined as
\begin{align}
r(y^1,y^2,u^1,u^2) =
\begin{cases}
1, & \text{if $u^1 \oplus u^2 = y^1 \cdot y^2$} \\
0, & \text{otherwise}.
\end{cases}\nonumber
\end{align}
For this problem, randomized strategic measures $L_C(\mu)$ can achieve $0.5$ \cite{ClHoToWa04}. Quantum-correlated strategic measures can achieve the maximum reward of $2 \sqrt{2}/4$ \cite[Section 6.3.2]{Wat18} (see also Example~\ref{CHSH}). However, non-signaling policies can achieve the maximum reward of $1$ using the following policy, which is called Popescu-Rohrlich (PR) box in the literature \cite{PoRo94}:
\begin{align}
P(u^1,u^2 | y^1,y^2) &=
\begin{cases}
1/2, &\text{if $u^1 \oplus u^2 = y^1 \cdot y^2$} \\
0, &\text{otherwise}
\end{cases}\nonumber
\end{align}
It is very easy to prove that $P$ is non-signaling; that is, $u^i$ is independent of $y^j$ for $i,j=0,1$ and $i \neq j$. The reward of $P$ is $1$, which is the maximum achievable reward by any policy as $0 \leq r \leq 1$. Hence, $P$ is the optimal non-signaling policy. Therefore, for CHSH team, we have 
$$\sup_{P \in L_{NS}(\mu)} \int P(ds) \, r(s) > \sup_{P \in L_Q(\mu)} \int P(ds) \, r(s) > \sup_{P \in L_C(\mu)} \int P(ds) \, r(s)$$ 
that is, non-signaling policies improve the optimal team cost of quantum-correlated policies and quantum correlated policies improve the optimal team cost of randomized policies.
\end{proof}

\syr{
\begin{remark}[On comparing quantum and non-signaling policies]
Note that to implement the optimal non-signaling policy in the CHSH team, agents should communicate their observations to a mediator, and then, mediator directs them to apply either the same actions or different actions based on the product of their observations. This kind of communication is, in general, infeasible for team decision problems. Therefore, although allowing non-signaling correlations among actions of agents enables us to formulate the team problem as a linear program (solvable in a polynomial time for finite case), it is in general not realistic to assume that agents can apply such policies in real-life applications due to communication constraints dictated by decentralized information structure. Therefore, non-signaling relaxation can only be used to provide a lower bound to the original team problem. 

However, quantum correlated relaxation is indeed an \textit{admissible} extension or relaxation of classical team problem since it does not require any communication between agents and a mediator (assuming that we live in the quantum world equipped with appropriate instrumentation!). Hence, the solution or approximation (as explained in the previous section) of quantum correlated team problems will be a significant contribution to team decision theory. 
\end{remark}

\begin{remark}[Information theoretic relaxations]
At the heart of information theory's success is the arrival at single-letter characterizations of optimal information transmission problems which are, in operational formulations, strictly non-convex optimization problems. As detailed out in \cite[Section 5.4]{YukselBasarBook} in the context of optimal quantization problems, information theory convexifies these problem by first relaxing the constraints (such as conditional independence) with mutual information constraints and randomized codes and then showing the attainability of such bounds as the dimension of the problem reaches infinity. Building on this insight, various efforts have been presented for team theoretic problems. 
Notably, Kulkarni and Coleman \cite{kulkarni2014optimizer} studied strategic measures in the context of product measures involving an encoder and a decoder, and a class of team problems, and have considered convexification of the strategic measures by characterizing information channels with information theoretic inequalities. Here, the idea is to abstractly view a channel by its mutual information properties, and thus avoiding the product-structure which ultimately makes $L_R(\mu)$ non-convex.
\end{remark}
}

\subsection{Linear Programming Formulation and its Dual}\label{dual}

Since the non-signaling constraints are linear in $P$, the optimization problem associated with non-signaling strategic measures can be written as a linear program over an appropriate vector spaces. The dual of this linear program and its approximation will provide a lower bound to the original problem. This lower bound, serving as a benchmark, can be quite useful in evaluating how a sub-optimal policy performs in the team problem.

One such linear programming formulation can be done as follows. For any metric space $\mathbb{E}$, let ${\cal M}(\mathbb{E})$ denote the set of finite signed measures on $\mathbb{E}$ and $C(\mathbb{E})$ denotes the vector space of some real measurable functions on $\mathbb{E}$. Depending on the topological properties of the cost function $c$, one should choose $C(\mathbb{E})$ appropriately. For instance, if $c$ is continuous, then $C(\mathbb{E})$ must be chosen as the set of all continuous functions. 

Consider the vector spaces 
$${\cal M}\left(\Omega_0\times\prod_{k=1}^N (\mathbb{Y}^k\times\mathbb{U}^k)\right), \,\, C\left(\Omega_0\times\prod_{k=1}^N (\mathbb{Y}^k\times\mathbb{U}^k)\right).$$ For each $i \in {\cal N}$, we also define 
$${\cal M}\left(\Omega_0\times\prod_{k=1}^N \mathbb{Y}^k\times \prod_{k\neq i} \mathbb{U}^k\right), \,\, C\left(\Omega_0\times\prod_{k=1}^N \mathbb{Y}^k\times\prod_{k\neq i}\mathbb{U}^k\right).$$ 
To ease the notation, we will denote these sets as ${\cal M}, \, C, \, {\cal M}_i, \, C_i$ for $i \in {\cal N}$. 
Let us define bilinear forms on
$({\cal M},C)$ and on $\left(\prod_{i=1}^N {\cal M}_i\times \mathbb{R},\prod_{i=1}^N C_i\times \mathbb{R}\right)$ by letting
\begin{align}
\langle \rho,v  \rangle_{1} &\coloneqq \int v(\omega_0,{\bf y},{\bf u}) \, \rho(d\omega_0,d{\bf y},d{\bf u}) \label{eqqq2}, \\
\left\langle \left(\prod_{i=1}^N \rho_i,a\right) ,\left(\prod_{i=1}^N v_i,b\right)  \right\rangle_{2}  
&\coloneqq \sum_{i=1}^N \int  v_i(\omega_0,{\bf y},{\bf u}^{-i}) \, \rho_i(d\omega_0,d{\bf y},d{\bf u}^{-i}) + a b \label{eqqq1}.
\end{align}
The bilinear forms in (\ref{eqqq2}) and (\ref{eqqq1}) constitute duality between spaces \cite[Chapter IV.3]{Bar02}. Hence, the topologies on these spaces should be understood as the weak topology of the duality induced by these bilinear forms. One should not confuse these topologies with the topologies induced by dual vector spaces.

We define the linear map 
$L: {\cal M} \rightarrow \prod_{i=1}^N {\cal M}_i\times \mathbb{R}$ 
by $L(\rho) = \left(\prod_{i=1}^N L_i(\rho), T(\rho)\right)$, where
\begin{align}
L_i: \rho(d\omega_0,d{\bf y},d{\bf u}) &\mapsto \rho(d\omega_0,d{\bf y},d{\bf u}^{-i}) - \rho(d{\bf y}^{-i},d{\bf u}^{-i}) \, \mu(d\omega_0,dy^i) \nonumber, \, \, i \in {\cal N} \\
T: \rho(d\omega_0,d{\bf y},d{\bf u}) &\mapsto \langle \rho,1  \rangle_{1} \nonumber.
\end{align}
Using $L$, the optimal team cost with non-signaling strategic measures can be written as a linear program as follows:
\begin{align}
(\textbf{NS}) \text{                         }&\text{minimize}_{\rho \in {\cal M}_+} \text{ } \langle \rho,c \rangle_1
\nonumber \\*
&\text{subject to  } L(\rho) = (0,\ldots,0,1). \label{aaaaa}
\end{align}

With this linear programming formulation, the optimal team cost with non-signaling strategic measures can be found in polynomial time if the measurement and action spaces are finite sets. This is not possible for teams with randomized policies \cite{PaTs86} and also class of teams with quantum-correlated relaxations  \cite{KeKoMaToVi11,Vid16,NaVi18,ji2020mip}. Note that if the measurement and action spaces are continuous, the linear program (\textbf{NS}) is infinite dimensional, and so, should be approximated. But, this approximation may not provide a lower bound to the original team problem. To achieve this, we should first formulate the dual of (\textbf{NS}) and approximate the resulting infinite dimensional dual linear program to obtain a lower bound to the original team problem. 

Note that the dual of $L$ is given by $L^*: \prod_{i=1}^N C_i \times \mathbb{R} \rightarrow C$, where 
$$
L^*(v_1,\ldots,v_N,b) \coloneqq \sum_{i=1}^N \left( v_i(\omega_0,{\bf y},{\bf u}^{-i}) - \int_{\Omega_0\times\mathbb{Y}^i} \mu(d\omega_0,dy^i) \, v_i(\omega_0,{\bf y},{\bf u}^{-i}) \right) + b.
$$
Then the dual program of $(\textbf{NS})$ can then be written as \cite[Chapter IV.6]{Bar02}
\begin{align}
(\textbf{NS}^*) \text{                         }&\text{maximize}_{(v_1,\ldots,v_N,b) \in \prod_{i=1}^N C_i \times \mathbb{R} } \text{ } \langle (0,\ldots,0,1), (v_1,\ldots,v_2,b) \rangle_2 = b
\nonumber \\*
&\text{subject to  } L^*(v_1,\ldots,v_N,b) \leq c. \label{aaaab}
\end{align}

Note that this dual linear program $(\textbf{NS}^*)$ is in infinite dimensional spaces that often computationally intractable. By weak duality, the maximum value of $(\textbf{NS}^*)$ is a lower bound to the original team decision problem. Therefore, its approximation is also a lower bound to the original team problem as we are maximizing the objective function as opposed to the linear program (\textbf{NS}). To approximate the infinite dimensional linear program $(\textbf{NS}^*)$, we can use the techniques developed in \cite{EsToDaLy18} in which approximations to infinite dimensional linear programming problems were introduced. Indeed, in this work, authors applied the findings of the paper to approximate the Markov decision processes (MDPs) using the linear programming formulation of the MDPs. Similar to the ($\textbf{NS}^*$), this linear programming problem is defined on a function space. In \cite{EsToDaLy18}, authors first approximate the infinite dimensional function space with a finite dimensional subspace spanned by finitely many independent functions. For instance, one can use finitely many Fourier basis functions to generate this finite dimensional subspace. With this reduction, the problem becomes semi-infinite since the functions should still satisfy the inequality constraint (i.e., $L^*(v_1,\ldots,v_N,b) \leq c$ in ($\textbf{NS}^*$)) for uncountably many variables (i.e., for all $(\omega_0,{\bf y}, {\bf u}) \in \Omega_0 \times \prod_{k=1}^N (\mathbb{Y}^k\times\mathbb{U}^k)$ in ($\textbf{NS}^*$)). The next step is the randomization step; that is, simulate i.i.d. samples using some distribution and let functions satisfy the inequality contraint only for these sampled points. As a result, we obtain a (random) finite dimensional approximation to the original linear program with some probabilistic convergence guarantee. In our setup, by solving this finite dimensional approximation, we can then obtain a lower bound to the original team problem. 

In the following, we will apply a similar reasoning, but to the primal problem to arrive at (as tight as desired) upper bounds, instead.

\section{Finite Approximations of Information Structures via Quantization}\label{finiteSupportApp}

In this section, we consider the finite approximation of static team problems. Since results of this section can also be applied to static reduction of dynamic teams, we suppose that the cost function $c$ also depends on the measurements ${\bf y}$ (which is not the case in the original problem formulation). Recall that, in the independent static reduction of a dynamic team, the reduced cost function $c_s$ is a function of $\omega_0, \, {\bf u}$, and ${\bf y}$. To obtain finite approximation result, the following assumptions are imposed on the components of the model.

\begin{assumption}\label{newas1}
\begin{itemize}
\item [(a)] The cost function $c$ is continuous in $({\bf u},{\bf y})$ for any fixed $\omega_0$. In addition, it is bounded on any compact subset of $\Omega_0 \times \prod_{k=1}^N (\mathbb{Y}^k \times \mathbb{U}^k)$.
\item [(b)] For each $k$, $\mathbb{U}^k$ is a closed and convex subset of a completely metrizable locally convex vector space.
\item [(c)] For each $k$, $\mathbb{Y}^k$ is locally compact.
\item [(d)] For any subset $G$ of $\prod_{k=1}^N \mathbb{U}^k$, the function $w_G(\omega_0,{\bf y}) \coloneqq \sup_{{\bf u} \in G} c(\omega_0,{\bf y},{\bf u})$ is integrable with respect to $\mu(d\omega_0,d{\bf y})$, for any compact subset $G$ of $\prod_{k=1}^N \mathbb{U}^k$ of the form $G = \prod_{k=1}^N G^k$.
\item [(e)] For any $\underline{\gamma} \in {\bf \Gamma}$ with $J(\underline{\gamma})<\infty$ and each $k$, there exists $u^{k,*} \in \mathbb{U}^k$ such that $J(\underline{\gamma}^{-k},\gamma^k_{u^{k,*}}) < \infty$, where $\gamma^k_{u^{k,*}} \equiv u^{k,*}$.
\end{itemize}
\end{assumption}

\smallskip

Note that Assumption~\ref{newas1}-(d),(e) hold if the cost function is bounded. Indeed, conditions in Assumption~\ref{newas1} are quite mild and hold for the celebrated counterexample of Witsenhausen.

In what follows, for any subset $G$ of $\prod_{k=1}^N \mathbb{U}^k$, we let
\begin{align}
{\bf \Gamma}_G &\coloneqq \biggl\{ \underline{\gamma} \in {\bf \Gamma}: \underline{\gamma}\left(\prod_{k=1}^N \mathbb{Y}^k\right) \subset G \biggr\} \nonumber
\end{align}
and
${\bf \Gamma}_{c,G} \coloneqq {\bf \Gamma}_c \cap {\bf \Gamma}_G$,
where ${\bf \Gamma}_c$ denotes the set of continuous strategies. Using these definitions, let us define the following set of strategic measures for any subset $G$ of $\prod_{k=1}^N \mathbb{U}^k$:
\begin{align}
&L_{A}^{G}(\mu) \nonumber \\
&\phantom{x}\coloneqq \bigg\{P \in L_A(\mu): P(B) = \int_{B^0 \times \prod_{k=1}^N A^k} \mu(d\omega_0, d{\bf y}) \, \prod_{k=1}^N 1_{\{u^k = \gamma^k(y^k) \in B^k\}}, \,\, \underline{\gamma} \in {\bf \Gamma}_{G} \bigg\}. \nonumber  
\end{align}
Let $L_{A}^{G,c}(\mu)$ denote the set of strategic measures in $L_{A}^{G}(\mu)$ induced by continuous policies. 

The following result states that, there exists a near optimal strategic measure whose support on the product of action spaces $\prod_{k=1}^N \mathbb{U}^k$ is convex and compact (and thus bounded) subset $G$ of it, and conditional distributions of actions given measurements are induced by continuous policies. 

\begin{proposition}\label{newprop1}
Suppose Assumption~\ref{newas1} holds. Then, for any $\varepsilon>0$ there exists a compact subset $G$ of $\prod_{k=1}^N \mathbb{U}^k$ of the form $G = \prod_{i=1}^N G^i$, where each $G^i$ is convex and compact, such that
$$\inf_{P \in L_{A}^{G,c}(\mu)} \int P(ds) \, c(s)  < J^* + \varepsilon.$$ 
\end{proposition}

\proofsketch
Given any strategic measure, using Assumption~\ref{newas1}-(e) and the fact that every measure on a Borel space is tight \cite[Theorem 3.2]{Par67}, one can construct a strategic measure in $L_A(\mu)$ whose support on the product of action spaces is convex and compact and whose cost is $\varepsilon/2$-close to the cost of the given strategic measure. 

For the new strategic measure, since it has a convex and compact support on the product of action spaces, using Lusin's theorem \cite[Theorem 7.5.2]{Dud02}, we can construct a strategic measure induced by continuous policies whose cost function is $\varepsilon/2$-close to the cost of bounded support strategic measure. 

We can complete the proof by combining these two results. 
\QEDA

\def\intr{\mathop{\rm int}}

\smallskip

Since each $\mathbb{Y}^i$ is a locally compact separable metric space, there exists an increasing sequence of compact subsets $\{K_l^i\}$ such that $K_l^i \subset \intr K_{l+1}^i$ and $\mathbb{Y}^i = \bigcup_{l=1}^{\infty} K_l^i$ \cite[Lemma 2.76]{AlBo06}, where $\intr D$ denotes the interior of the set $D$.

Let $d_i$ denote the metric on $\mathbb{Y}^i$. For each $l\geq1$, let $\mathbb{Y}_{l,n}^i \coloneqq \bigl\{ y_{i,1}, \ldots, y_{i,i_{l,n}} \bigr\} \subset K_l^i$ be an $1/n$-net in $K_l^i$.
Recall that if $\mathbb{Y}_{l,n}^i$ is an $1/n$-net in $K_l^i$, then  for any $y \in K_l^i$ we have
\begin{align}
\min_{z \in \mathbb{Y}_{l,n}^i} d_i(y,z) < \frac{1}{n}. \nonumber
\end{align}
For each $l$ and $n$, let $q_{l,n}^i: K_l^i \rightarrow \mathbb{Y}_{l,n}^i$ be a nearest neighborhood quantizer given by
\begin{align}
q_{l,n}^i(y) = \arg \min_{z \in \mathbb{Y}_{l,n}^i} d_{i}(y,z), \nonumber
\end{align}
where ties are broken so that $q_{l,n}^i$ is measurable.
If $K_l^i=[-M,M] \subset \mathbb{Y}^i = \mathbb{R}$ for some $M \in \mathbb{R}_+$, the finite set $\mathbb{Y}_{l,n}^i$ can be chosen such that $q_{l,n}^i$ becomes a uniform quantizer. We let $Q_{l,n}^i: \mathbb{Y}^i \rightarrow \mathbb{Y}_{l,n}^i$ denote the extension of $q_{l,n}^i$ to $\mathbb{Y}^i$ given by
\begin{align}
Q_{l,n}^i(y) \coloneqq \begin{cases}
q_{l,n}^i(y), &\text{ if } y \in K_l^i, \\
y_{i,0}, &\text{ otherwise},
\end{cases} \nonumber
\end{align}
where $y_{i,0} \in \mathbb{Y}^i$ is some auxiliary element.

Define $\Gamma_{l,n}^i = \Gamma^i \circ Q_{l,n}^i \subset \Gamma^i$; that is, $\Gamma_{l,n}^i$ is defined to be the set of all strategies $\tilde{\gamma}^i \in \Gamma^i$ of the form $\tilde{\gamma}^i = \gamma^i \circ Q_{l,n}^i$, where $\gamma^i \in \Gamma^i$. Define also ${\bf \Gamma}_{l,n} \coloneqq \prod_{i=1}^N \Gamma_{l,n}^i \subset {\bf \Gamma}$. Note that, for any $i=1,\ldots,N$, $\Gamma_{l,n}^i$ is the set of policies for DM~$i$ which can only use the output levels of the quantizer $Q_{l,n}^i$. In other words, in addition to measurement channel $g^i(dy^i|\omega_0)$ between DM~$i$ and the Nature, there is also an analog-to-digital converter (quantizer) between them.

Using these definitions, let us define the following set of strategic measures for any $l$ and $n$:
\begin{align}
&L_{A}^{l,n}(\mu) \nonumber \\
& \phantom{x} \coloneqq \bigg\{P \in L_A(\mu): P(B) = \int_{B^0 \times \prod_{k=1}^N A^k} \mu(d\omega_0, d{\bf y}) \, \prod_{k=1}^N 1_{\{u^k = \gamma^k(y^k) \in B^k\}}, \,\, \underline{\gamma} \in {\bf \Gamma}_{l,n}\bigg\}.\nonumber
\end{align}
The following theorem states that an optimal (or almost optimal) strategic measure can be approximated with arbitrarily small approximation error for the induced costs by strategic measures in $L_{A}^{l,n}(\mu)$ for sufficiently large $l$ and $n$. 

\begin{theorem}\label{newthm1}\cite{saldiyuksellinder2017finiteTeam}
For any $\varepsilon>0$, there exist $(l,n(l))$, a compact subset $G$ of $\prod_{k=1}^N \mathbb{U}^k$ of the form $G = \prod_{i=1}^N G^i$, where each $G^i$ is convex and compact, and $P \in L_{A}^{l,n(l)}(\mu) \bigcap L_{A}^{G}(\mu)$ such that
$$\int P(ds) \, c(s) < J^* + \varepsilon$$ 
\end{theorem}

\proofsketch
Fix any strategic measure $P \in L_A(\mu)$. Then, by Proposition~\ref{newprop1}, there exists a strategic measure $P_G \in L_A^{G,c}(\mu)$ for some compact subset $G$ of $\prod_{k=1}^N \mathbb{U}^k$ of the form $G = \prod_{i=1}^N G^i$, where each $G^i$ is convex and compact, such that
$$
\int P_G(ds) \, c(s) \leq \int P(ds) \, c(s) +\varepsilon/2 .
$$ 
Let $\underline{\gamma} \in {\bf \Gamma}_{G,c}$ be the strategy that induces $P_G$. Then define $\gamma^i_{l,n} = \gamma^i \circ Q^i_{l,n}$ for each $i=1,\ldots,N$. Note that on any compact subset of $\mathbb{Y}^i$, we have 
$
\gamma^i_{l,n} \to \gamma^i
$ 
as $\gamma^i$ is continuous. Let $P_{l,n}$ the strategic measure induced by the strategy $(\gamma^1_{l,n},\ldots,\gamma^N_{l,n})$. 

Note that since $\mu(d\omega_0,d{\bf y})$ is tight, for any $\delta>0$, there exists a compact subset $K$ of $\prod_{k=1}^N \mathbb{Y}^k$ of the form $K = \prod_{i=1}^N K^i$, where each $K^i$ is compact such that $\mu(\Omega_0\times K) \geq 1-\delta$. By choosing a sequence of $\{\delta(l,n)\}$ and the corresponding compact sets $\{K_{l,n}\}$ appropriately,  we can prove that 
$$
\lim_{(l,n)\rightarrow\infty} \int P_{l,n}(ds) \, c(s) = \int P_G(ds) \, c(s),
$$
as $\gamma^i_{l,n} \to \gamma^i$ on each compact set $K_{l,n}$ and $\mu(\Omega_0\times K_{l,n}) \geq 1-\delta(l,n)$. Hence, there exists $(l,n(l))$ such that 
$$
\int P_{l,n(l)}(ds) \, c(s) \leq \int P_G(ds) \, c(s) +\varepsilon/2 .
$$
This completes the proof. 
\QEDA

The above result implies that to compute a near optimal strategic measure for the team problem it is sufficient to consider quantized measurements $\bigl( Q_{l,n}^1(y^1), \ldots, Q_{l,n}^N(y^N) \bigr)$ for sufficiently large $l$ and $n$. Furthermore, this nearly optimal strategic measure can have a compact support of the form $G = \prod_{i=1}^N G^i$ on $\prod_{k=1}^N \mathbb{U}^k$, where $G^i$ is convex and compact for each $i=1,\ldots,N$.

For ease of reference, we define
\[L_{A(q)}(\mu) = \bigcup_{(l,n)} L_{A}^{l,n}(\mu) \]
and we define $L_{A(c)}(\mu)$
to be the subset of $L_{A}(\mu)$ defined in (\ref{measurableStrategic}) where the policies are restricted to be in $\Gamma_c$.

\subsection{Finite measurement approximate models}

In this section, for each $(l,n)$, we define a team model with finite measurement spaces. We prove that, for sufficiently large $l$ and $n$, optimal strategic measure of the team model corresponding to $(l,n)$ will provide a strategic measure to the original team model which is nearly optimal. 

To this end, fix any $(l,n)$. For the pair $(l,n)$, the corresponding finite measurement team model has the following measurement spaces: $\mathbb{Z}_{l,n}^i \coloneqq \{y_{i,0},y_{i,1},\ldots,y_{i,i_{l,n}}\}$ (i.e., the output levels of $Q_{l,n}^i$), $i \in {\cal N}$. The stochastic kernels $g_{l,n}^i(\,\cdot\,|\omega_0)$ from $\Omega_0$ to $\mathbb{Z}_{l,n}^i$ denotes the measurement constraints and
given  by:
\begin{align}
g_{l,n}^i(\,\cdot\,|\omega_0) &\coloneqq \sum_{j=0}^{i_{l,n}} g(S_{i,j}^{l,n}|\omega_0) \, \delta_{y_{i,j}}(\,\cdot\,), \nonumber
\end{align}
where $S_{i,j}^{l,n} \coloneqq \bigl\{ y \in \mathbb{Y}^i: Q_{l,n}^i(y) = y_{i,j} \bigr\}$. Indeed, $g_{l,n}^i(\,\cdot\,|\omega_0)$ is the push-forward of the measure $g^i(\,\cdot\,|\omega_0)$ with respect to the quantizer $Q_{l,n}^i$.

Let $\Phi_{n,l}^i \coloneqq \bigl\{\phi^i: \mathbb{Z}_{l,n}^i \rightarrow \mathbb{U}^i, \text{$\phi^i$ measurable}\bigr\}$ denote the set of measurable policies for DM~$i$ and  let ${\bf \Phi}_{l,n} \coloneqq \prod_{i=1}^N \Phi_{l,n}^{i}$. The cost of this team model is $J_{l,n}: {\bf \Phi}_{l,n} \rightarrow \mathbb{R}_+$ and defined as
\begin{align}
J_{l,n}(\underline{\phi}) \coloneqq \int_{\Omega_0 \times \prod_{i=1}^N \mathbb{Z}_{l,n}^i} c(\omega_0,{\bf y},{\bf u}) \, P_{l,n}(d\omega_0,d{\bf y}), \nonumber
\end{align}
where $\underline{\phi} = (\phi^1,\ldots,\phi^N)$, ${\bf u} = \underline{\phi}({\bf y})$, and 
$$P_{l,n}(d\omega_0,d{\bf y}) = P(d\omega_0) \prod_{i=1}^N g_{l,n}^i(dy^{i}|\omega_0)\eqqcolon \mu_{l,n}(d\omega_0,d{\bf y}).$$ 
For any compact subset $G$ of $\prod_{k=1}^N \mathbb{U}^k$, we also define ${\bf \Phi}_{l,n}^G \coloneqq \{\underline{\phi} \in {\bf \Phi}_{l,n}: \underline{\phi}(\prod_{i=1}^N \mathbb{Z}_{l,n}^i) \subset G\}$.

In order to obtain the approximation result, we need to impose the following additional assumption.

\begin{assumption}\label{nnewas1}
For any compact subset $G$ of $\prod_{k=1}^N \mathbb{U}^k$ of the form $G = \prod_{i=1}^N G^i$, we assume that the function $w_G$ is uniformly integrable with respect to the measures $\{\mu_{l,n}\}$; that is,
\begin{align}
\lim_{R\rightarrow\infty} \sup_{l,n} \int_{\{w_G > R\}} w_G(\omega_0,{\bf y}) \text{ } d\mu_{l,n} = 0. \nonumber
\end{align}
\end{assumption}

This assumption is quite mild and satisfied by the celebrated counterexample of Witsenhausen. 

Let $T_A(\mu_{l,n})$ denote the set of strategic measures induced by measurable policies for this team model. Let $T_A^G(\mu_{l,n})$ be the strategic measures in $T_A(\mu_{l,n})$ whose support on $\prod_{k=1}^N \mathbb{U}^k$ is a subset $G$.   

The following theorem is the main result of this section. It states that to compute a near optimal strategic measure for the original team problem, it is sufficient to compute an optimal strategic measure for the finite measurement team problem corresponding to sufficiently large $l$ and $n$.

\begin{theorem}\label{newthm2}\cite{saldiyuksellinder2017finiteTeam}
Suppose Assumptions~\ref{newas1} and \ref{nnewas1} hold. Then, for any $\varepsilon>0$, there exists a pair $(l,n(l))$ and a compact subset $G = \prod_{i=1}^N G^i$ of $\prod_{k=1}^N \mathbb{U}^k$ such that an optimal (or almost optimal) strategic measure $P^{l,n(l)}$ in the set $T_A^G(\mu_{l,n(l)})$ for the $(l,n(l))$ team is $\varepsilon$-optimal for the original team problem when $P^{l,n(l)}$ is extended to $\Omega \times \prod_{k=1}^N (\mathbb{Y}^k\times \mathbb{U}^k)$ via quantizers $Q_{l,n(l)}^i$; that is,
$$
P_{\mathrm{ex}}^{l,n(l)}(\cdot) = \int_{\cdot} \mu(d\omega_0, d{\bf y}) \, \prod_{k=1}^N 1_{\{u^k = \gamma^k\circ Q_{l,n(l)}^k(y^k) \in \cdot\}}
$$
where 
$$
P^{l,n(l)}(B) = \int_{\cdot} \mu_{l,n(l)}(d\omega_0, d{\bf y}) \, \prod_{k=1}^N 1_{\{u^k = \gamma^k(y^k) \in \cdot\}}
$$
\end{theorem}

\begin{proof}
We first prove the following fact. Let $\{\underline{\phi}_{l,n}\}$ be a sequence of strategies such that $\underline{\phi}_{l,n} \in {\bf \Phi}_{l,n}^G$, where $G= \prod_{i=1}^N G^i$ and each $G^i$ is convex and compact. For each $l$ and $n$, define $\underline{\gamma}_{l,n} \coloneqq \underline{\phi}_{l,n} \circ Q_{l,n}$, where $Q_{l,n} \coloneqq (Q_{l,n}^1,\ldots,Q_{l,n}^N)$. Then, we have
\begin{align}
\lim_{l,n\rightarrow\infty} |J_{l,n}(\underline{\phi}_{l,n}) - J(\underline{\gamma}_{l,n})| =  0. \nonumber
\end{align}
Indeed, let us introduce the following finite measures on $\Omega_0 \times \prod_{k=1}^N\mathbb{Y}^k$:
\begin{align}
\mu_G(S) &\coloneqq \int_{S} w_G(\omega,{\bf y}) \text{ } d\mu, \nonumber \\
\mu_G^{l,n}(S) &\coloneqq \int_{S} w_G(\omega,{\bf y}) \text{ } d\mu_{l,n}.  \nonumber
\end{align}
Since $\mu_{l,n}$ converges to $\mu$ weakly, by \cite[Theorem 3.5]{Ser82} and Assumption~\ref{nnewas1} we have $\mu_G^{l,n} \rightarrow \mu_G$ weakly as $l,n \rightarrow \infty$. Hence, the sequence $\{\mu_G^{l,n}\}$ is tight. Therefore, there exists a compact subset $K$ of $\Omega_0 \times \prod_{k=1}^N\mathbb{Y}^k$ such that $\mu_G(K^c)< \varepsilon/2$ and $\mu_G^{l,n}(K^c) < \varepsilon/2$ for all $l,n$. Then, we have
\begin{align}
&|J_{l,n}(\underline{\phi}_{l,n}) - J(\underline{\gamma}_{l,n})| \nonumber \\
&\phantom{xx}= \bigg| \int_{\Omega_0 \times \prod_{k=1}^N\mathbb{Y}^k} c(\omega_0,Q_{l,n}({\bf y}),\underline{\gamma}_{l,n}({\bf y})) \text{ } d\mu - \int_{\Omega_0 \times \prod_{k=1}^N\mathbb{Y}^k} c(\omega_0,{\bf y},\underline{\gamma}_{l,n}({\bf y})) \text{ } d\mu \biggr|  \nonumber \\
&\phantom{xx}\leq \int_{K} \bigl| c(\omega_0,Q_{l,n}({\bf y}),\underline{\gamma}_{l,n}({\bf y})) - c(\omega_0,{\bf y},\underline{\gamma}_{l,n}({\bf y})) \bigr| \text{ } d\mu \nonumber \\
&\phantom{xxxxxxxxxxxx}+ \int_{K^c} w_G(\omega_0,{\bf y}) \text{ } d\mu + \int_{K^c} w_G(\omega_0,{\bf y}) \text{ } d\mu_{l,n}. \nonumber
\end{align}
The first term in the last expression goes to zero as $l,n\rightarrow\infty$ by dominated convergence theorem and the fact that $c$ is bounded and continuous in ${\bf y}$. The second term is less than $\varepsilon$ by tightness. Since $\varepsilon$ is arbitrary, this completes the proof of the fact. Note that this fact implies that for any $P^{l,n} \in T_A^G(\mu_{l,n})$, we have 
$$
\lim_{l,n \rightarrow \infty} \left|\int P^{l,n}(ds) \, c(s) - \int P^{l,n}_{\mathrm{ex}}(ds) \, c(s) \right| = 0. 
$$

Note that using this fact, it is straightforward prove that 
$$
\liminf_{l,n \rightarrow \infty} \inf_{P \in T_A^G(\mu_{l,n})}  \int P(ds) \, c(s) \geq J^*. 
$$
Since the extension of the strategic measures $T_A^G(\mu_{l,n})$ to $\Omega_0 \times \prod_{k=1}^N(\mathbb{Y}^k\times\mathbb{U}^k)$ forms the set of strategic measures $L_A^{l,n}(\mu) \bigcap L_A^{G}(\mu)$, by Theorem~\ref{newthm1}, we can complete the proof. Indeed, let $P^{l,n(l)}_{\mathrm{ex}} \in L_A^{l,n}(\mu) \bigcap L_A^{G}(\mu)$ be the $\varepsilon$ optimal policy for the original problem. Then we have 
\begin{align}
J^* + \varepsilon &\geq \limsup_{l,n\rightarrow\infty} \int P^{l,n(l)}_{\mathrm{ex}}(ds) \, c(s) \nonumber \\
&= \limsup_{l,n\rightarrow\infty} \int P^{l,n(l)}(ds) \, c(s) \nonumber \\
&\geq \limsup_{l,n \rightarrow \infty} \inf_{P \in T_A^G(\mu_{l,n})}  \int P(ds) \, c(s). \nonumber 
\end{align}
Hence 
$$
J^*+\varepsilon \geq \limsup_{l,n \rightarrow \infty} \inf_{P \in T_A^G(\mu_{l,n})}  \int P(ds) \, c(s) \geq \liminf_{l,n \rightarrow \infty} \inf_{P \in T_A^G(\mu_{l,n})}  \int P(ds) \, c(s) \geq J^*.
$$
This completes the proof. 
\end{proof}

We interpret this result as saying that the space of information structures obtained by quantizing the original information structure is approximately optimal: denseness under weak convergence implies denseness of strategic measures as far as optimality is concerned. 

\begin{remark}
Results proved in this section about approximation of static team problems can be directly applied to static reductions of dynamic teams if static reductions satisfy Assumption~\ref{newas1} and Assumption~\ref{nnewas1}. In particular, above approximation results hold for the celebrated counterexample of Witsenhausen.  
\end{remark}


\section{Product Topology induced by Individual Policies and Relations with Information Structures}\label{TopologyOnPolicies}

Instead of strategic measures, one may choose to directly induce topologies on control policies or strategies alone. We will see that such an approach will let us arrive at complementary conditions compared with what we have studied earlier. In the following, we first revisit some classical results in optimal control theory via this approach. 

\subsection{Some remarks on classical deterministic and stochastic control and Young measures}\label{YoungM}

It is instructive to revisit various control topologies that are already well-known in classical control theory (when there is a single controller who has access to the state variable). In deterministic nonlinear, geometric, and continuous-time control, properties on stabilizability, controllability, and reachability are drastically impacted by the restrictions on the classes of allowed controls (e.g., continuous, Lipschitz, finitely differentiable, or smooth control functions in the state or time when control is open-loop \cite{ryan1994brockett,brockett1983asymptotic,sontag2013mathematical,jafarpour2016locally}) and naturally the control topology induced is dictated by the class of admissible controls. 

For optimal control, to allow for continuity/compactness arguments, apriori imposing compactness over spaces of measurable functions would be an artificial restriction, and the use of powerful theorems such as the Arzela-Ascoli theorem which necessarily entail (usually very restrictive and suboptimal) conditions on continuity properties of the considered policies. In deterministic optimal control theory, relaxed controls \cite{young1937generalized}  \cite{warga2014optimal} allow for this machinery to be applied with no artificial restrictions on the classes of control policies considered; these are known as Young measures.

Let us consider an open-loop controller, where the control is only a function of the time variable. We let $\nu(dt,du)$ be a measure on $[0,T] \times \mathbb{U}$ where the first marginal $\lambda(dt)$ is the normalized Lebesge measure on time interval $[0,T]$ and let $\nu(du | t) = 1_{\{\gamma(t) \in du\}}$ be the conditional  measure induced by deterministic open loop control. So, any deterministic open-loop control is embedded via:
\[\nu(dt,du) = \lambda(dt) \, 1_{\{\gamma(t) \in du\}}.\]
If allow for randomized policies, we obtain the set ${\cal P}_{\lambda}([0,T]\times\mathbb{U})$ of all probability measures with fixed marginal on $[0,T]$. This set is weakly closed, whose extreme points are those induced by deterministic policies (as also was discussed in the context of (\ref{extremePointQuan0})). Thus, any deterministic optimal control problem, which can be written in an integral form and have lower semi-continuous cost functions in actions, will have an optimal solution, which will then be deterministic as these  form the extreme points of randomized controls. It can also in fact be shown that such policies are dense in the space of randomized policies, in addition to these policies forming the extreme points in the set of randomized policies (see e.g., \cite[Proposition 2.2]{beiglbock2018denseness} \cite{lacker2018probabilistic}, \cite[19, Theorem 3]{milgrom1985distributional}, but also many texts in optimal stochastic control where denseness of deterministic controls have been established inside the set of relaxed controls \cite{borkar1988probabilistic}). We refer the reader to \cite{mascolo1989relaxation} (see also the review paper \cite{lacker2018probabilistic}) for further discussion. 

The following example builds on these, with somewhat different arguments. Let $\mathbb{X} = \mathbb{R}, \mathbb{U} = [0,1]$, and let $f: \mathbb{X} \times \mathbb{U} \to [0,1]$ and $c:\mathbb{X} \times \mathbb{U} \to [0,1]$ be measurable functions continuous in the control action variable. Consider the following optimal control problem:
\begin{eqnarray}
\inf_{\substack{\gamma: \mathbb{X} \to \mathbb{U}\\ u_t = \gamma(x_t)}} \int_0^1 c(x_t,u_t) \, \lambda(dt) \label{cost1}
\end{eqnarray}
subject to 
\begin{eqnarray}
\frac{dx}{dt} = f(x_t,u_t)\label{const1}
\end{eqnarray}
The natural space to consider is the set of all control functions which depends on the current state, where the only restriction is measurability. However, allowing for measurability only does not facilitate continuity/compactness arguments since, as noted above, imposing compactness on a space of functions is an unnecessarily restrictive condition. \sy{Accordingly, one often cites appropriate but tedious measurable selection theorems building on optimality equations through dynamic programming}. 

 On the other hand, every deterministic function of state can be expressed as a deterministic function of time, and so, be considered open-loop. Accordingly, we consider open loop controls and those which are relaxed. Let ${\cal P}_{\lambda}([0,T]\times\mathbb{U})$ be the set of relaxed open loop policies (Young measures). Now consider the space $C([0,1];\mathbb{X}) \times {\cal P}_{\lambda}([0,T]\times\mathbb{U})$, where $C([0,1];\mathbb{X})$ is the space of continuous functions from $[0,1]$ to $\mathbb{X}$. 
We endow this space with the product topology with the first component being under the supremum norm and the second under Prohorov metric (or any weak convergence inducing metric). Note now that the cost (\ref{cost1}) is continuous on $C([0,1];\mathbb{X}) \times {\cal P}_{\lambda}([0,T]\times\mathbb{U})$. Note that since $f$ is uniformly bounded, we have that the set ${\cal A}$ of all admissible sample paths of the state $x: [0,1] \rightarrow \mathbb{X}$ is equicontinuous, and so, by the Arzela-Ascoli theorem, ${\cal A}$ is relatively compact in $C([0,1];\mathbb{X})$. Accordingly, our space of interest ${\cal A} \times {\cal P}_{\lambda}([0,T]\times\mathbb{U})$ is a relatively compact subset of $C([0,1];\mathbb{X}) \times {\cal P}_{\lambda}([0,T]\times\mathbb{U})$.

Define now $$H = \left\{(x,m) \in C([0,1];\mathbb{X}) \times {\cal P}_{\lambda}([0,T]\times\mathbb{U}): x_t - \int_0^t f(x_s,u) \, m_s(du) \, \lambda(ds) = 0\right\},$$ where $m_s(du) = m(du | s)$. This set is closed under the topology defined on $C([0,1];\mathbb{X}) \times {\cal P}_{\lambda}([0,T]\times\mathbb{U})$ and is a subset of $C([0,1];\mathbb{X}) \times {\cal P}_{\lambda}([0,T]\times\mathbb{U})$. Hence, $H$ is compact. Now, the problem then is to find an optimal $(x,m) \in H$ which minimizes (\ref{cost1}), reformulated as: 
\[\inf_{(x,m) \in H} \int_0^1 c(x_t,u) \, m(dt,du) \]
This is continuous in $(x,m)$ by an application of the generalized weak convergence theorem under continuous convergence \cite[Theorem 3.5]{serfozo1982convergence} or \cite[Theorem 3.5]{Lan81}. Therefore, there exists an optimal solution to the problem.

Note that as has been reported in various literatures (e.g. in optimal quantization \cite[p. 878]{YukselOptimizationofChannels}) that the set of deterministic controls is not weakly closed (or setwise) under Young's relaxation. In fact, there exist problems where an optimal solution exists among relaxed controls but not in deterministic controls \cite{lacker2018probabilistic}. 

On the other hand, in the continuous-time stochastic context, the analysis can be quite subtle due to the fact that the control policy (only restricted to be measurable in general) may violate conditions needed for the existence of strong solutions for a given stochastic differential equation since the control policy may couple the state dynamics with the past in an arbitrarily complicated, though measurable, way and hence violating the existence conditions for strong solutions to stochastic differential equations \cite{lindquist1973feedback} \cite{KushnerControlStochasticSC} (e.g., in \cite{wonham1968separation} Lipschitz continuity conditions are imposed, where Lipschitz property holds in the control when viewed as a map from the normed linear space of continuous functions on measurements to control actions). To avoid such technical issues on strong solutions, relaxed solution concepts were introduced and studied in the literature based on the measure transformation technique due to Girsanov \cite{benevs1971existence,davis1972information,davis1973dynamic} which allows the control to be a function of an independent Brownian innovations process. These approaches require strong absolute continuity conditions on the measurement process which may not be always applicable. See \cite{lacker2018probabilistic} for a detailed analysis for controlled stochastic differential systems. As we will see later in the paper, in the partially observed setup, a further interpretation of relaxed controls, called {\it wide sense admissible controls} \cite{FlPa82,Bor00,Bor07} has been utilized to arrive at existence results on optimal control policies.

As we will see in the following, adding more general, non-classical, information structures will entail further intricacies but also facilitate additional creativity for the analysis in optimal stochastic control. 

\subsection{Trading-off continuity with compactness: a product space approach on individual control policies}\label{tradeOffContComp}

To make the ideas in this section more explicit, consider the set of individually randomized strategic measures $L_R(\mu)$ studied in Section \ref{setsOfStrategicMeasures}. In this set, a policy of a DM~$k$ is a stochastic kernel $\Pi^k(du^k|y^k)$ from $\mathbb{Y}^k$ to $\mathbb{U}^k$. Instead of studying continuity and compactness properties regarding strategic measures (joint probability measures on the state, measurement and action spaces), we can focus on individual control policies and measurements, and their Cartesian products. Let us denote, as in Section \ref{distributionalPolicy},
\[\Theta^i = \{v \in {\cal P}(\mathbb{Y}^i\times \mathbb{U}^i): v(dy^i \times \mathbb{U}^i) = \xi^i(dy^i)\}\]
to be the distributional strategy for DM~$i$, which is the joint probability measure with a fixed marginal $\xi^i$ on the measurement variable $y^i$. Note that if the problem is an independent static reduction of a dynamic team, then one can choose $\xi^i$ as the reference measure $Q^i$ in the static reduction.  

By our earlier analysis, we know that $\Theta^i$ (a subset of $\Theta^i$ that is sufficient for optimality) is compact under the weak topology or the $w$-$s$ topology if $\mathbb{U}^i$ is compact (if this subset satisfies some moment condition). Hence, since the marginal on measurements are fixed, we can view the weak topology or $w$-$s$ topology on $\Theta^i$ as a topology on the set of policies $\Gamma^i$. In particular, any sequence of policies:
\[\underline{\gamma}_n := (\gamma^1_n,\cdots,\gamma^N_n)\]
will have a converging subsequence, $\underline{\gamma}_{n_k}$ to a limit $\underline{\gamma}$, and therefore, closedness (in particular compactness) is not an issue under this topology. In other words, the IS is always preserved under this topology, which is not the case in strategic measure approach as it has been shown in Theorem~\ref{counterEx}. 

Note that if the measurements are independent, then strong existence results can be obtained using this approach. This is precisely the same condition obtained via the strategic measures approach studied in Theorem \ref{existenceRelaxed3}.

On the other hand, what is not clear is, in general, whether we have
\[\lim_{n \to \infty} J(\underline{\gamma}_n) =J(\underline{\gamma})\]
under this topology (this is always true when we place a topology on the set of strategic measures). In particular, if we have a general $\mu(d\omega_0,d{\bf y})$, it is not clear if we have this {\it joint continuity} condition. The following example will demonstrate this subtlety through a negative implication. Therefore, there is a trade-off between continuity and compactness in these approaches. Strategic measure approach provides continuity but may not preserve IS (so compactness), and conversely, topology on individual policies leads to compact strategy space but lose the continuity.

\begin{example}
Consider the counterexample presented in Theorem \ref{counterEx}. By considering each individual decision makers' policy separately, we have that each will converge individually; that is,
$$
P_n(dy^1,du^1) \rightarrow P(dy^1,du^1), \,\, P_n(dy^2,du^2) \rightarrow P(dy^2,du^2). 
$$
Now, consider a cost function $c(\omega_0,u^1,u^2)=(u^1-u^2)^2$. Note that, for each $n$, under $P_n$, actions $u^1$ and $u^2$ are the same (i.e., $u^1 = u^2$) given $y^1=y^2=y$, and so, 
$J(\underline{\gamma}_n) = \int P_n(ds) \, c(s) = 0$ for every $n$. However, in the limit where $u^i$ and $y=y^i$ ($i=1,2$) are independent under $P$, we have $J(\underline{\gamma}_n) = \int P(ds) \, c(s) = 1/2$. This is a consequence of viewing individual policies separately without their joint convergence properties. Thus, we do not even have lower semi-continuity. Hence, while we have established compactness, we have lost lower semi-continuity by placing the topologies directly on control policies and not on the strategic measures.
\end{example}


%
%

\subsection{Policies defined by conditional independence given measurements}\label{univDP}

\cite{YukselWitsenStandardArXiv} considered the following topology on control policies, while developing a universal dynamic programming algorithm applicable to any sequential decentralized stochastic control problem, generalizing Witsenhausen's program \cite{WitsenStandard} which was tailored primarily for countable probability spaces.

Define
\begin{itemize}
\item[(i)] State: $x_t =  \{\omega_0,u^1,\cdots,u^{t-1},y^1,\cdots,y^{t}\}$, $1 \leq t \leq N$.
\item[(i')] Extended State: $\pi_t \in {\cal P}(\Omega_0 \times \prod_{i=1}^{t} \mathbb{Y}^i \times \prod_{i=1}^{t-1} \mathbb{U}^i)$ where, for Borel $B \in \Omega_0 \times \prod_{i=1}^{t} \mathbb{Y}^i \times \prod_{i=1}^{t-1} \mathbb{U}^i
$, \[\pi_t(B) := E_{\pi_t}[1_{ \{ (\omega_0,y^1,\cdots,y^t; u^1,\cdots,u^{t-1}) \in B  \} } ].\]
Thus, $\pi_t \in {\cal P}(\Omega_0 \times \prod_{i=1}^{t} \mathbb{Y}^i \times \prod_{i=1}^{t-1} \mathbb{U}^i)$ where the space of probability measures is endowed with the weak convergence topology.
\item[(ii)]
 Control Action: Given $\pi_t$, $\hat{\gamma}^t$ is a probability measure in ${\cal P}(\Omega_0 \times \prod_{k=1}^{t} \mathbb{Y}^k \times \prod_{k=1}^{t} \mathbb{U}^k)$ that satisfies the conditional independence relation: 
\[u^t \leftrightarrow y^t \leftrightarrow x_t = (\omega_0,y^{1},\cdots,y^{t}; u^1,\cdots,u^{t-1})\]
(that is, for every Borel $B \in \mathbb{U}^i$, almost surely under $\hat{\gamma}^t$, the following holds:
\[P(u^t \in B| y^t , (\omega_0,y^{1},\cdots,y^{t}; u^1,\cdots,u^{t-1})) = P(u^t \in B| y^t)\]
with the restriction
\[x_t \sim \pi_t.\]
Denote with $\Gamma^t(\pi_t)$ the set of all such probability measures. Any $\hat{\gamma}^t \in \Gamma^t(\pi_t)$ defines, for almost every realization $y^t$, a conditional probability measure on $\mathbb{U}^t$. 
When the notation does not lead to confusion, we will denote the action at time $t$ by $\gamma^t(du^t|y^t)$, which is understood to be consistent with $\hat{\gamma}^t$. 
\item[(ii')] Alternative Control Action for Static Teams with Independent Measurements: Given $\pi_t$, $\hat{\gamma}^t$ is a probability measure on $\mathbb{Y}^t \times \mathbb{U}^t$ with a fixed marginal $P(dy^t)$ on $\mathbb{Y}^t$, that is $\pi^{\mathbb{Y}^t}_t(dy^t)=P(dy^t)$. Denote with $\Gamma^t(\pi_t^{\mathbb{Y}^t})$ the set of all such probability measures. As above, when the notation does not lead to confusion, we will denote the action at time $t$ by $\gamma^t(du^t|y^t)$, which is understood to be consistent with $\hat{\gamma}^t$. In particular, $(y^t,u^t)$ is independent of $(y^k,u^k)$ for $k \neq t$.
\end{itemize}

With the control actions defined as in the above \cite{YukselWitsenStandardArXiv} developed a universal dynamic program for any sequential decentralized stochastic control and established, as a corollary of the program, further existence results, one which is essentially identical to that presented in \ref{existenceRelaxed3}, but slightly more restrictive in that the cost function is assumed to be continuous in all of its arguments. 

\begin{theorem}\cite{YukselWitsenStandardArXiv}
\begin{itemize}
\item[(i)] Under the kernel (\ref{kernelDefn}) and controlled Markov construction presented, the optimal team problem admits a well-defined backwards-induction (dynamic programming) recursion. 
\item[(ii)] In particular, if the problem is independent static-reducible, actions are compact-valued and the cost function is continuous, an optimal policy exists and the value function is continuous in the prior (that is, in the distribution of primitive noise variables) under weak convergence. 
\end{itemize}
\end{theorem}

\begin{remark}The above construction is related to an interpretation put forward by Witsenhausen in his standard form \cite{WitsenStandard} where all the uncertainly is embedded into the initial state and the controlled system evolves deterministically. Witsenhausen had considered only countable probability spaces for an optimality analysis. 
\end{remark}

\subsection{Weak-$^*$ topology on Randomized Policies} 


In this section, inspired by the topology on policies for classical stochastic control problems constructed in Borkar \cite{borkar1989topology} (see also \cite{arapostathis2012ergodic}), we study a topology on the set of randomized policies for static teams or dynamic teams that admit static reduction introduced in \cite{SaldiArXiv2017}. With an abuse of notation, we denote the set of randomized policies for DM~$i$ by $\Gamma^i$, which can be written as
%

\begin{align}
\Gamma^i = \bigg\{\gamma^i: \text{$\gamma^i$ is a measurable function from $\mathbb{Y}^i$ to ${\cal P}(\mathbb{U}^i)$} \bigg\}, \nonumber
\end{align}
where ${\cal P}(\mathbb{U}^i)$ is endowed with Borel $\sigma$-algebra generated by weak convergence topology. We will use the latter characterization of $\Gamma^i$ when introducing the topology. 

For a metric space $\mathbb{E}$, let $C_0(\mathbb{E})$ be the Banach space of all continuous real functions on $\mathbb{E}$ vanishing at infinity endowed with the norm
\begin{align}
\|g\|_{\infty} = \sup_{e \in \mathbb{E}} |g(e)|. \nonumber
\end{align}
Recall that ${\cal M}(\mathbb{E})$ and ${\cal P}(\mathbb{E})$ denote the set of all finite signed measures and probability measures on $\mathbb{E}$. Let $\mathbb{E}_1$ and $\mathbb{E}_2$ be two metric spaces. For any $\mu \in {\cal M}(\mathbb{E}_1 \times \mathbb{E}2)$, 
$\Proj_{\mathbb{E}_1}(\mu)(\,\cdot\,) = \mu(\,\cdot\times \mathbb{E}_2)$ be the marginal of $\mu$ on $\mathbb{E}_1$. For a Banach space $\mathbb{K}$, let $\mathbb{K}^*$ denote its topological dual. For any $i$, we now introduce a topology on $\Gamma^i$. To this end, we impose the following assumption.

\begin{assumption}\label{top-as1}
\begin{itemize}
\item [(a)] $\Omega_0$, $\mathbb{Y}^i$, and $\mathbb{U}^i$ ($i=1,\ldots,N$) are locally compact.
\item [(b)] For all $i$, the measurement channel $g^i(dy_i|\omega_0) = g^i(y_i,\omega_0) \, \mu^i(dy^i)$ for some probability measure $\mu^i$ on $\mathbb{Y}^i$; that is, $g^i(dy_i|\omega_0)$ is absolutely continuous with respect to $\mu^i(dy^i)$ for all $\omega_0$.
\end{itemize}
\end{assumption}

\def\ess{{\mathrm{ess}}}

Note that by Riesz representation theorem \cite[Theorem 7.17]{Fol99}, the topological dual of $(C_0(\mathbb{U}^i),\|\,\cdot\,\|_{\infty})$ is $({\cal M}(\mathbb{U}^i),\|\,\cdot\,\|_{TV})$, where $\|\cdot\|_{TV}$ is the total variation norm on ${\cal M}(\mathbb{U}^i)$.

We now define the set of $w^*$-measurable functions from $\mathbb{Y}^i$ to ${\cal M}(\mathbb{U}^i)$. Later, we will prove that the set of randomized policies $\Gamma^i$ for DM~$i$ is a bounded subset of this set. A function $\gamma: \mathbb{Y}^i \rightarrow {\cal M}(\mathbb{U}_i)$ is called $w^*$-measurable if the mapping
\begin{align}
\mathbb{Y}^i \ni y \mapsto \int d\gamma(y)(du) \, g(u) \in \mathbb{R} \nonumber
\end{align}
is measurable for all $g \in C_0(\mathbb{U}^i)$ \cite[p. 18]{CeMe97}. Let ${\cal L}\bigl(\mu^i,{\cal M}(\mathbb{U}^i)\bigr)$ denote the set of all such functions. With this definition, we now define the following set
\begin{align}
&{\cal L}_{\infty}\bigl(\mu^i,{\cal M}(\mathbb{U}^i)\bigr) \coloneqq \biggl\{ \gamma \in {\cal L}\bigl(\mu^i,{\cal M}(\mathbb{U}^i)\bigr): \|\gamma\|_{\infty} \coloneqq \ess \sup_{y \in \mathbb{Y}^i} \|\gamma(y)\|_{TV} <\infty \biggr\}, \nonumber
\end{align}
where $\ess \sup$ is taken with respect to the measure $\mu^i$. Recall that $\mu_i$ is the reference probability measure in Assumption~\ref{top-as1}-(b) for the measurement  channel $g^i$.


Let $L_1\bigl(\mu^i,C_0(\mathbb{U}^i)\bigr)$ denote the set of all Bochner-integrable \cite{DiUh77} functions from $\mathbb{Y}^i$ to $C_0(\mathbb{U}^i)$ endowed with the norm
\begin{align}
\|f\|_1 \coloneqq \int_{\mathbb{Y}^i} \|f(y)\|_{\infty} \, \mu^i(dy^i). \nonumber
\end{align}

By using the fact $C_0(\mathbb{U}^i)^* = {\cal M}(\mathbb{U}^i)$, one can prove that the topological dual of $\left(L_1\bigl(\mu^i,C_0(\mathbb{U}^i)\bigr),\|\cdot\|_1\right)$ can be identified by $\left({\cal L}_{\infty}\bigl(\mu^i,{\cal M}(\mathbb{U}^i)\bigr),\|\cdot\|_{\infty}\right)$ \cite[Theorem 1.5.5, p. 27]{CeMe97}; that is,
$$
L_1\bigl(\mu^i,C_0(\mathbb{U}^i)\bigr)^* = {\cal L}_{\infty}\bigl(\mu^i,{\cal M}(\mathbb{U}^i)\bigr).
$$
Using this identification, we equip ${\cal L}_{\infty}\bigl(\mu^i,{\cal M}(\mathbb{U}^i)\bigr)$ with $w^*$-topology induced by $L_1\bigl(\mu^i,C_0(\mathbb{U}^i)\bigr)$. We write $\gamma_n \rightharpoonup^* \gamma$, if $\gamma_n$ converges to $\gamma$ in ${\cal L}_{\infty}\bigl(\mu^i,{\cal M}(\mathbb{U}^i)\bigr)$ with respect to $w^*$-topology. As noted earlier, for this topology, we have been inspired by the topology introduced in \cite{arapostathis2012ergodic,borkar1989topology}, where in these works, a similar topology is introduced for randomized Markov policies to study continuous-time stochastic control problems with average cost optimality criterion.

Note that we can identify the set of randomized policies $\Gamma^i$ as a bounded subset of ${\cal L}_{\infty}\bigl(\mu^i,{\cal M}(\mathbb{U}^i)\bigr)$:
\begin{align}
\Gamma^i = \biggl\{ \gamma \in {\cal L}_{\infty}\bigl(\mu^i,{\cal M}(\mathbb{U}^i)\bigr): \gamma(y) \in {\cal P}(\mathbb{U}^i) \text{ } \mu^i-\text{a.e.} \biggr\}. \nonumber
\end{align}
Hence, we can equip $\Gamma^i$ with the relative $w^*$-topology inherited by $w^*$-topology on ${\cal L}_{\infty}\bigl(\mu^i,{\cal M}(\mathbb{U}^i)\bigr)$.

Now, we derive some properties of this topology on $\Gamma^i$. Firstly, since ${\cal P}(\mathbb{U}^i)$ is bounded in ${\cal M}(\mathbb{U}^i)$, by Banach-Alaoglu Theorem \cite[Theorem 5.18]{Fol99}, $\Gamma^i$ is  relatively compact with respect to $w^*$-topology. Since $C_0(\mathbb{U}^i)$ is separable, then by \cite[Lemma 1.3.2]{HeLa03}, $\Gamma^i$ is also relatively sequentially compact. However, $\Gamma^i$ is not closed with respect to $w^*$-topology. Indeed, let $\mathbb{Y}^i = \mathbb{U}^i = \mathbb{R}$. Define $\gamma_n(y^i)(\,\cdot\,) \coloneqq \delta_n(\,\cdot\,)$ and $\gamma(y^i)(\,\cdot\,) \coloneqq 0(\,\cdot\,)$, where $\delta_n$ denotes Dirac-delta measure at point $n$ and $0(\,\cdot\,)$ denotes the degenerate measure on $\mathbb{U}^i$; that is, $0(D)=0$ for all $D \in {\cal B}(\mathbb{R})$. Let $g \in L_1\bigl(\mu^i,C_0(\mathbb{U}^i)\bigr)$. Then we have
\begin{align}
\lim_{n\rightarrow\infty}& \int_{\mathbb{Y}^i} \left( \int_{\mathbb{U}^i}  \gamma_n(y)(du) \, g(y)(u) \right) \mu^i(dy) = \lim_{n\rightarrow\infty} \int_{\mathbb{Y}^i} g(y)(n) \, \mu^i(dy) \nonumber \\
&= \int_{\mathbb{Y}^i} \lim_{n\rightarrow\infty} g(y)(n) \, \mu^i(dy) \text{ (as $\|g(y)\|_{\infty}$ is $\mu^i$-integrable)} \nonumber \\
&= 0 \text{ (as $g(y) \in C_0(\mathbb{U}^i)$)}. \nonumber
\end{align}
Hence, $\gamma_n \rightharpoonup^* \gamma$. But, $\gamma \notin \Gamma^i$, and so, $\Gamma^i$ is not closed. 

Thus, as opposed to the topology introduced in Section~\ref{TopologyOnPolicies} on randomized policies, here the set of policies is not closed under $w^*$-topology; that is, IS may not be preserved. However, as we will see in the sequel, under this topology, the cost $J$ is lower semi-continuous, which is in general not the case for the topology introduced in Section~\ref{TopologyOnPolicies}. However, by imposing additional assumptions on the components of the problem, we may first ensure the closedness of the set policies or subset of it that is sufficient for optimality under this topology using lower semi-continuity of the cost. Then, we can prove the existence of optimal policy. This result will be complementary to the existence results presented in Section~\ref{regularity}. 

In the remainder of this section, $\Gamma^i$ is equipped with this topology. In addition, ${\bf \Gamma}$ has the product topology induced by these $w^*$-topologies; that is, $\underline{\gamma}_n$ converges to $\underline{\gamma}$ in ${\bf \Gamma}$ with respect to the product topology if and only if $\gamma_{n}^i \rightharpoonup^* \gamma^i$ for all $i=1,\ldots,N$. In this case we write $\underline{\gamma}_n \rightharpoonup^* \underline{\gamma}$. Note that ${\bf \Gamma}$ is sequentially relatively compact under this topology.

To this end, for any $L>0$, we define
\begin{align}
{\bf \Gamma}_L &\coloneqq \biggl\{\underline{\gamma} \in {\bf \Gamma}: J(\underline{\gamma}) < J^* + L  \biggr\} \nonumber \\
\intertext{and}
S_{L} &\coloneqq \biggl\{\lambda \in {\cal P}\left(\Omega_0 \times \prod_{i=1}^N (\mathbb{Y}^i \times \mathbb{U}^i)\right): \nonumber \\
&\phantom{xxxxxxxxxxxxxxxxxx}\lambda(d\omega_0,d{\bf y},d{\bf u}) = \mu(d\omega_0) \prod_{i=1}^N \gamma^i(du^i|y^i) \mu^i(dy^i), \underline{\gamma} \in {\bf \Gamma}_L \biggr\}. \nonumber
\end{align}
For each $i=1,\ldots,N$, we define $S_L^i \coloneqq \Proj_{\mathbb{Y}^i \times \mathbb{U}^i}(S_L)$. In order to prove that ${\bf \Gamma}_L$ is closed with respect to the $w^*$-topology, we should impose the below assumption. 

\begin{assumption}\label{top-as2}
For some $L>0$, $S_L^i$ is tight for $i=1,\ldots,N$.
\end{assumption}

We provide give several conditions that imply Assumption~\ref{top-as2}.

\begin{theorem}[{\cite[Theorem 4]{SaldiArXiv2017}}]
Suppose \textit{either }of the following conditions hold:
\begin{itemize}
\item [(i)] $\mathbb{U}^i$ is compact for all $i$.
\item [(ii)] For non-compact case, we assume
\begin{itemize}
\item [(a)] The cost function $c$ satisfies the following condition: for each $j$, for any $M \geq 0$, and for any compact $K \subset \Omega_0 \times \prod_{k=1}^N \mathbb{Y}^k \times \prod_{k=1}^{j-1} \mathbb{U}^k$, there exists a compact $L \subset \mathbb{U}^j$ such that 
$$
\inf_{K \times L^c \times \prod_{k=j+1}^N \mathbb{U}^k} c(\omega_0,{\bf y},{\bf u}) \geq M. 
$$
\item [(b)] For all $j$, $g^j > 0$ and $g^j$ is lower semi-continuous.
\end{itemize}
\end{itemize}
Then, Assumption~\ref{top-as2} holds.
\end{theorem}

%
%
The following theorem establishes the existence of optimal team decision rule using $w^*$-topology. 

\begin{theorem}[{\cite[Theorem 6]{SaldiArXiv2017}}]\label{top-main2}
Suppose Assumptions~\ref{top-as1} and \ref{top-as2} hold. Moreover, we assume that $c$ is lower semi-continuous and measurement channels $g^i(dy^i|\omega_0)$ are continuous with respect to the total variation distance. Then, there exists $\underline{\gamma}^{*} \in {\bf \Gamma}_L$ which is optimal. 
\end{theorem}

\smallskip

\proofsketch
Note that when the cost function $c$ is compactly supported continuous function, by Stone-Weierstrass Theorem \cite[Lemma 6.1]{Lan93}, $c$ can be uniformly approximated by functions of the form
\begin{align}
\sum_{j=1}^k r_j \prod_{i=1}^N f_{j,i} g_{j,i}, \label{newww-eq}
\end{align}
where $r_j \in C_0(\Omega_0)$, $f_{j,i} \in C_0(\mathbb{Y}^i)$, and $g_{j,i} \in C_0(\mathbb{U}^i)$ for each $j=1,\ldots,k$ and $i=1,\ldots,N$, are compactly supported. One can prove that the cost $J$ is continuous with respect to $w^*$-topology if $c$ is of the form in (\ref{newww-eq}) \cite[Theorem 5]{SaldiArXiv2017}. Since any compactly supported function can be uniformly approximated by such functions, $J$ is also continuous if $c$ is compactly supported. Note that if $c$ is lower semi-continuous, then $c$ can be approximated pointwise from below by compactly supported functions (see the proof of \cite[Proposition 1.4.18]{HeLa03}). Hence, if $c$ is lower semi-continuous, then $J$ is also lower semi-continuous with respect to $w^*$-topology by monotone convergence theorem.

Now, we prove that ${\bf \Gamma}_L$ is closed under $w^*$-topology, which will complete the proof since ${\bf \Gamma}_L$ is relatively compact. To this end, let
\begin{align}
\underline{\gamma}_{n} \rightharpoonup^* \underline{\gamma}, \nonumber
\end{align}
for some $\underline{\gamma} \in \prod_{i=1}^N {\cal L}_{\infty}\bigl(\mu^i,{\cal M}(\mathbb{U}^i)\bigr)$ (recall that ${\bf \Gamma}$ is not closed with respect to $w^*$-topology). If $\gamma^i \in {\cal P}(\mathbb{U}^i)$ $\mu^i$-a.e. for all $i$, then $\underline{\gamma}$ is in ${\bf \Gamma}_L$, and so, ${\bf \Gamma}_L$ is closed. 

Fix any $i$. Note that the sequence $\{\gamma^i_{n}\otimes\mu^i\}$ is tight as it is a subset of $S_L^i$. Thus, there exists a further subsequence, denoted for simplicity by $\{\gamma^i_{l}\otimes\mu^i\}$, that converges weakly to some $\lambda \in {\cal P}(\mathbb{Y}^i \times \mathbb{U}^i)$.
This implies that $\gamma^i\otimes\mu^i = \lambda$, and so, $\gamma^i\otimes\mu^i(\mathbb{Y}^i \times \mathbb{U}^i) = 1$. Hence, $\gamma^i \in {\cal P}(\mathbb{U}^i)$ $\mu^i$-a.e. Thus, ${\bf \Gamma}_L$ is closed. This completes the proof. 
\QEDA

\syr{
\begin{remark}
Kulkarni \cite{Kul16,Kul18} uses weak topology in a similar context. In these works, \cite[Lemma 4.9]{Kul16} and \cite[Lemma A.1]{Kul18} can be exploited to prove the existence of optimal policies. Indeed, if one can prove that actions are square integrable under any strategy, \syr{the set of actions is weakly relatively compact}, and the cost function is lower semi-continuous, then the existence of optimal policies can be deduced by Weierstrass Extreme Value Theorem since \cite[Lemma 4.9]{Kul16} and \cite[Lemma A.1]{Kul18} guarantees that the limiting actions preserve the information structure of the problem.
\end{remark}
}

%

%
\sy{
\subsection{Exchangeability, Infinite Products of Individual Policies and Optimality of Symmetric Randomized Policies for Mean-Field Teams}

As a final example on the utility of placing a product topology on individual policies, we consider stochastic team problems with infinitely many decision makers. Such problems have seen a significant activity in the context of mean field theory 
\cite{CainesMeanField2,CainesMeanField3,LyonsMeanField} (see also more recent papers \cite{fischer2017connection, biswas2015mean, arapostathis2017solutions,lacker2016general}) and in mean-field team problems \cite{huang2012social, wang2017social}\cite{arabneydi2015team} \cite{mahajan2013static} \cite{sanjariyukseltac}. In the context of mean-field team problems \cite{mahajan2013static} and \cite{sanjariyukseltac} have shown that, under sufficient convexity conditions, a sequence of optimal policies for teams with $N$ number of decision makers as $N \to \infty$ converges to a team optimal policy for the static team with countably infinite number of decision makers, where the latter establishes the optimality of symmetric (i.e., identical for each DM) policies as well as existence of optimal team policies for both finite and infinite DM setups. In the following, we demonstrate how the control topology approach can be utilized to relax such strong convexity conditions.


Consider the following:

\begin{itemize}
\item[ ]{\bf\text{Problem}} {\bf ($\mathcal{P}_{N}$)}:
Let ${\cal N}=\{1,\dots, N\}$. Let $\underline{\gamma}_{N} = (\gamma^1, \cdots, \gamma^N)$ and ${\Gamma}_{N} = \prod_{i=1}^{N} \Gamma^i$. Define an expected cost function of $\underline{\gamma}_{N}$ as
\begin{equation}\label{eq:1.1}
J_{N}(\underline{\gamma}_{N}) = \mathbb{E}^{\underline{\gamma}_{N}}[c(\omega_{0},\underline{u}_{N})]:= \mathbb{E}[c(\omega_{0},\gamma^1(y^1),\cdots,\gamma^N(y^N))],
\end{equation}
for some Borel measurable cost function $c: \Omega_{0} \times \prod_{k=1}^{N} \mathbb{U}^k \to \mathbb{R}_{+}$. We define $\omega_{0}$ as the $\Omega_{0}$-valued cost function relevant exogenous random variable as $\omega_{0}:(\Omega,\mathcal{F}, \mathbb{P}) \to (\Omega_{0},\mathcal{F}_{0})$, where $\Omega_{0}$ is a Borel space with its Borel $\sigma$-field $\mathcal{F}_{0}$. Here, we have the notation $\underline {u}_{N}:=\{u^i, i \in {\cal N}\}$.
\end{itemize}


\begin{itemize}
\item[ ]{\bf\text{Problem}} {\bf ($\mathcal{P}_{\infty}$)}:
Consider a stochastic team with countably infinite number of decision makers, that is, $\mathcal{N}=\mathbb{N}$. Let ${\Gamma}=\prod_{i \in \mathbb{N}} \Gamma^{i}$ and $\underline{\gamma}=(\gamma^{1},\gamma^{2},\dots)$. Let $c:\Omega_{0} \times \mathbb{U}\times \mathbb{U} \rightarrow \mathbb{R}_{+}$. Define the expected cost of $\underline{\gamma}$ as
\begin{equation}\label{eq:2.5.5}
J(\underline{\gamma})=\limsup\limits_{N\rightarrow \infty} \frac{1}{N} \mathbb{E}^{\underline{\gamma}}\bigg[\sum_{i=1}^{N}c\bigg(\omega_{0},u^{i},\frac{1}{N}\sum_{p=1}^{{N}}u^{p}\bigg)\bigg].
\end{equation}
With slight abuse of notation, we use the same notation for the cost function $c$ as in \eqref{eq:1.1}.
\end{itemize}

\begin{assumption}\label{assump:ind}
Assume for any DM$^{i}$, there exists a probability measure $Q^{i}$ on $\mathbb{Y}^{i}$ and a function $f^{i}$ such that for all Borel sets $S=S^1 \times \dots \times S^N$, we have
\begin{align}
&\mathbb{P}((y^{1},\dots,y^{N}) \in S \big|\omega_{0})=\prod_{i=1}^{N}\int_{S^i} f^{i}(y^{i}, \omega_{0},y^{1},\dots, y^{i-1})Q^{i}(dy^{i})\label{eq:abscon}.
\end{align}
\end{assumption}

As observed in Section \ref{EquivIS}, Assumption \ref{assump:ind} allows us to reduce the problem as a static team problem where the observation of each DM is independent of observations of other DMs and also independent of $\omega_{0}$. Hence, under Assumption \ref{assump:ind}, we can focus on each DM$^{i}$ separately and identify $\Gamma^{i}$ via the set of probability measures
 \begin{align}
 &\Theta^{i}:=\bigg\{P \in \mathcal{P}(\mathbb{U}^{i}\times \mathbb{Y}^{i}) \bigg| P(du^{i},dy^{i})=1_{\{\gamma^{i}(y^{i})\in du^{i}\}}Q^{i}(dy^{i}), \gamma^{i} \in \Gamma^{i}\bigg\}.
 \end{align}
As noted earlier, the above set is the set of extreme points of the set of probability measures on $(\mathbb{U}^{i}\times \mathbb{Y}^{i})$ with fixed marginals $Q^{i}$ on $\mathbb{Y}^{i}$. Hence it inherits Borel measurability and topological properties of  that Borel measurable set \cite{BorkarRealization}. As before, we define convergence on policies as $\gamma^{i}_{n} \to \gamma^{i}$ iff $1_{\{\gamma^{i}_{n}(y^{i})\in du^{i}\}}Q^{i}(dy^{i}) \to 1_{\{\gamma^{i}(y^{i})\in du^{i}\}}Q^{i}(dy^{i})$ (in the weak convergence topology) as $n\to \infty$ for each DM. We will also allow for randomized (relaxed) policies. Accordingly, each individual control policy $\gamma^i \in \Gamma^i$ is an element in the set of probability measures $\mathcal{P}(\mathbb{U}^{i}\times \mathbb{Y}^{i})$ with a fixed marginal, $Q^i$, on $\mathbb{Y}^i$.
    

Now that we have a standard Borel space formulation for policies, we can define the set of probability measures on policies with product topology on $\Gamma_{N}=\prod_{i=1}^{N}\Gamma^{i}$. We define the following set of Borel probability measures on admissible relaxed policies $\Gamma_{N}$ as follows:
\begin{align}
&L^{N}:=\mathcal{P}(\Gamma_{N})\label{eq:L},
\end{align}
where Borel $\sigma$-field $\mathcal{B}(\Gamma^{i})$ is induced by the topology defined above.

Recall the definition of {\it exchangeability} for random variables.

\begin{definition}
Random variables $x^{1},x^{2},\dots,x^{N}$ defined on a common probability space are $N$-\it{exchangeable} if for any permutation $\sigma$ of the set $\{1,\dots,N\}$, 
\begin{align}
&{P}\bigg(x^{\sigma(1)} \in A^{1},x^{\sigma(2)} \in A^{2},\dots,x^{\sigma(N)}\in A^{N}\bigg)={P}\bigg(x^{1} \in A^{1},x^{2} \in A^{2},\dots,x^{N}\in A^{N}\bigg) \nonumber
\end{align}
for any measurable $\{A^{1},\dots, A^N\}$, and $(x^{1},x^{2},\dots)$ is {\it infinitely-exchangeable} if it is $N$-\it{exchangeable} for all $N\in \mathbb{N}$.
\end{definition}

Now, we define the set of exchangeable probability measures on policies as:
\begin{align}
L_{\text{EX}}^{N}:=\bigg\{&P_{\pi} \in L^{N}\bigg{|}\text{for all}~A_{i} \in \mathcal{B}(\Gamma^{i})~\text{and for all}~\sigma \in S_{N}:\nonumber\\
&P_{\pi}(\gamma^{1} \in A_{1},\dots,\gamma^{N}\in A_{N})=P_{\pi}(\gamma^{\sigma(1)} \in A_{1},\dots,\gamma^{\sigma(N)}\in A_{N})\bigg\}, \nonumber
\end{align}
where $S_{N}$ is the space of permutations of $\{1,\dots,N\}$. We note that $L_{\text{EX}}^{N}$ is a convex subset of $L^{N}$. Define the set of probability measures on policies induced by a common randomness as: 
\begin{eqnarray}
L_{\text{CO}}^{N}&&:=\bigg\{ P_{\pi} \in L^{N}\bigg{|}\text{for all}~A_{i} \in \mathcal{B}(\Gamma^{i}): \nonumber\\
&&P_{\pi}(\gamma^{1} \in A_{1},\dots,\gamma^{N}\in A_{N})=\int_{z\in [0,1]}\prod_{i=1}^{N}P_{\pi}^{i}(\gamma^{i}\in A_{i}|z)\eta(dz),\qquad \eta \in \mathcal{P}([0, 1])\bigg\}, \nonumber
\end{eqnarray} 
where $\eta$ is the distribution of common, but independent (from intrinsic exogenous system variables), randomness. Note that conditioned on $z$, policies are independent. We also define the set $L_{\text{CO,SYM}}^{N}$ as the set of identical probability measures on  policies induced by a common randomness:
\begin{eqnarray}
L_{\text{CO,SYM}}^{N}&&:=\bigg\{ P_{\pi} \in L^{N}\bigg{|}\text{for all}~A_{i} \in \mathcal{B}(\Gamma^{i}): \nonumber\\
&& P_{\pi}(\gamma^{1} \in A_{1},\dots,\gamma^{N}\in A_{N})=\int_{z\in [0,1]}\prod_{i=1}^{N}P_{\pi}(\gamma^{i} \in A_{i}|z)\eta(dz), \qquad \eta \in \mathcal{P}([0, 1])\bigg\}, \nonumber
\end{eqnarray} 
where we drop the index $i$ in $P_{\pi}$ to indicate that the independent randomization is identical through DMs. Also, define the set of probability measures on policies with only private independent randomness as:
\begin{eqnarray}
L_{\text{PR}}^{N}&:=&\bigg\{P_{\pi} \in L^{N}\bigg{|}\text{for all}~A_{i} \in \mathcal{B}(\Gamma^{i}): P_{\pi}(\gamma^{1} \in A_{1},\dots, \gamma^{N} \in A_{N})=\prod_{i=1}^{N}P_{\pi}^{i}(\gamma^{i}\in A_{i})\bigg\}\nonumber.
\end{eqnarray} 
Finally, define the set of probability measures on policies with identical and independent randomness:
\begin{align}
&L_{\text{PR,SYM}}^{N} \nonumber \\
&\phantom{xx}:=\bigg\{P_{\pi} \in L^{N}\bigg{|}\text{for all}~A_{i} \in \mathcal{B}(\Gamma^{i}): P_{\pi}(\gamma^{1} \in A_{1},\dots, \gamma^{N} \in A_{N})=\prod_{i=1}^{N}P_{\pi}(\gamma^{i}\in A_{i})\bigg\}\nonumber.
\end{align} 

For a team with countably infinite number of decision makers, we define sets of probability measures $L, L_{\text{EX}}, L_{\text{CO}}, L_{\text{CO,SYM}}, L_{\text{PR}}, L_{\text{PR,SYM}}$ similarly using Ionescu Tulcea
extension theorem by iteratively adding new coordinates for our probability measure (see e.g., \cite{InfiniteDimensionalAnalysis, HernandezLermaMCP}). We define the set of probability measures $L$ on the infinite product Borel spaces $\Gamma=\prod_{i\in \mathbb{N}}\Gamma^{i}$ as: 
\begin{eqnarray}
L:=\mathcal{P}(\Gamma)\label{eq:Linf}.
\end{eqnarray} 
Now, we define the set of infinitely exchangeable probability measures on policies as:
\begin{eqnarray}
L_{\text{EX}}&:=&\bigg\{P_{\pi} \in L\bigg{|}\text{for all}~A_{i} \in \mathcal{B}(\Gamma^{i})~\text{and for all $N\in\mathbb{N}$, and for all $\sigma \in S_{N}$:} \nonumber\\
&&P_{\pi}(\gamma^{1} \in A_{1},\dots,\gamma^{N}\in A_{N})=P_{\pi}(\gamma^{\sigma(1)} \in A_{1},\dots,\gamma^{\sigma(N)}\in A_{N})\bigg\}, \nonumber
\end{eqnarray} 
and we define
\begin{eqnarray}
L_{\text{CO}}&:=&\bigg\{P_{\pi} \in L\bigg{|}\text{for all}~A_{i} \in \mathcal{B}(\Gamma^{i}): \nonumber\\
&&P_{\pi}(\gamma^{1} \in A_{1},\gamma^{2} \in A_{2},\dots)=\int_{z\in [0,1]}\prod_{i\in \mathbb{N}}P_{\pi}^{i}(\gamma^{i}\in A_{i}|z)\eta(dz), \qquad \eta \in \mathcal{P}([0, 1])\bigg\}. \nonumber
\end{eqnarray} 
 Note that $L_{\text{CO}}$ is a convex subset of $L$ and its extreme points are in the set of probability measures on policies with private independent randomness:
\begin{eqnarray}
L_{\text{PR}}&:=\bigg\{P_{\pi} \in L\bigg{|}\text{for all}~A_{i} \in \mathcal{B}(\Gamma^{i}): P_{\pi}(\gamma^{1} \in A_{1},\gamma^{2} \in A_{2},\dots)=\prod_{i\in \mathbb{N}}P_{\pi}^{i}(\gamma^{i}\in A_{i})\bigg\}\nonumber.
\end{eqnarray} 
Also, we define
\begin{eqnarray}
L_{\text{CO,SYM}}&&:=\bigg\{ P_{\pi} \in L\bigg{|}\text{for all}~A_{i} \in \mathcal{B}(\Gamma^{i}): \nonumber\\
&& P_{\pi}(\gamma^{1} \in A_{1},\gamma^{2}\in A_{2}, \dots)=\int_{z\in [0,1]}\prod_{i\in \mathbb{N}}P_{\pi}(\gamma^{i} \in A_{i}|z)\eta(dz), \qquad \eta \in \mathcal{P}([0, 1])\bigg\}, \nonumber
\end{eqnarray} 
and we define
\begin{eqnarray}
L_{\text{PR,SYM}}&:=\bigg\{P_{\pi} \in L\bigg{|}\text{for all}~A_{i} \in \mathcal{B}(\Gamma^{i}): P_{\pi}(\gamma^{1} \in A_{1},\gamma^{2}\in A_{2}, \dots)=\prod_{i\in \mathbb{N}}P_{\pi}(\gamma^{i} \in A_{i})\bigg\}. \nonumber
\end{eqnarray} 

\begin{theorem}\label{the:defin}\cite{sanjari2020optimality}
Suppose that Assumption \ref{assump:ind} holds. Then, any $P_{\pi}\in L_{\text{EX}}$ satisfying the following condition:
\begin{itemize}
\item For every $i \in \mathbb{N}$, $\mathbb{E}(\phi_{i}(u^{i}))\leq K$ for some finite $K$, where $\phi_{i}:\mathbb{U}^{i} \to \mathbb{R}_{+}$ is a lower semi-continuous moment function.
\end{itemize}
 is in $L_{\text{CO,SYM}}$, i.e., for any $P_{\pi}\in L_{\text{EX}}$ satisfying the above moment condition, there exists a $[0,1]$-valued random variable $z \sim \eta$ such that for any $A_{i} \in \mathcal{B}(\Gamma^{i})$
\begin{align}
&P_{\pi}(\gamma^{1} \in A_{1},\gamma^{2} \in A_{2},\dots)=\int_{z\in [0,1]}\prod_{i\in \mathbb{N}}P_{\pi}(\gamma^{i}\in A_{i}|z)\eta(dz)
\end{align} 
\end{theorem}

\begin{assumption}\label{assump:exccost}
The cost function $c$ in problem $({\cal P}_N)$ is exchangeable with respect to actions for all $\omega_{0}$, i.e., for any permutation $\sigma$ of $\{1,\dots,N\}$ $c(\omega_{0}, u^{1},\dots, u^{N})=c(\omega_{0}, u^{\sigma(1)},\dots, u^{\sigma(N)})$ for all $\omega_{0}$.
\end{assumption}

\begin{lemma}\label{lem:exc}
For a fixed $N$, consider an $N$-DM static team. Assume $\bar{L}^{N}$ is an arbitrary convex subset of $L^{N}$. Under Assumption \ref{assump:ind} and Assumption \ref{assump:exccost}, if observations of DMs are exchangeable conditioned on $\omega_{0}$, then 
\begin{eqnarray}\label{eq:3.5}
\inf \limits_{P_{\pi} \in \bar{L}^{N}}\int P_{\pi}(d\underline{\gamma})\mu^{N}(d\omega_{0},d\underline{y})c^{N}(\underline{\gamma},\underline{y}, \omega_{0})=\inf\limits_{P_{\pi} \in \bar{L}^{N} \cap L_{\text{EX}}^{N}}\int P_{\pi}(d\underline{\gamma})\mu^{N}(d\omega_{0},d\underline{y})c^{N}(\underline{\gamma},\underline{y}, \omega_{0}), \nonumber
\end{eqnarray} 
where $c^{N}(\underline{\gamma},\underline{y}, \omega_{0}) \coloneqq c(\omega_0,\gamma^1(y^1),\ldots,\gamma^N(y^N)).$
\end{lemma}


\begin{assumption}\label{assump:cont}
The cost function $c:\Omega_{0} \times \mathbb{U} \times \mathbb{U} \rightarrow \mathbb{R}_{+}$ in problem $({\cal P}_{\infty})$ is continuous in its second and third arguments for all $\omega_{0}$.
\end{assumption}

Under mild conditions, we can show that the optimal expected cost function induced by $L^{N}_{\text{EX}}$ and $L_{\text{EX}}$ are equal as $N$ goes to infinity. Hence, by Lemma \ref{lem:exc}, under symmetry, this allows us to show that  without loss of global optimality, optimal policies of static mean-field teams with countably infinite number of DMs can be considered to be an infinitely exchangeable type. 

\begin{lemma}\label{lem:findef}
Suppose that Assumption \ref{assump:ind} and Assumption \ref{assump:cont} hold. Assume further that $\mathbb{U}$ is compact  and the cost function is bounded. If observations of DMs are i.i.d. random vectors conditioned on $\omega_{0}$, then 
\begin{align}
&\limsup\limits_{N \to \infty}\inf\limits_{P_{\pi}^{N} \in L^{N}_{\text{EX}}}\int P_{\pi}^{N}(d\underline{\gamma})\mu^{N}(d\omega_{0},d\underline{y}) c^{N}(\underline{\gamma},\underline{y}, \omega_{0})\nonumber\\
&=\limsup\limits_{N \to \infty}\inf\limits_{P_{\pi} \in L_{\text{EX}}}\int P_{\pi, N}(d\underline{\gamma})\mu^{N}(d\omega_{0},d\underline{y})c^{N}(\underline{\gamma},\underline{y}, \omega_{0})\label{eq:ert},
\end{align}
where $c^{N}(\underline{\gamma},\underline{y}, \omega_{0}) \coloneqq \frac{1}{N} \sum_{i=1}^N c\left(\omega_0,\gamma^i(y^i),\frac{1}{N} \sum_{p=1}^N \gamma^p(y^p)\right)$, $P_{\pi,N}$ is the marginal of the $P_{\pi}\in L_{\text{EX}}$ to the first $N$ components and $\mu^{N}$ is the marginal of the fixed probability measure on $(\omega_{0}, y^{1},y^{2},\dots)$ to the first $N+1$ components.
\end{lemma}

We now establish an existence of a randomized optimal policy for ($\mathcal{P}_{\infty}$), which is symmetric. 

\begin{theorem}\label{the:2}\cite{sanjari2020optimality}
Consider a static team problem ($\mathcal{P}_{\infty}$) where Assumption \ref{assump:ind} and Assumption \ref{assump:cont} hold. Assume further that $\mathbb{U}$ is compact. 
If observations of DMs are i.i.d. random vectors conditioned on $\omega_{0}$, then there exists a randomized optimal policy $P^{*}_{\pi}$ for {($\mathcal{P}_{\infty}$)} which is in $L_{\text{PR,SYM}}$:
\begin{align}
&\min \limits_{P_{\pi} \in L_{\text{PR,SYM}}} \limsup\limits_{N \to \infty} \int P_{\pi, N}(d\underline{\gamma})\mu
^{N}(d\omega_{0},d\underline{y})c^{N}(\underline{\gamma},\underline{y}, \omega_{0})\\
&=\inf\limits_{P_{\pi} \in L_{\text{PR}}} \limsup\limits_{N \to \infty} \int P_{\pi, N}(d\underline{\gamma})\mu
^{N}(d\omega_{0},d\underline{y})c^{N}(\underline{\gamma},\underline{y}, \omega_{0}).
\end{align}
\end{theorem}

Following from Lemma \ref{lem:exc}, Lemma \ref{lem:findef}, and our analysis in the proof of Theorem \ref{the:2}, thanks to Theorem \ref{the:defin}, we can show that without losing global optimality, optimal policies for mean-field teams can be considered to be symmetric and privately randomized ($L_{\text{PR,SYM}}$).

We note that if one also has convexity in the cost as well as action sets $\mathbb{U}^i$, then one can also establish that for every finite $N$, the optimal policies are symmetric and deterministic, but in the infinite limit, randomization may be required \cite{sanjariyukseltac}. Similar results also hold for dynamic team problems \cite{sanjari2019optimal}. We emphasize that a {\it strategic measures} approach would not be feasible for arriving at this solution since exchangeability in the actions is not sufficient to ensure that the {\it dominating} random variable (in the de Finetti representation) is independent of the intrinsic randomness in the system. 
 }
 
\subsection{Extended Weak Convergence, Topology of Information, and Adapted Topologies}
\sy{
A versatile topology, which has evidently been used in a variety of contexts in stochastic analysis, is essentially given by the following convergence notion. Consider a stochastic process converging to another one in the following sense: all finite dimensional marginals converge weakly and the conditional kernels on the future random variables given the past converge weakly as well when conditional kernels are viewed as measure valued stochastic processes. 

This notion has been applied for different problems: Aldous has termed it {\it extended weak convergence} \cite{aldous1981weak} and Hellwig has named it \textit{the information topology} \cite{hellwig1996sequential}; these have recently shown to be equivalent in discrete-time \cite[Theorem 1.1]{backhoff2019all}. In addition to the conditional independence preservation \cite{hellwig1996sequential} \cite{gupta2020existence}, applications in robust stochastic control also follow from this discussion, with details being very context-specific: \cite{kara2020robustness}, \cite{Lan81}, \cite{bayraktar2020continuity,julio2020adapted}.

In relation to our context, such a convergence notion again requires strong continuity conditions under all admissible policies. In the same spirit of the discussion in Section \ref{tradeOffContComp}, though continuity and closedness hold under convergence with this notion, compactness will require more restrictive conditions. 

With this approach, the recent work \cite{gupta2020existence} (see also \cite{barbie2014topology}) has established existence results for a setup where either the measurements are countable or there is a common information among decision makers which is countable space-valued with the private information satisfying an absolute continuity condition. As noted earlier in the paper, static reduction applies in both such setups and the results presented in this paper (notably Theorem \ref{existenceRelaxed3}) generalize those reported in \cite{gupta2020existence}. We note also that the use of $w$-$s$ topology in Theorem \ref{existenceRelaxed3} significantly relaxes the requirements of continuity.}

\section{Revisiting (and avoiding a subtle potential error in) Relaxed Policies for Partially Observed Stochastic Control}\label{relaxedPOMDPs}

In this section, we will revisit the concept of {\it relaxed} control policies for classical stochastic control problems, with a further relaxation known as {\it wide sense admissible policies} introduced by Fleming and Pardoux \cite{FlPa82} and prominently used to establish the existence of optimal solutions for partially observed stochastic control problems. Borkar \cite{Bor00,Bor07,Bor03} (see also \cite{BoBu04}) has utilized these policies for a coupling/simulation method to arrive at optimality results for average cost partially observed stochastic control problems. 

The main goal of this section is to show that if wide sense admissible control policies are not defined in a correct form, this can lead to a significant error in reasoning: the controllers may be allowed to have access to have information that they should not. 

Relaxed control policies are extremely useful concepts as shown in Section~\ref{setsOfStrategicMeasures} and allows one to use topologies on the sets of probability measures to study existence, optimality, and structural results. A key aspect of such relaxations is that, the relaxation should not allow for optimal expected cost values to be improved; they should only be means to facilitate stochastic analysis. Our goal here is, building on the insights developed in Section~\ref{nonsignal} (in particular by the analysis on the CHSH team \cite{CHSH1969} reviewed in Theorem~\ref{relation-thm}), to highlight a subtlety which may lead to incorrect conclusions if the relaxation is not cautiously constructed. Accordingly, we first present a brief overview of relaxed control policies in continuous-time or discrete-time stochastic control. 

Consider a continuous-time Markov decision process $\{x_t\}$ on an Euclidean space $\mathbb{R}^{N}$, controlled by a control process $\{u_t\}$ taking values in a convex and compact Borel action space $\mathbb{U} \subset \mathbb{R}^L$, and with an associated observation process $\{y_t\}$ taking values in $\mathbb{R}^M$, where $0 \leq t \leq T$. The evolution of $\{x_t,y_t\}$ is given by stochastic differential equations
\begin{align}
dx_t &= b(x_t,y_t,u_t) dt + \sigma(x_t,y_t) dW_t, \label{sde1} \\
dy_t &= h(x_t) dt + dB_t. \label{sde2}
\end{align}
Here, $W$ and $B$ are independent standard Wiener processes with values in $\mathbb{R}^D$ and $\mathbb{R}^M$, respectively (hence, $\sigma$ is a  $N \times D$-matrix). The objective is to minimize the following cost function
\begin{align}
E \bigg[ \int_0^T F(x_t,u_t) dt + G(x_T) \bigg], \nonumber
\end{align}
where $F: \mathbb{R}^N \times \mathbb{U} \rightarrow [0,\infty)$ and $G:\mathbb{R}^N \rightarrow [0,\infty)$. In the literature, it is customary to require that control process $\{u_t\}$ be adapted to the filtration generated by the observation process $\{y_t\}$; that is, for each $t \in [0,T]$, $u_t$ is $\sigma\left(y_s, 0 \leq s \leq t \right)$-measurable. We will call such policies (strict-sense or precise) admissible policies. In \cite{FlPa82}, Fleming and Pardoux introduced another class of policies which they named to be \textit{wide-sense admissible policies}. Using this relaxed class of policies, they study the existence of optimal policies to the above problem.

To define wide sense admissible policies, we first reproduce the above processes on a canonical probability space
$$\Omega = \Omega_0 \times \Omega_1 \times \Omega_2 \times \Omega_3,$$
where $\Omega_0$, $\Omega_1$, and $\Omega_2$ are $C([0,T];\mathbb{R}^m)$ with $m = D, N, M$, respectively, and $\Omega_3 = L^2([0,T];\mathbb{U})$. We identify the Wiener process $W$, state process $x$, observation process $y$, and action process $u$ as follows. If $\omega = (W,x,y,u) \in \Omega$, then
\begin{align}
w(t) = (w_t(w),x_t(w),y_t(w),u_t(w)), \, \, 0 \leq t \leq T. \nonumber
\end{align}
Here, $\Omega_1$, $\Omega_2$, and $\Omega_3$ are endowed with usual sup-norm topology and $\Omega_3$ is endowed with weak topology. Let $\Omega^2 = \Omega_2 \times \Omega_3$ and define
\begin{align}
{\cal F}_t(y) = \sigma\left(y_s, 0 \leq s \leq t\right), \, \,  {\cal F}_t(u) = \sigma\left(u_s, 0 \leq s \leq t\right), \,\, \G_t^2 = {\cal F}_t(y) \times {\cal F}_t(u). \nonumber
\end{align}
Here, $\G_t^2$ is the product $\sigma$-field generated by ${\cal F}_t(y)$ and ${\cal F}_t(u)$. Note that $\{\G_t^2\}$ is a filtration on $\Omega^2$. 

We can now define the class of wide-sense admissible policies. A \textit{wide-sense admissible policy} $\pi$ is a probability measure on $(\Omega^2,\G_T^2)$ such that $y$ is a $(\pi,\{G_t^2\})$ Wiener process. Note that this definition requires that the projection $(y,u) \rightarrow y$ maps $\pi$ onto Wiener measure and $\{u_s, 0 \leq s \leq t\}$ is independent of the increment $y_r - y_t$ for all $t \leq r \leq T$. The latter condition states that actions upto time $t$ is independent of the observations after time $t$ given past observations and actions. In other words, instead of saying that actions should be dependent on current and past observations, this condition states that actions should be independent of future observations given past observations and actions.

Given a distribution $\mu$ of $x_0$, each wide-sense admissible policy induces a joint measure $P^{\pi}$ of $(w,x,y,u)$ as follows. Let $P^{y,u}$ denote the conditional probability of $(w,x)$ given $(y,u)$ induced by stochastic differential equations (\ref{sde1}),(\ref{sde2}), which is independent of the policy as the action process $u$ is given. Then define $P_0^{\pi}$ as follows
\begin{align}
P_0^{\pi}(dw,dx,dy,du) = P^{y,u}(dw,dx) \otimes \pi(dy,du). \nonumber
\end{align}
Let us define the following density function
\begin{align}
Z_T = \exp \biggl[ \int_0^T h(x_s) dy_s - \frac{1}{2} \int_0^T |h(x_s)|^2 ds \biggr]. \nonumber
\end{align}
Then, $P^{\pi}$ is given by
\begin{align}
\frac{dP^{\pi}}{dP_0^{\pi}} = Z_T. \nonumber
\end{align}
\syr{Here, we apriori assume that this is integrable under the new measure,} we are indeed applying Girsanov's transformation. Let $B_t = y_t - \int_0^t h(x_s) ds$. Then, under $P^{\pi}$, $W$ and $B$ are independent standard Wiener processes and stochastic differential equations (\ref{sde1}) and (\ref{sde2}) hold. Hence, the cost function of the wide-sense admissible policy $\pi$ is given by
\begin{align}
J(\pi) = E^{\pi} \biggl[ \int_0^T F(x_t,u_t) dt + G(x_T) \biggr], \nonumber
\end{align}
where $E^{\pi}$ denotes the expectation with respect to $P^{\pi}$. 

Note that in this setup, a policy $\pi$ is admissible in the classical sense if $\pi$ is wide-sense admissible and there exists ${\cal F}_T(y)/{\cal F}_T(u)$-measurable $f:\Omega_2 \rightarrow \Omega_3$ such that $\pi(dy,du) = \delta_{f(y)}(du) \otimes w(dy)$, where $w$ is a Wiener measure. \sy{Likewise, we can say that the policy is relaxed if $f: y \mapsto {\cal P}(\mathbb{U})$ is probability measure valued. We refer the reader to \cite{FlPa82} for a more explicit construction of wide-sense admissible policies in the continuous time setup.}

In \cite[Theorem 7.2]{FlPa82}, Fleming and Pardoux proved the existence of optimal wide-sense policy by converting the original problem to a fully-observed continuous-time Markov decision process on the belief-space $P(\mathbb{R}^N)$; that is, the state of the belief-space MDP is the conditional distribution of the state $x_t$ given the past observations and actions $\G_t^2$. Then, by adding some mild conditions on the stage-wise cost functions $F$ and $G$, they also proved that the infimum achieved by classical admissible policies is the same as the infimum achieved by wide-sense admissible policies \cite[Theorem 6.1]{FlPa82}. Therefore, without loss of generality, one can work with wide-sense admissible policies instead of classical policies in order to further analyze such problems. 

\begin{remark}
It may be important to note that Bismut \cite{bismut1982partially} arrived at further existence results, through an approach which avoids separation (and the construction of a belief-MDP), in discrete-time a similar approach is given in \cite[Section 5.4.2]{YukselWitsenStandardArXiv}.
\end{remark}

\subsection{Discrete-time Case}

Inspired by the work of Fleming and Pardoux \cite{FlPa82}, Borkar introduced wide-sense control policies to study discrete-time partially-observed finite state-observation Markov decision processes with average cost criterion in \cite[p.675 ]{Bor00}, \cite[item 1]{Bor07}. Using coupling methods, Borkar proved that fully-observed belief-space Markov decision process under wide-sense admissible policies admits the solution of average cost optimality equation and any stationary policy for the belief-space MDP that solves this equation is optimal. Borkar later extended this result to continuous state-observation MDPs in \cite{Bor03,BoBu04}. While it is evident that the discussion in \cite[p.675]{Bor00}, \cite[item 1]{Bor07} is just an oversight (since this was corrected later, as we will discuss below), nonetheless we wish to note that if a reader applies this as written, this may lead to a consequential error, as we note in the following and as this was not explicitly noted.


To this end, we first review the construction of Borkar, which is very similar to the above construction, and we wish to recognize also that Borkar achieves what is in essence equivalent to Witsenhausen's static reduction reviewed Section \ref{EquivIS}. For simplicity, we only consider here the case where state and observation spaces are finite. We consider a discrete-time Markov decision process $\{x_n\}$ on a finite state space $\mathbb{X}$, controlled by a control process $\{u_n\}$ taking values in a compact Borel action space $\mathbb{U}$, and with an associated observation process $\{y_n\}$ taking values in a finite observation space $\mathbb{Y}$, where $n=0,1,2,\ldots$. The evolution of $\{x_n,y_n\}$ is given by
\begin{align}
P\big(x_{n+1}, y_{n+1} \in \cdot \big| x_m, x_m, u_m, m \leq n\big) = \rho(x_{n+1}, y_{n+1} \in \cdot | x_n,u_n), \nonumber
\end{align}
where $\rho: \mathbb{X} \times \mathbb{U} \rightarrow {\cal P}(\mathbb{X}) \times {\cal P}(\mathbb{Y})$ is some transition kernel. To ease the exposition, we assume that $\rho$ is of the following form:
\begin{align}
\rho(x_{n+1},y_{n+1} | x_n,u_n) = r(y_{n+1}|x_{n+1}) \otimes p(x_{n+1}|x_n,u_n), \nonumber
\end{align}
where $p$ is the state transition kernel and $r$ is the observation kernel. The initial distribution of $x_0$ is $\mu$.

A control process $\{u_n\}$ is admissible in classical sense if it is adapted to the filtration $\{\sigma(y_m,m\leq n)\}$ generated by observations $\{y_n\}$. In this case, one can write
\begin{align}
u_n = \pi_n(y_0,\ldots,y_n), n \geq 0,
\end{align}
for some $\pi_n: \prod_{k=0}^n \mathbb{Y} \rightarrow \mathbb{U}$. Let us denote $\pi = \{\pi_n\}$. 


Note that one can always write the evolution of the state process $\{x_n\}$ as a noise-driven dynamical system
\begin{align}
x_{n+1} = F(x_n,u_n,w_n), \label{eq1}
\end{align}
where $F:\mathbb{X} \times \mathbb{U} \times [0,1] \rightarrow \mathbb{X}$ is measurable and $\{w_n\}$ are independently and identically distributed uniformly on $[0,1]$. Using this dynamical system, we now reproduce the above process on a more convenient probability space. This will then enable us to define \textit{wide-sense admissible policies}.


\syr{In the following, we reduce the problem to an independent static one via Witsenhausen/Girsanov/Borkar, see Borkar's \cite{Bor00,Bor07} explicit analysis or Witsenhausen's method presented in Section \ref{EquivIS}.}

\syr{Under this reduction, we obtain a new probability space $P_0^{\pi}$ under which:}
\begin{itemize}
\item [(a)] $\{y_n\}$ is i.i.d. uniform on $\mathbb{Y}$ and independent of $x_0$ and $\{w_n\}$,
\item [(b)] $\{u_n,y_0,\ldots,y_n\}$ is independent of $\{w_n\}$, $x_0$, and $\{y_m, m>n\}$, for all $n$.
\end{itemize}
Using these properties, Borkar defined wide sense admissible policies in \cite{Bor00,Bor07} as follows. A policy $P_0$ is \textit{wide sense admissible} if $P_0$ satisfies (a) and (b). 
Note that condition (b) is very similar to non-signaling condition introduced in Section~\ref{nonsignal}. It assumes that action at time $n$ is independent of the observations after time $n$ given past observations (but {\bf not} necessarily past actions). In other words, instead of saying that action $u_n$ should be dependent on current and past observations $\{y_0,\ldots,y_n\}$, this condition states that action $u_n$ should be independent of future observations $\{y_m, m>n\}$ given past observations. \syr{Borkar and Budhiraja thankfully realized this seemingly simple, but consequential as we will see later, typo as this was corrected in further publications: in \cite{Bor03,BoBu04}, when Borkar and Budhiraja extend this definition to the continuous space case, they have slightly changed the condition: This condition is denoted by (b') and stated as follows:
\begin{itemize}
\item [(b')] $\{u_0,\ldots,u_n,y_0,\ldots,y_n\}$ is independent of $\{w_n\}$, $x_0$, and $\{y_m, m>n\}$, for all $n$.
\end{itemize}
In (b'), in addition to $\{u_n,y_0,\ldots,y_n\}$, we also suppose that past actions $u_0,\ldots,u_{n-1}$ are independent of $\{w_n\}$, $x_0$, and $\{y_m, m>n\}$. In other words, $u_n$ is independent of the observations after time $n$ given past observations and \textit{actions} (note that in condition (b) past actions are missing). This is indeed the right relaxation since in the next section, we establish via a counterexample that the optimal value achieved by wide sense admissible policies (with condition (b)) is strictly better than the optimal value achieved by classically admissible policies. Moreover, we show that an optimal wide-sense admissible policy under (b) evidently has access to control policy that it should not have; that is, it violates the causality of the problem, which is prohibited in general.}

\subsection{A counterexample}

%

We have partially observed MDP with the components $\mathbb{X} = \{0,1\} \times \{0,1\} \times \{0,1\}$, $\mathbb{U} = \{0,1\}$, and $\mathbb{Y} = \{y^*\}$ contains only one element. Let $x_0 \sim \pi_0 \otimes \gamma_0$, where $\pi_0$ is a distribution on $\{0,1\} \times \{0,1\}$ and $\gamma_0$ is a distribution on $\{0,1\}$. A typical element of $\mathbb{X}$ is denoted by $x = [x^1,x^2,x^3]$. The transition and the observation kernels are defined as follows:
\begin{align}
p\big([x_{n+1}^1,x_{n+1}^2,x_{n+1}^3] \big| [x_{n}^1,x_{n}^2,x_{n}^3],u_n \big) &= \lambda(x_{n+1}^1,x_{n+1}^2|x_{n}^1,x_{n}^2) \otimes \delta_{u_n}(x_{n+1}^3) \nonumber \\
r(y_n|x_n) &= \delta_{y^*}(y_n), \nonumber
\end{align}
where $\lambda : \{0,1\} \times \{0,1\} \rightarrow P(\{0,1\} \times \{0,1\})$ is a stochastic kernel such that $\pi_0$ is an invariant probability measure of $\lambda$. The reward\footnote{We note that all results in this paper apply with straightforward modifications for the case of maximizing reward instead of minimizing cost.} function is given by
\begin{align}
c(x,u) =
\begin{cases}
1, & \text{if $x^3 \oplus u = x^1 \cdot x^2$} \\
0, & \text{otherwise}.
\end{cases}\nonumber
\end{align}
To construct this counterexample, we have been inspired by the CHSH team \cite{CHSH1969} reviewed earlier in Theorem~\ref{relation-thm}, which establishes that non-signaling policies are not admissible relaxations as they strictly improve the performance. A similar conclusion will be obtained in this counterexample. 

In this problem, the observation process $\{y_n\}$ is non-informative as, for each $n$, $y_n \sim \delta_{y^*}$. Since any admissible control process $\{u_n\}$ is adapted to the filtering generated by observation process; that is, one can write
\begin{align}
u_n = \pi_n(y_0,\ldots,y_n), n \geq 0,
\end{align}
for some $\pi_n: \prod_{k=0}^n \mathbb{Y} \rightarrow \mathbb{U}$, any admissible control process $\{u_n\}$ can be represented as a deterministic sequence $\{a_n\}_{n=0}^{\infty} \subset \{0,1\}^{\infty}$, i.e., $u_n \sim \delta_{a_n}$ for all $n$. Note that, by the definition of the state transition kernel, the distribution $\pi_0$ of the first two components of the state $(x_n^1,x_n^2)$ remains as it is during the evolution of the state. With these observations, the maximum reward an admissible policy can attain is
\begin{align}
\max\big\{\pi_0(x_0^1 \cdot x_0^2 = 0),\pi_0(x_0^1 \cdot x_0^2 = 1)\big\}. \label{cls}
\end{align}
For instance, if
$$\max\big\{\pi_0(x_0^1 \cdot x_0^2 = 0),\pi_0(x_0^1 \cdot x_0^2 = 1)\big\} = \pi_0(x_0^1 \cdot x_0^2 = 0),$$
then this reward can be obtained by picking the control sequence as follows:
$$u_0 = 0, u_1 = 1, u_2 = 0, u_3 = 1, \ldots$$

We now construct an optimal wide sense admissible policy with the reward function $1$, which is in general strictly larger than (\ref{cls}). Let the probability distribution $P_0$ on $\Omega$ has the following properties:
\begin{itemize}
\item [(a)] $\{y_n\}$ is i.i.d. uniform on $\mathbb{Y}$ and independent of $x_0$ and $\{w_n\}$,
\item [(b)] For any $n \geq 0$,
\begin{align}
P_0(u_n,u_{n+1} | x_m,y_m,w_m, m \leq n) &=
\begin{cases}
1/2, &\text{if $u_n \oplus u_{n+1} = x_n^1 \cdot x_n^2$} \\
0, &\text{otherwise}
\end{cases}\nonumber
\end{align}
\end{itemize}
One can prove that $P_0$ is wide sense admissible; that is, $\{u_n,y_0,\ldots,y_n\}$ is independent of $\{w_n\}$, $x_0$, and $\{y_m, m>n\}$, for all $n$, under $P_0$. Indeed, $P_0(u_n | x_m,y_m,w_m, m \leq n) = {\mathop{\rm U}}_{\{0,1\}}(u_n)$ and $P_0(u_{n+1} | x_m,y_m,w_m, m \leq n) = {\mathop{\rm U}}_{\{0,1\}}(u_{n+1})$, where ${\mathop{\rm U}}_{\{0,1\}}$ is the uniform distribution on $\{0,1\}$. The average reward of $P_0$ is $1$, which can be the maximum achievable by any policy as $0 \leq c \leq 1$. Hence, $P_0$ is the optimal wide-sense policy.

Note that if
$$\max\{\pi_0(x_0^1 \cdot x_0^2 = 0),\pi_0(x_0^1 \cdot x_0^2 = 1)\} = \pi_0(x_0^1 \cdot x_0^2 = 0) = \pi_0(x_0^1 \cdot x_0^2 = 1) = 1/2,$$
then classical policies at most attain the reward of $1/2$, while wide sense admissible policies can get reward of $1$ which is twice as big as $1/2$. Hence, this proves the fact that wide sense admissible policies are much bigger than the classical policies. Furthermore, $P_0$ cannot be implemented in real life control applications as the action $u_n$ at any time $n$ depends on the knowledge of the action $u_{n+1}$ at time $n+1$, which violates the causality of the problem. Hence, this shows that for certain class of problems, wide-sense admissible policies are not a legitimate class of policies to study in discrete-time setup. In order to work with such policies, it is necessary to establish a theorem like \cite[Theorem 6.1]{FlPa82} that establishes equivalency of these two classes of policies in terms of achievable optimal value. This can indeed be done if condition (b) is replaced by condition (b'). Note that the wide sense policy $P_0$ defined above violates the condition (b'). Therefore, it is not wide sense admissible under condition (b'). 

In summary, if the relaxation is done with the interpretation that $u_n$ is conditionally independent from future observations given the past observations, an incorrect conclusion can be made. If the relaxations is so that the control action $u_n$ is conditionally independent from future observations, given all the past, then relaxation is valid.

\section{Conclusion, Topics Left Out and Some Open Problems}
\subsection{Concluding Remarks}
The way information is decentralized is a key attribute determining how to approach a problem in various areas of applied mathematics. In this review article, we studied information structures in a probability theoretic and topological context. We defined information structures, placed various topologies on them, and study closedness and compactness properties. We presented existence and approximation results for optimal decision/control policies. We discussed various upper and lower bounding techniques, through relaxations and convex programs ranging from classically realizable and classically non-realizable (such as quantum and non-signaling) relaxations. Figure \ref{HierInfoStrategicMeasure} depicts a summary of our findings on strategic measures.

We later presented various topologies on decision/control strategies defined independent of information structures, but for which information structures determine whether the topologies have utility in arriving at existence, compactness, convexification or approximation results. 

We showed that viewing decentralization with regard to the induced strategic measures and viewing decentralization with respect to the applied control policies, lead to different operational conclusions and tools to arrive at optimality conditions and results for optimal decentralized decision making.

We showed that externally provided randomness should be well motivated and when one defines relaxed control; which is a common solution technique in discrete-time and continuous-time classical stochastic control, this should be introduced cautiously since such randomness may indeed lead to unacceptable performance improvement. 

%

In the paper we considered the set of correlations {\it given} an information structure. A further related set of results involve the problem of optimal {\it design} of information structures. This subject is beyond the scope of this review article with some related results in \cite[Thms. 3.2, 3.3, 3.4]{YukselOptimizationofChannels} and \cite{YukselBasarBook}, which also directly apply to the multi-agent setting.

%
%
%
%

\subsection{Discussion on Existence and Some Related Results}

As noted earlier, while the topological constructions on policies are quite useful, Theorem \ref{existenceRelaxed3} (for static teams or dynamic teams with an independent-measurements reduction) and Theorems \ref{SuffCon2} and \ref{SuffCon2'''} (for sequential teams that do not allow an independent-measurements reduction) are the most general existence results, to our knowledge, for the problems considered here. However, some slightly weaker versions of these results can be arrived at through different methods, as laid out in the paper. 
We also note, for completeness, that existence of optimal policies for static and a class of sequential dynamic teams had been studied in \cite{WuVer11}, \cite{wit68} \cite{gupta2014existence,YukselSaldiSICON17}, \cite{YukselOptimizationofChannels} \cite{YukselWitsenStandardArXiv}. Conditions for optimality have been established in \cite{charalambous2017centralizedI} for a class of continuous-time decentralized stochastic control problems. We also noted in the article, for (classical) partially observed models, existence results in \cite{FlPa82}, \cite{bismut1982partially} with a discrete-time review in \cite{YukselWitsenStandardArXiv}. For a class of teams which are convex, one can reduce the search space to a smaller parametric class of policies, such as linear policies for quasi-classical linear quadratic Gaussian problems \cite{rad62,KraMar82,HoChu}.

Note that standard dynamic programming is known to be a useful tool for a class of dynamic teams since 1970s. Clearly, if all the information at any given decision maker is common knowledge between all decision makers, then the system can be viewed to be a centralized system and standard dynamic programming is applicable. However, if only some of the system variables are common knowledge, the remaining unknowns may or may not lead to a computationally tractable program generating an optimal solution. A possible approach toward establishing a tractable program is through the construction of a controlled Markov chain where the controlled Markov state may now live in a larger state space (for example a space of probability measures) and the actions are elements in possibly function spaces. This controlled Markov construction may lead to a computation of optimal policies. Such a {\it dynamic programming approach} has been adopted extensively in the literature (see for example, \cite{Athans}, \cite{yos75}, \cite{ChongAthans}, \cite{AicardiDavoli}, \cite{YukTAC09}, \cite{lamperski2015optimal} and significantly generalized and termed as the {\it common information approach} in \cite{NayyarMahajanTeneketzis} and \cite{NayyarBookChapter}) through the use of a team-policy which uses common information to generate partial functions for each DM to generate their actions using local information. This construction requires a common knowledge among decision makers, which is a unifying assumption in the aforementioned contributions in the literature. Witsenhausen \cite{WitsenStandard} and \cite{YukselWitsenStandardArXiv} developed universal dynamic programming algorithms, which are conceptually useful and mathematically consequential (on existence and recursive analysis) but practically of limited algorithmic use with our current knowledge.

\subsection{Some Open Problems}


\noindent In the following, we present a number of open problems:

\begin{itemize}
\item[(a)] Note that one can prove, using the same argument in Theorem~\ref{relation-thm}-(ii), that 
$$\inf_{P \in L_{Q(1)}(\mu)} \int P(ds) \, c(s)\leq J^*,$$
since there exists an optimal individually randomized policy that achieves $J^*$. Therefore, for any $d$, the team problem $\inf_{P \in L_{Q(d)}(\mu)} \int P(ds) \, c(s)$ is indeed an \textit{admissible} extension or relaxation of classical team problem since it does not require any communication between agents and a mediator. Hence, the solution of this problem will be a significant contribution to the team decision theory as the optimal quantum-correlated policy can be realizable in real life in view of recent and potential advances in quantum technology. However, it is important to note that $L_{Q(d)}(\mu)$ cannot be convex if the dimension constraint $d$ is small \cite{MaWo15} (see also Theorem~\ref{quantum-convex-closed}). Therefore, with dimensionality constraint, the optimization problem corresponding to the quantum-correlated strategic measures can be non-convex. To convexify the problem, we can either add unlimited common randomness without changing $d$ or increase $d$ and allow for a limited common randomness \cite{MaWo15}. We can now state the following open problem: 

\begin{itemize}
\item[\textbf{(OP1)}] Under what conditions on the components of the team, there exists $d \geq 1$ such that the optimization problem
\begin{align}
\inf_{P \in L_{Q(d)}(\mu)} \int P(ds) \, c(s) \label{optQd}
\end{align}
can be written or can be approximated by a semi-definite program? 
\end{itemize}
One way to solve \textbf{(OP1)} might be to use the \textit{so-called} NPA hierarchy \cite{NaPiAc08}, which provides an infinite hierarchy of SDP outer approximations to the set of quantum-correlated strategic measures. However, although it gives SDP outer approximations, there are no bounds on the rate of convergence quantifying how the approximation improves as the level in the hierarchy increases. On the other hand, if we can bound the rate of convergence for the NPA hierarchy, it will be possible to use it for solving \textbf{(OP1)}.

\item[(b)]
Recall the dual program introduced in Section~\ref{dual} using non-signaling policies. In this section, we establish dual problem for the celebrated Witsenhausen's counterexample and pose several open problems. To this end, we let $\mathbb{Y}^1 = \mathbb{Y}^2 = \mathbb{U}^1 = \mathbb{U}^2 = \mathbb{R}$, which denote observation and action spaces of agents. In Witsenhausen's celebrated counterexample \cite{WitsenhausenCounter}, depicted in Fig.~\ref{fig1}, there are two decision makers: Agent~$1$ observes a zero mean and $\sigma$-variance Gaussian random variable $y^1 \in \mathbb{Y}^1$ and decides its action $u^1 \in \mathbb{U}^1$. Agent~$2$ observes $y^2 \coloneqq u^1 + v \in \mathbb{Y}^2$, where $v \in \mathbb{V} = \mathbb{R}$ is zero mean and unit variance Gaussian noise independent of $y^1$, and decides its action $u^2 \in \mathbb{U}^2$

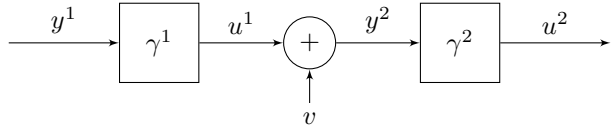
\begin{figure}[h]
\centering
\tikzstyle{int}=[draw, fill=white!20, minimum size=3em]
\tikzstyle{init} = [pin edge={to-,thin,black}]
\tikzstyle{sum} = [draw, circle]
\begin{tikzpicture}[node distance=2cm,auto,>=latex']
    \node [int] (a) {$\gamma^1$};
    \node (b) [left of=a,node distance=2cm, coordinate] {a};
    \node [sum] (d) [right of=a] {$+$};
    \node [int] (c) [right of=d] {$\gamma^2$};
    \node [coordinate] (end) [right of=c, node distance=2cm]{};
    \node (e) [below of=d, node distance=1cm] {$v$};
    \path[->] (e) edge node {} (d);
    \path[->] (b) edge node {$y^1$} (a);
    \path[->] (a) edge node {$u^1$} (d);
    \path[->] (d) edge node {$y^2$} (c);
    \path[->] (c) edge node {$u^2$} (end) ;
\end{tikzpicture}
\caption{Witsenhausen's counterexample.}
\label{fig1}
\end{figure}

The cost function of the team is given by
\begin{align}
c(y^1,u^1,u^2) = k^2 (u^1 - y^1)^2 + (u^2 - u^1)^2, \nonumber
\end{align}
where $k > 0$. Let $g(y) \coloneqq \frac{1}{\sqrt{2\pi}}\exp{\{-y^2/2\}}$ be probability density function of zero mean and unit variance Gaussian random variable. Then we have
\begin{align}
P(y^2 \in S\,|\,u^1) = \int_{S} g(y^2-u^1) \, m(dy^2), \nonumber
\end{align}
where $m$ denotes the Lebesgue measure on $\mathbb{Y}^2$. Let
\begin{align}
f(u^1,y^2) \coloneqq \exp{\biggl\{-\frac{(u^1)^2-2y^2u^1}{2}\biggr\}} \label{eq20}
\end{align}
so that $g(y^2-u^1) =f(u^1,y^2)\frac{1}{\sqrt{2\pi}}\exp{\{-(y^2)^2/2\}} = f(u^1,y^2) \, g(y^2)$. The independent static reduction of Witsenhausen counterexample proceeds as follows: for any strategic measure \syr{$P \in L_R(\mu)$}, the expected cost can be written as 
\begin{align}
J(P) &= \int c(y^1,u^1,u^2) \, P(du^2|y^2)\, P(dy^2|u^1) \, P(du^1|y^1) \, \mu_{\sigma}(dy^1) \nonumber \\
&= \int c(y^1,u^1,u^2) \, f(u^1,y^2) \, \mu_{1}(dy^2) \mu_{\sigma}(dy^1), \nonumber
\end{align}
where $\mu_{\rho}$ denotes zero mean and $\rho$-variance Gaussian distribution. Hence, by defining 
$
c_s(y^1,y^2,u^1,u^2) = c(y^1,u^1,u^2) \, f(u^1,y^2) $
and 
$\mu(dy^1,dy^2) = \mu_{\sigma}(dy^1) \, \mu_{1}(dy^2),$
we can write $J(P)$ as
\begin{align}
J(P) = \int c_s(y^1,y^2,u^1,u^2) \, P(du^2|y^2)\, P(du^1|y^1) \, \mu(dy^1,dy^2). \label{eq4}
\end{align}
Therefore, in the static reduction of Witsenhausen's counterexample, the agents observe independent zero mean Gaussian random variables. In the remainder of this note, we consider static reduction of Witsenhausen's counterexample. Note that a strategic measure $P \in {\cal P}(\mathbb{Y}^1\times\mathbb{Y}^2\times\mathbb{U}^1\times\mathbb{U}^2)$ is \textit{non-signaling} if
\begin{align}
P(du^1 | y^1,y^2) &= P(du^1 | y^1), \nonumber \\
P(du^2 | y^1,y^2) &= P(du^2 | y^2). \label{nnnonsignaling}
\end{align}
As noted in Section~\ref{dual}, since the constraints in (\ref{nnnonsignaling}) for $P$ are linear, the optimal team cost with non-signaling policies can be written as a linear program over an appropriate vector spaces as follows.
Recall that, for any metric space $\mathbb{E}$, ${\cal M}(\mathbb{E})$ denotes the set of finite signed measures on $\mathbb{E}$ and $C(\mathbb{E})$ denotes the set of continuous real functions. Consider the vector spaces ${\cal M}(\mathbb{U}^1 \times \mathbb{U}^2 \times \mathbb{Y}^1 \times \mathbb{Y}^2)$, $C(\mathbb{U}^1 \times \mathbb{U}^2 \times \mathbb{Y}^1 \times \mathbb{Y}^2)$, ${\cal M}(\mathbb{U}^1 \times \mathbb{Y}^1 \times \mathbb{Y}^2)$, $C(\mathbb{U}^1 \times \mathbb{Y}^1 \times \mathbb{Y}^2)$, ${\cal M}(\mathbb{U}^2 \times \mathbb{Y}^1 \times \mathbb{Y}^2)$, and $C(\mathbb{U}^2 \times \mathbb{Y}^1 \times \mathbb{Y}^2)$. Let us define bilinear forms on
$$\bigl({\cal M}(\mathbb{U}^1 \times \mathbb{U}^2 \times \mathbb{Y}^1 \times \mathbb{Y}^2),C(\mathbb{U}^1 \times \mathbb{U}^2 \times \mathbb{Y}^1 \times \mathbb{Y}^2)\bigr)$$ and on $$\bigl({\cal M}(\mathbb{U}^1 \times \mathbb{Y}^1 \times \mathbb{Y}^2) \times {\cal M}(\mathbb{U}^2 \times \mathbb{Y}^1 \times \mathbb{Y}^2) \times \mathbb{R},C(\mathbb{U}^1 \times \mathbb{Y}^1 \times \mathbb{Y}^2) \times C(\mathbb{U}^2 \times \mathbb{Y}^1 \times \mathbb{Y}^2) \times \mathbb{R}\bigr)$$ by letting
\begin{align}
&\langle \rho,v  \rangle_{1} \coloneqq \int_{\mathbb{U}^1 \times \mathbb{U}^2 \times \mathbb{Y}^1 \times \mathbb{Y}^2} v(u^1,u^2,y^1,y^2) \, \rho(du^1,du^2,dy^1,dy^2) \label{eqqq2}, \\
&\langle (\rho_1,\rho_2,a) ,(v_1,v_2,b)  \rangle_{2}  
\coloneqq \int_{\mathbb{U}^1 \times \mathbb{Y}^1 \times \mathbb{Y}^2} v_1(u^1,y^1,y^2) \, \rho_1(du^1,dy^1,dy^2) \nonumber \\
&\phantom{xxxxxxxxxxxx} + \int_{\mathbb{U}^2 \times \mathbb{Y}^1 \times \mathbb{Y}^2} v_2(u^2,y^1,y^2) \, \rho_2(du^2,dy^1,dy^2) + a b \label{eqqq1}.
\end{align}
The bilinear forms in (\ref{eqqq2}) and (\ref{eqqq1}) constitute duality between spaces \cite[Chapter IV.3]{Bar02}. Hence, the topologies on these spaces should be understood as the weak topology of the duality induced by these bilinear forms. We define the linear map $L: {\cal M}(\mathbb{U}^1 \times \mathbb{U}^2 \times \mathbb{Y}^1 \times \mathbb{Y}^2) \rightarrow {\cal M}(\mathbb{U}^1 \times \mathbb{Y}^1 \times \mathbb{Y}^2) \times {\cal M}(\mathbb{U}^2 \times \mathbb{Y}^1 \times \mathbb{Y}^2) \times \mathbb{R}$ by
\begin{align}
L(\rho) = \left(L_1(\rho), L_2(\rho), L_3(\rho)\right), \nonumber
\end{align}
where
\begin{align}
L_1: \rho(du^1,du^2,dy^1,dy^2) &\mapsto \rho(du^1,dy^1,dy^2) - \rho(du^1,dy^1) \, \mu_1(dy^2) \nonumber \\
L_2: \rho(du^1,du^2,dy^1,dy^2) &\mapsto \rho(du^2,dy^1,dy^2) - \rho(du^2,dy^2) \, \mu_{\sigma}(dy^1) \nonumber \\
L_3: \rho(du^1,du^2,dy^1,dy^2) &\mapsto \langle \rho,1  \rangle_{1} \nonumber.
\end{align}
Using $L$, the optimal value of the team with non-signaling policies can be written as a linear program as follows:
\begin{align}
(\textbf{NS}) \text{                         }&\text{minimize}_{\rho \in {\cal M}_+(\mathbb{U}^1 \times \mathbb{U}^2 \times \mathbb{Y}^1 \times \mathbb{Y}^2)} \text{ } \langle \rho,c_s \rangle_1
\nonumber \\*
&\text{subject to  } L(\rho) = (0,0,1). \label{aaaaa}
\end{align}
Since $L_C(\mu) \subset L_{NS}(\mu)$, the solution of above linear program gives a lower bound to original formulation of Witsenhausen's counterexample. Note that the dual of $L$ is given by $L^*: C(\mathbb{U}^1 \times \mathbb{Y}^1 \times \mathbb{Y}^2) \times C(\mathbb{U}^2 \times \mathbb{Y}^1 \times \mathbb{Y}^2) \times \mathbb{R} \rightarrow C(\mathbb{U}^1 \times \mathbb{U}^2 \times \mathbb{Y}^1 \times \mathbb{Y}^2)$, where 
\begin{align}
L^*(v_1,v_2,b) &= v_1(u^1,y^1,y^2) - \int_{\mathbb{Y}^2} v_1(u^1,y^1,y^2) \, d\mu_1(y^2) \nonumber \\ 
&+ v_2(u^2,y^1,y^2) - \int_{\mathbb{Y}^1} v_2(u^2,y^1,y^2) \, d\mu_{\sigma}(y^1) + b .
\end{align}
Then the dual program of $(\textbf{NS})$ can be written as \cite[Chapter IV.6]{Bar02}
\begin{align}
(\textbf{NS}^*) \text{                         }&\text{maximize}_{(v_1,v_2,b) \in C(\mathbb{U}^1 \times \mathbb{Y}^1 \times \mathbb{Y}^2) \times C(\mathbb{U}^2 \times \mathbb{Y}^1 \times \mathbb{Y}^2) \times \mathbb{R}} \text{ } b
\nonumber \\*
&\text{subject to  } L^*(v_1,v_2,b) \leq c_s. \label{aaaab}
\end{align}
We can now state the following problems:
\begin{itemize}
\item [\textbf{(OP2)}] Approximate numerically the dual linear program $(\textbf{NS}^*)$ to obtain a lower bound to $J^*$. 
\syr{\item [\textbf{(OP3)}] Analyze the difference 
$$\inf_{P \in L_{NS}(\mu)} \int P(ds) \, c_s(s) - \inf_{P \in L_C(\mu)} \int P(ds) \, c_s(s),$$ 
or at least with an upper bound obtained with numerical methods (see \cite{saldiyuksellinder2017finiteTeam} for a review of numerical results in the literature).
\item [\textbf{(OP4)}] We know that when there is no $u^1 \cdot u^2$ term in the original cost function $c$, there is an affine optimal policy which can be obtained analytically using information theoretic tools \cite{BaBa87}. In this case, what is the relation between
\[\inf_{P \in L_{NS}(\mu)} \int P(ds) \, c_s(s)\]
and
\[\inf_{P \in L_C(\mu)} \int P(ds) \, c_s(s)?\]
More generally, for which cost functions, can we establish that
 $$\inf_{P \in L_{NS}(\mu)} \int P(ds) \, c_s(s) = \inf_{P \in L_C(\mu)} \int P(ds) \, c_s(s)?$$
}
\end{itemize}
\item[(c)] In the paper, we considered only sequential decentralized stochastic control. As noted earlier, if there is a pre-defined order in which the decision makers act, then we say that a system is {\it sequential}; otherwise, the system is non-sequential. Such non-sequential systems are substantially more difficult to study, since the ambiguities in the order of actions lead to challenges on the interpretation of local information. Optimal design of such non-sequential models requires the systems to be deadlock-free, that is the actions of a given DM should not depend on the actions of DMs acting in the future, for any realized random ordering. Furthermore, the optimization problem for such systems should be well posed/solvable, since for some designs the expected cost may not be well-defined. We refer the reader to Witsenhausen \cite{WitsenhausenSIAM71}, Andersland and Teneketzis \cite{AnderslandTeneketzisI}, \cite{AnderslandTeneketzisII} and Teneketzis \cite{Teneketzis2} for a comprehensive study of non-sequential systems, see \cite{YukselBasarBook} for a brief review. We also note that with the observation that the information fields generated by local measurements lead to subtle conditions on solvability and causality, an alternative probabilistic model, based on quantum mechanics, for describing such problems has been proposed by Baras in \cite{Baras1} and \cite{Baras2}. These papers also present an accessible review of related developments on quantum information prior to the publications. In summary, the study of non-sequential decentralized stochastic control systems in the context of what we studied throughout this paper is an open problem. 
\end{itemize}

\section{Acknowledgements}
The authors would like to acknowledge discussions and collaborations with Tamer Ba\c{sar}, Tam\'as Linder, Sina Sanjari, Ian Hogeboom-Burr, Ankur Kulkarni, Abhishek Gupta, Aditya Mahajan, Nuno Martins, and Demos Teneketzis. 



\end{document}